\documentclass[11pt]{amsart}   	
\usepackage{geometry}             
\usepackage[dvipsnames]{xcolor}   		             		
\usepackage{graphicx}				
\usepackage{mathrsfs}
\usepackage{amssymb,amsthm,amsmath,mathtools}
\usepackage{tikz}
\usepackage{tikz-cd}
\usepackage{accents,upgreek,enumerate}
\usepackage[headings]{fullpage}
\usepackage{bm}
\usepackage[all]{xy}
\usepackage{caption}
\usepackage{floatrow}       
\usepackage{csquotes}
\usepackage{enumitem}
\usepackage{comment}
\usepackage{rotating}       
\usepackage{makecell}       
\usepackage{multirow}       
\usepackage{crossreftools}  
\usepackage{mathdots}
\usepackage{hyphenat}

\usetikzlibrary{calc}
\usetikzlibrary{fadings}
\usetikzlibrary{decorations.pathmorphing}
\usetikzlibrary{decorations.pathreplacing}

\DeclareMathAlphabet{\mathpzc}{OT1}{pzc}{m}{it}

\usepackage[backref]{hyperref}
\hypersetup{
colorlinks   = true,          
urlcolor     = violet,          
linkcolor    = Purple,          
citecolor   = RubineRed             
}
\usepackage[capitalize,noabbrev]{cleveref}

\newtheorem{theorem}{Theorem}[section]
\newtheorem{corollary}[theorem]{Corollary}
\newtheorem{lemma}[theorem]{Lemma}
\newtheorem{proposition}[theorem]{Proposition}

\newtheorem{quasi-theorem}[theorem]{Quasi-Theorem}

\theoremstyle{definition}
\newtheorem{definition}[theorem]{Definition}

\newtheorem{example}[theorem]{Example}
\newtheorem{remark}[theorem]{Remark}

\newtheorem{notation}[theorem]{Notation}
\newtheorem{setting}[theorem]{Setting}

\newcommand{\direction}{\vec{\mathbf{v}}} 

\newcommand{\quotient}[2]{\raisebox{0.5ex}{$#1$}/\raisebox{-0.5ex}{$#2$}}

\newcommand{\bigmid}{\mathrel{\big|}}

\newcommand{\A}{{\mathbb{A}}}                     
            
\newcommand{\NN} {{\mathbb N}}		
\newcommand{\PP}{\mathbb{P}}

\newcommand{\RR} {{\mathbb R}}	
\newcommand{\CC} {{\mathbb C}}
\newcommand{\ZZ} {{\mathbb Z}}	
\renewcommand{\AA}{\mathbb{A}}	

\def\L{\mathrm{L}}

\def\setminus{\smallsetminus}

\newcommand{\Puiseux}[2][t]{#2\{\!\{#1\}\!\}}

\newcommand{\conj}[1]{\overline{#1}}
\newcommand{\conjb}[1]{\widehat{#1}}

\DeclareMathOperator{\Tr}{Tr}               
\DeclareMathOperator{\tr}{tr}               

\DeclareMathOperator{\chara}{char} 
 
\DeclareMathOperator{\conv}{conv}

\DeclareMathOperator{\val}{val}
\DeclareMathOperator{\ev}{ev}
\DeclareMathOperator{\floor}{floor}

\newcommand{\gw}[1]{\left\langle #1 \right\rangle}
\newcommand{\qinv}[1]{\left\langle#1\right\rangle}
\newcommand{\qqinv}[1]{\left(\left\langle#1\right\rangle\right)}

\newcommand{\h}{h}                              

\DeclareMathOperator{\GW}{GW}                   
\DeclareMathOperator{\Wel}{\operatorname{Wel}^{\mathbb{A}^1}}
\DeclareMathOperator{\Hess}{Hess}						
\DeclareMathOperator{\mass}{mass}               

\newcommand{\rk}{\operatorname{rk}}
\newcommand{\sgn}{\operatorname{sgn}}

\newcommand{\W}{\operatorname{W}}

\newcommand{\cal}{\mathcal}

\def\cD{{\cal D}}

\def\cT{{\cal T}}

\def\trop{\mathrm{trop}}
\def\floor{\mathrm{floor}}
\def\mult{\mathrm{mult}}
\newcommand{\elev}{\mathrm{elev}}

\newcommand{\Spec}{\operatorname{Spec}}

\makeatletter
\def\blfootnote{\xdef\@thefnmark{}\@footnotetext}
\makeatother

\definecolor{darkspringgreen}{rgb}{0.09, 0.6, 0.1}
\definecolor{dsg}{rgb}{0.09, 0.6, 0.1}

\def\vertexsize {1.5pt}       
\def\edgewidth {0.5 pt}      
\def\halfedgelength {0.7}    
\def\doubleedgesep {0.03}      

\newcommand{\thinpoint}[1]{\fill (#1) circle [radius = \vertexsize];}

\newcommand{\fatpoint}[1]{\fill (#1) circle [radius = 2*\vertexsize];}

\newcommand{\edgeweight}[1]{{\color{blue}#1}}
\newcommand{\coordinates}[1]{{\color{darkspringgreen}#1}}
\newcommand{\quadextension}[1]{{\color{red}#1}}

\newcommand{\fatPointOnVertex}{{\ref{fatPointOnVertex}}}
\newcommand{\fatPointOnParallelogram}{{\ref{fatPointOnParallelogram}}}
\newcommand{\fourValentVertex}{{\ref{fourValentVertex}}}
\newcommand{\parallelogramWithOneDoubleEdge}{{\ref{parallelogramWithOneDoubleEdge}}}
\newcommand{\parallelogramWithTwoDoubleEdges}{{\ref{parallelogramWithTwoDoubleEdges}}}
\newcommand{\triangleWithMergedEdge}{{\ref{triangleWithMergedEdge}}}
\newcommand{\fourValentVertexWithMergedEdge}{{\ref{fourValentVertexWithMergedEdge}}}
\newcommand{\triangleWithTwoMergedEdges}{{\ref{triangleWithTwoMergedEdges}}}
\newcommand{\allDoubleVertex}{{\ref{allDoubleVertex}}}
\newcommand{\mergedTriangles}{{\ref{mergedTriangles}}}

\title{\bf Quadratically enriched plane curve counting via tropical geometry} 

\date{}

\author {Andr\'es Jaramillo Puentes}
\address {Dipartimento di Matematica e Informatica,
Universit\`a degli Studi di Catania,
Viale Andrea Doria~6,
95125 Catania, Italy}
\email {ga.jaramillopuentes@unict.it}
\author {Hannah Markwig}
\address {Universit\"at T\"ubingen, Fachbereich Mathematik, Auf der Morgenstelle 10, 72076 T\"ubingen, Germany }
\email {hannah@math.uni-tuebingen.de}
\author {Sabrina Pauli}
\address {TU Darmstadt, Fachbereich Mathematik, Schlossgartenstra\ss{}e 7, 64289 Darmstadt, Germany}
\email {pauli@mathematik.tu-darmstadt.de}
\author{Felix R\"ohrle}
\address {Universit\"at T\"ubingen, Fachbereich Mathematik, Auf der Morgenstelle 10, 72076 T\"ubingen, Germany }
\email{roehrle@math.uni-tuebingen.de}


\begin{document}

\begin{abstract}
We prove that the quadratically enriched count of rational curves in a smooth toric del Pezzo surface passing through $k$-rational points and pairs of conjugate points in quadratic field extensions $k\subset k(\sqrt{d_i})$ can be determined by counting certain tropical stable maps through vertically stretched point conditions with a suitable multiplicity. Building on the floor diagram technique in tropical geometry, we provide an algorithm to compute these numbers.

Our tropical algorithm computes not only these new quadratically enriched enumerative invariants, but simultaneously also the complex Gromov-Witten invariant, the real Welschinger invariant counting curves satisfying real point conditions only, the real Welschinger invariant of curves satisfying pairs of complex conjugate and real point conditions, and the quadratically enriched count of curves satisfying $k$-rational point conditions. 

\end{abstract}

\maketitle

\vspace{-0.2in}

\setcounter{tocdepth}{1}
\tableofcontents

\section{Introduction}

Enumerative geometry poses fascinating research problems since the days of the Ancient \linebreak 
Greeks. Among the interesting enumerative questions which continue to pose problems in today's research is the problem of counting rational plane curves of degree $d$ that satisfy the correct number of point incidence conditions. If we consider this problem over an algebraically closed field, we denote the number of such curves by $N_d$ and this is an \emph{invariant}, i.e.\ $N_d$ does not depend on the configuration of the chosen points conditions, as long as it is generic. 

If we work over the field of real numbers $\RR$ instead, the count of real rational plane curves of degree $d$ satisfying real point conditions, is \emph{not} an invariant \cite{DK00}. However, Welschinger suggested in a landmark paper \cite{Wel05} that we can weigh real curves with a sign inherent to the geometry of the curve and obtain the so-called Welschinger invariant.
This idea still works for the generalized version of the real counting problem where some of the real point conditions are replaced with pairs of complex conjugate point conditions. We denote the Welschinger invariant counting degree $d$ rational plane real curves through $r=3d-1-2s$ real points and $s$ pairs of complex conjugate points by $W_{d,s}$.

More recently, the count of rational plane curves is being studied over even more general base fields and this uses tools from $\mathbb{A}^1$-homotopy theory. Just as signed counts of real objects have generalized the idea of \enquote{counting} for the sake of obtaining an invariant, these new counts now assign quadratic forms as multiplicities. The resulting enumerative count lives in the Grothendieck-Witt ring $\GW(k)$ over the ground field $k$. 
In this spirit, Levine \cite{LevineWelschinger} and Kass-Levine-Solomon-Wickelgren \cite{KassLevineSolomonWickelgren} define a quadratically enriched count $N^{\AA^1}_d(\sigma)$ of rational plane curves. Here, $\sigma=(L_1,\ldots,L_m)$ is a list of finitely many finite field extensions of $k$ and $N^{\AA^1}_d(\sigma)$ counts rational plane curves of degree $d$ passing through a configuration of $m$ points with residue fields listed in $\sigma$ counted with a weight in $\GW(k)$.
Setting $k = \RR$ and $\sigma=(\RR^r,\CC^s)$, taking the signature of $N^{\AA^1}_d(\sigma)$, one gets the Welschinger invariant $W_{d,s}$, and taking the rank one gets the number $N_d$.

\subsection{Results}
The work of \cite{LevineWelschinger} and \cite{KassLevineSolomonWickelgren} defines the quadratically enriched curve counts with point conditions in a more general setting, namely for rational curves on a del Pezzo surface $S$ of a fixed a degree $\beta$ with some additional assumptions. We denote this invariant by $N^{\AA^1}_{S,\beta}(\sigma)$.
In this paper, we restrict ourselves to the case that $S$ is a toric del Pezzo surface, which already satisfies the assumptions from \cite{KassLevineSolomonWickelgren}.

\begin{setting} \label{setting}    
    Let $\Delta$ be the Newton polygon of a smooth toric del Pezzo surface $S$, i.e. one of the polygons shown in \cref{fig-degrees}.
    Let $k$ be a perfect field, either of characteristic 0, or of characteristic greater than the diameter of $\Delta$ and greater than 3. 
    Define $n (\Delta)\coloneqq \#(\partial\Delta \cap \mathbb{Z}^2)-1$ and let $r, s$ be natural numbers such that~$n(\Delta)=r + 2s$.
     Let $Q_i$ be a point of $S$ defined over a degree~$2$ field extension $k(\sqrt{d_i}) \supset k$, where $d_i \in k$ is an element which has no square root in $k$, for each $i = 1, \ldots, s$, and let~$P_j$ be a $k$-rational point in $S$,  for~$j = 1, \ldots, r$.  
\end{setting}

In \cref{setting} we write $N^{\mathbb{A}^1}_\Delta(r,(d_1,\ldots,d_s))$ for the quadratically enriched count $N^{\AA^1}_{S,\beta}(\sigma)$ of rational curves of class $\Delta$ in $S$ passing through the $Q_i$ and $P_j$ for $i=1,\ldots,s$ and $j = 1, \ldots, r$, that is for~$\sigma=\big(k^r,k(\sqrt{d_1}),\ldots,k(\sqrt{d_s})\big)$. Such $\sigma$ consisting only of trivial and quadratic field extensions is called \emph{multiquadratic}. 

After reviewing some preliminaries on tropical and quadratically enriched curve counting in \cref{sec-preliminaries}, we set up a tropical analogue of this counting problem in \cref{sec-tropdoublepoints}. 
On the tropical side, we count tropical stable maps passing through \emph{simple} and \emph{double} point conditions in $\RR^2$. 
The double points locally impose two conditions, and one could locally produce tropical curves satisfying two simple point conditions by a small perturbation. 
Since we know that the algebraic count $N^{\mathbb{A}^1}_\Delta(r,(d_1,\ldots,d_s))$ is an invariant, we restrict our study of the tropical enumerative problem to vertically stretched point conditions as this simplifies the computations in the entire article greatly. Under this vertically stretched assumption a complete classification of the combinatorial features of tropical curves which appear in this counting problem is achieved in \cref{lem-casesfatpoint}. Based on this classification, we define a quadratically enriched multiplicity $\mult^{\AA^1}(\Gamma)$ for any tropical curve $\Gamma$ which may appear in this context in \cref{section:computations}. The definition of the multiplicity may appear very involved, but it still takes the shape of a product over local contributions, which is very familiar to tropical geometers.
The main result of this paper is the following correspondence theorem, which identifies the quadratically enriched algebraic count with a tropical count.

\begin{theorem}[Main Theorem]\label{thm-main}
    In \cref{setting}, the quadratically enriched count of rational curves in the toric surface $S$ of class $\Delta$ passing through the points $Q_i$ and $P_j$ is equal to its tropical counterpart, summing the quadratically enriched multiplicity over tropical curves satisfying vertically stretched point conditions:
    \[N_\Delta^{\AA^1}(r,(d_1,\ldots,d_s))=N_\Delta^{\AA^1,\trop}(r,(d_1,\ldots,d_s))\coloneqq \sum_\Gamma \operatorname{mult}^{\A^1}(\Gamma).\]
\end{theorem}

Similar to $N_d^{\AA^1}(r,(d_1,\ldots,d_s))$ generalizing $N_d$ and $W_{s, d}$, we show in \cref{sec-compatible} that our tropical count generalizes various existing counts in tropical enumerative geometry (see \cref{sec-history} below for an overview of these counts).

The main theorem above also generalizes already existing tropical correspondence theorems for vertically stretched points (see section \ref{sec-history}) by Mikhalkin \cite{Mi03} and Shustin \cite{Shu06b}: $$N_\Delta=N_\Delta^{\trop}\; \text{ and }\;W_{\Delta,s}=W_{\Delta,s}^\trop$$
(where $N_\Delta$  and $W_{\Delta,s}$ denote the generalizations of $N_d$ and $W_{d,s}$ for the toric surface $\Delta$); 
and by the first and third author \cite{JPP23}:
$$N^{\AA^1}_\Delta=N^{\AA^1,\trop}_\Delta,$$
where $N_\Delta^{\AA^1}=N_\Delta^{\AA^1}(\sigma)$ and $N_\Delta^{\AA^1,\trop}=N_\Delta^{\AA^1,\trop}(\sigma)$ with $\sigma=(k^n)=(k,\ldots,k)$.

\begin{theorem} \label{thm-generalizing}
    The tropical invariant $N_\Delta^{\AA^1, \trop}(r, (d_1, \ldots, d_s))$ has the following properties.
    \begin{enumerate}
        \item (\cref{prop-specializetokrational_s_0}) If $s=0$, then our count equals the quadratically enriched count from \cite{JPP23}, 
        $$N_\Delta^{\AA^1,\trop}(r,(d_1,\ldots,d_s))=N_\Delta^{\mathbb{A}^1,\trop}.$$
        
        \item (\cref{prop-specializetokrational}) If $d_s$ is a square in $k$, then
        $$N_\Delta^{\AA^1,\trop}(r,(d_1,\ldots,d_s))=N_\Delta^{\mathbb{A}^1,\trop}(r+2,(d_1,\ldots,d_{s-1})).$$
        In particular, if all the $d_i$'s are squares in $k$, out count equals the quadratically enriched count $N^{\AA^1,\trop}_\Delta$ from \cite{JPP23}.
        
        \item (\cref{cor-Mikhalkincount}) When $k=\mathbb{C}$, the rank of our count equals the complex count of tropical curves pioneered by Mikhalkin \cite{Mi03}, see \cref{eq:MikhalkinC}, 
        $$\operatorname{rk}\left( N_\Delta^{\AA^1,\trop}(r,(d_1,\ldots,d_s))\right)=N_\Delta^{\trop}.$$
        
        \item (\cref{prop-Shustincount}) When $k=\mathbb{R}$ and all the $d_i<0$, the signature of our count equals the tropical Welschinger invariant for $s$ complex conjugate points defined by Shustin \cite{Shu06b}, see \cref{eq:Welschingerinvariant2},
        $$\operatorname{sgn}\left(N_\Delta^{\AA^1,\trop}(r,(d_1,\ldots,d_s))\right)=W_{\Delta,s}^{\trop}.$$
        
        \item (\cref{cor-MikhalkincountR}) For $k = \RR$ and all the $d_i>0$, the signature of our count equals the Welschinger invariant from \cref{eq:MikhalkinR}, 
        $$\operatorname{sgn}\left(N_\Delta^{\AA^1,\trop}(r,(d_1,\ldots,d_s))\right)=W_\Delta^{\trop}.$$
    \end{enumerate}
    In particular, our tropical count $N_\Delta^{\AA^1,\trop}(r,(d_1,\ldots,d_s))$ simultaneously recovers both $W_{\Delta, s}$ for all $s$ and~$N_\Delta$.
\end{theorem}

Shustin's correspondence theorem \cite{Shu06b} is in some sense the special case of \cref{thm-main} for~$k=\mathbb{R}$ and $d_i<0$ and our strategy to prove \cref{thm-main} is very similar to \cite{Shu06b}. However, compared to Shustin's theorem, we have many more local cases to deal with. This is due to certain tropical curves providing a total signed count of zero for their \emph{real} preimages under degeneration, but not a trivial quadratically enriched count. These local computations take place in Sections~\ref{section:computations local pieces}, \ref{sec-edges} and~\ref{sec-twintree} and the proof of \cref{thm-main} is presented in \cref{sec-corres}. 
We note that some of the local cases in Sections~\ref{section:computations local pieces} and~\ref{sec-edges} have already been studied in the context of the first quadratically enriched correspondence theorem \cite{JPP23}, but most of them are new.

One of the reasons for the success of Mikhalkin's original correspondence theorem \cite{Mi03} is that it provides a very efficient algorithm for computing $N_d$ and $W_d$. In particular the development of floor diagram techniques has contributed to this. With \cref{thm-main}, a similar approach to computing $N_\Delta^{\AA^1}(r,(d_1,\ldots,d_s))$ becomes possible, see \cref{sec-floor}.

\begin{theorem}[\cref{thm-flooreqtrop}] \label{thm-flooreqtrop_intro}
    The count of rational floor diagrams of degree $\Delta$ with $s$ merged points with quadratically enriched multiplicity equals the count of rational tropical stable maps of degree $\Delta$ satisfying $r$ simple and $s$ double point conditions, where $r+2s=n(\Delta)-1$:
    $$N^{\A^1,\trop}_\Delta(r,(d_1,\ldots,d_s)) = N^{\A^1,\floor}_\Delta(r,(d_1,\ldots,d_s)) \coloneq \sum_\cD \mult^{\AA^1} (\cD).$$
\end{theorem}

Using the floor diagram approach, we have computed $N_\Delta^{\AA^1}(r,(d_1,\ldots,d_s))$ for small degrees~$\Delta$, see \cref{sec-computations}.
For the case of $S = \mathbb{P}^2$ and degree $d$, the numbers can be computed with a wall-crossing formula as proved recently by the first author \cite{JP24}, but for other $S$ the numbers we provide are new. In particular, the count for bidegree $(2,4)$ in $\mathbb{P}^1\times\mathbb{P}^1$ and any number of conjugate pairs of points was not known before.

\subsection{A general formula $N^{\mathbb{A}^1}_\Delta(\sigma)$ for arbitrary $\sigma$}
\label{sec-generalformula}

Combining our Main Theorem \ref{thm-main} with very recent work of Brugallé--Rau--Wickelgren
\cite{BrugalleRauWickelgren}, one can compute $N^{\A^1}_\Delta(\sigma)$ for
arbitrary $\sigma$, not only for multiquadratic ones.
More precisely, Brugallé--Rau--Wickelgren show that
$N^{\A^1}_\Delta(\sigma)$ defines a \emph{Witt invariant} in the sense of
Serre \cite{GMS}. We briefly recall this notion.

The assignment sending a field $k$ to its Witt ring $\W(k)$ may be viewed as a
functor
\[
\W \colon \mathbf{Fields} \longrightarrow \mathbf{Set},
\]
which in fact factors through the category of rings.
Indeed, a morphism of fields
$\phi \colon k \to L$
induces a homomorphism
\[
\W(k) \longrightarrow \W(L)
\]
via extension of scalars.

Fix an integer $n \ge 1$. Define another functor
\[
Et_n \colon \mathbf{Fields} \longrightarrow \mathbf{Set}
\]
by sending a field $k$ to the set of isomorphism classes of rank $n$ étale
$k$-algebras. On morphisms, the functor is defined by tensoring with $L$ over $k$.

An \emph{étale Witt invariant of degree $n$} (in the sense of Serre) is a
natural transformation
\[
Et_n \longrightarrow \W .
\]

Being a Witt invariant imposes strong restrictions.
\begin{enumerate}
\item By \cite[Theorem~29.1]{GMS}, a Witt invariant is uniquely determined by
its values on multiquadratic étale algebras, i.e.\ algebras of the form
\[
\frac{k[x]}{(x^2-d_1)}\times \cdots \times \frac{k[x]}{(x^2-d_s)}\times k
\quad\text{if $n$ is odd},
\]
and
\[
\frac{k[x]}{(x^2-d_1)}\times \cdots \times \frac{k[x]}{(x^2-d_s)}
\quad\text{if $n$ is even},
\]
where $s=\lfloor \frac{n}{2} \rfloor$.

\item Serre also shows that the module $\mathrm{Inv}_k(n)$ of Witt invariants
is a free $\W(k)$-module of rank $s+1$.
\end{enumerate}

By \cite[Theorem~2.5]{BrugalleRauWickelgren} (see also
\cite[Theorem~27.16]{GMS}), the symmetric polynomials
$\beta_s^{(0)}$, $\ldots$, $\beta_s^{(s)}$
in the elements
\[
\beta_i=\Tr_{k[x]/(x^2-d_i)/k}(\gw{1}), \qquad i=1,\ldots,s,
\]
form a $\W(k)$-basis of $\mathrm{Inv}_k(n)$.

Brugallé--Rau--Wickelgren prove that $N^{\A^1}_\Delta(\sigma)$ defines a Witt
invariant in this sense \cite{BrugalleRauWickelgren}: to an étale $k$-algebra
determined by $\sigma$ one associates the image of
$N^{\A^1}_\Delta(\sigma)$ in $\W(k)$.

Combining this with Serre's theorem implies that the image of
$N^{\A^1}_\Delta(\sigma)$ in $\W(k)$ is completely determined, for arbitrary
$\sigma$, by its values on multiquadratic étale algebras. These are precisely
the cases computed by our Main Theorem \ref{thm-main}.

Moreover, Brugallé--Rau--Wickelgren show that the image of
$N^{\A^1}_\Delta(\sigma)$ in $\W(k)$ is not merely a $\W(k)$-linear
combination of the basis elements
$\beta_s^{(0)},\ldots,\beta_s^{(s)}$, but in fact a $\mathbb{Z}$-linear
combination. This phenomenon is also reflected in all of our explicit
computations in Section~\ref{sec-computations}, and their proof of $\mathbb{Z}$-linearity in fact relies on our computations in this paper.

Since an element of $\GW(k)$ is determined by its rank together with its image
in $\W(k)$, and since the rank equals the corresponding enumerative count over
$\mathbb{C}$ (which does not depend on the choice of $\sigma$), it follows that
our Main Theorem \ref{thm-main} determines
$N^{\A^1}_\Delta(\sigma)$ for arbitrary $\sigma$.

In particular, in Section~\ref{sec-computations} we provide explicit formulas
for $N^{\A^1}_\Delta(\sigma)$ that are valid for arbitrary $\sigma$ in all of
our examples. These formulas are correct because they agree with the values
given by our Main Theorem for all multiquadratic choices of $\sigma$.

As an illustration of the dependence of $N^{\A^1}_\Delta(\sigma)$ on $\sigma$,
consider the example of degree~$3$ rational plane curves. The quadratically
enriched count of such curves with arbitrary point conditions $\sigma$ was
computed in \cite{KassLevineSolomonWickelgren}. Using different methods, the
authors show that
\[
N^{\A^1}_{\Delta_3}(\sigma)
=
\Tr_{k(\sigma)/k}(\gw{1}) + 2\h
\]
for every choice of $\sigma$. Here $k(\sigma)$ denotes the product of the
fields appearing in $\sigma$, and the formula holds for any étale $k$-algebra
$\sigma$ with $\dim_k \sigma = 8$.

\subsection{Possible generalizations}

Conceptually, it would be interesting to extend the Main Theorem~\ref{thm-main} to more general sequences $\sigma$, beyond just multiquadratic ones. One possible approach is to extend the scope of the correspondence theorem by merging more than two points on the tropical side and develop the correct multiplicities for the new tropical curves that arise. The first steps involve considering small-degree field extensions in $\sigma$ where only a few points are moved together. 

 A less involved approach to compute $N^{\AA^1}_\Delta(\sigma)$ for more general $\sigma$ is to find ``wall-crossing formulas" as in \cite{JP24}.

Throughout the paper we work with vertically stretched point conditions on the tropical side. This does not restrict generality, as the quadratically enriched count is an invariant, but it allows to reduce the combinatorial challenges on the tropical side. Especially the computations in \cref{section:computations local pieces} are simplified by the vertically stretched assumption. In principle, we do not believe that the assumption is strictly necessary for the proof.


For base fields $k$ of positive characteristic, our proof only holds if the characteristic of $k$ is greater than the diameter of $\Delta$. This restriction is due to the structure of the proof and cannot be lifted easily.

Finally, the tropical curves we consider in this paper have unbounded ends of weight 1 only, corresponding to the algebraic curves meeting the toric boundary of $S$ transversely. Empirically, we expect the techniques of logarithmic geometry to give a path for generalizing this. However, for this one would first have to lift \cite{KassLevineSolomonWickelgren} to log geometry and establish invariance of the generalized algebraic count.

\subsection{Historical developments in counting rational plane (tropical) curves}
\label{sec-history}
Let $k$ be an algebraically closed base field and recall that $N_d$ denotes the count of rational plane curves of degree~$d$ through the correct number of points.
The easiest case is, of course, the number $N_1=1$, counting the number of lines through $2$ points. Until the 1990s, the numbers $N_d$ were known only for small values of $d$. Kontsevich \cite{KM94}, inspired by connections to theoretical physics and mirror symmetry, then provided a recursive formula that computes all $N_d$.

Around 2002, the development of tropical geometry provided a new tool for enumerative problems, pioneered by Mikhalkin's celebtrated correspondence theorem. Mikhalkin showed that for any smooth toric surface $\Delta$, $N_\Delta$ is equal to its counterpart in tropical geometry \cite{Mi03, Shu04, Tyo09, AB22}, that is
\begin{equation}
    \label{eq:MikhalkinC}
    N_\Delta=N_\Delta^{\trop}\coloneqq \sum_\Gamma \operatorname{mult}_\CC(\Gamma),
\end{equation}
where we sum over certain \emph{plane tropical curves} satisfying point conditions.
Tropical geometry can be viewed as a refined degeneration technique and the multiplicity $\mult_\CC(\Gamma)$ used in Mikhalkin's correspondence theorem \eqref{eq:MikhalkinC} reflects how many algebraic curves degenerate to a given tropical curve $\Gamma$.  
The correspondence theorem promises an easier solution to the algebro-geometric enumerative problem, since tropical plane curves can simply be viewed as piecewise linear graphs in the plane. 
Indeed, the right hand side of Equation~\eqref{eq:MikhalkinC} can be computed efficiently.
Moreover, for~$\mathbb{P}^2$, not only the numbers $N_d$, but also their recursive structure described by Kontsevich, can be translated into the tropical world \cite{GM053}. 
 
Moving on to base field $k = \RR$, Welschinger observed the following. A real curve in a surface can have two types of real nodes. A node is called \emph{hyperbolic}, if its branches are defined over $\RR$, and \emph{elliptic}, if its branches are not defined over $\RR$. Welschinger defined
\begin{equation}
    \label{eq:Welschingerinvariant}
    W_\Delta\coloneqq \sum_{\text{real curves}} (-1)^{\#\text{elliptic nodes}}
\end{equation}
and proved that $W_\Delta$ is an invariant, which is independent of the chosen point configuration (as long as it is generic). Here, $\Delta$ is the Newton polygon of a smooth toric del Pezzo surface as in Figure \ref{fig-degrees}.
The tropical approach to counting curves allows to compute $W_\Delta$ as well, i.e. there is another correspondence theorem: 
\begin{equation}
    \label{eq:MikhalkinR}
    W_\Delta=W_\Delta^{\trop}\coloneqq\sum_\Gamma\mult_\RR(\Gamma).
\end{equation}
We point out that the summation on the right hand side of Equations~\eqref{eq:MikhalkinC} and~\eqref{eq:MikhalkinR} runs over the same set of tropical curves and it is only through the use of different multiplicities that we compute~$N_\Delta$ in one case and $W_\Delta$ in the other. 
Important properties of Welschinger invariants, such as their logarithmic equivalence with the numbers $N_\Delta$, were proved by means of tropical geometry \cite{IKS03, IKS04}.
 
Generalizing the real curve counting problem, we let $s$ be the number of pairs of complex conjugate point conditions. Welschinger proved that 
\begin{equation}
    \label{eq:Welschingerinvariant2}
    W_{\Delta,s}\coloneqq \sum_{\text{real curves}}(-1)^{\#\text{elliptic nodes}}
\end{equation} 
is still an invariant and this recovers $W_\Delta$ for $s = 0$. 
Shustin investigated the tropical approach to this counting problem \cite{Shu06b} and it turns out that the combinatorial part is considerably more involved in the generalized setting. 
Nonetheless, it is still possible to describe the tropical curves which occur in the tropical analogue of the counting problem and to define a tropical multiplicity $\mult_{\RR}$ for the generalized problem which leads to a correspondence theorem
\begin{equation}
    \label{eq:Shustin}
    W_{\Delta,s}=W_{\Delta,s}^{\trop}\coloneqq \sum_\Gamma \operatorname{mult}_\RR(\Gamma).
\end{equation}
We emphasize that in \cref{eq:Shustin} the sum is taken over a different set of tropical curves compared to Mikhalkin's correspondence theorems in Equations~\eqref{eq:MikhalkinC} and~\eqref{eq:MikhalkinR}. In particular, the tropical curves which appear in \cref{eq:Shustin} are not necessarily \emph{simple} anymore.
The reason is that two complex conjugate points are mapped to the same tropical point under the tropical degeneration procedure. This leads to point conditions in the tropical counting problem which behave as \emph{double points}. 
 
All of the above curve counting theories are unified in the quadratically enriched counting theory in the sense that the quadratically enriched counts solve the geometric counting problem simultaneously over a large class of fields, including $\CC$ and $\RR$. 

Working towards a tropical approach to the quadratically enriched count, the first and third author developed a quadratic enrichment of Mikhalkin's correspondence theorem \cite{JPP23}. 
More precisely, in the special case of where all point conditions are defined over $k$, they provide a Grothendieck-Witt valued multiplicity $\mult^{\A^1}(\Gamma)$ for a tropical curve $\Gamma$ satisfying the point conditions and they show 
\begin{equation}
    N_\Delta^{\AA^1}=N_\Delta^{\AA^1,\trop}\coloneqq\sum_\Gamma\mult^{\AA^1}(\Gamma).\label{eq-quadraticallyenrichedcount}
\end{equation}
Moreover, the new quadratically enriched multiplicity for tropical curves specializes to the real multiplicity contributing to the Welschinger invariant when the ground field is $k = \mathbb{R}$ and to the complex multiplicity contributing to $N_\Delta$ when the ground field is $k = \mathbb{C}$. 
Note that the sum in the right hand side of \cref{eq-quadraticallyenrichedcount} runs again over the same tropical curves $\Gamma$ as in Equations~\eqref{eq:MikhalkinC} and~\eqref{eq:MikhalkinR}. This can be attributed to all of these theorems dealing with $k$-rational point conditions.
The structure of the result in \eqref{eq-quadraticallyenrichedcount} is similar to Shustin's version of the correspondence theorem for the complex count $N_\Delta$ \cite{Shu04, IMS09}, in the sense that the case-by-case analysis of the local pieces of tropical curves which have to be considered is the same; however, the construction of the quadratically enriched multiplicity is naturally more involved. 
Our current \cref{thm-main} generalizes \eqref{eq-quadraticallyenrichedcount} and \eqref{eq:Shustin} simultaneously and follows these ideas as well, however the combinatorial and computational complexity increases dramatically.

\subsection{The tropical approach to quadratically enriched counting}

Since both tropical geometry as well as quadratic enrichment allow to treat the count of complex and real rational plane curves simultaneously, a possible connection of tropical and arithmetic methods seems possible.
Indeed, three of the four authors were involved in first publications mingling tropical and arithmetic tools: Payne and Shaw with the second author investigated the effect of tropical degeneration for the quadratically enriched count of bitangents to a plane quartic \cite{MPS22}; and the first and the third author studied the tropical version of the quadratically enriched B\'ezout's theorem \cite{puentes2022quadratically}. At the base of combining tropical degenerations with quadratic enrichments lies the fact that the Grothendieck-Witt ring of a non-Archimedean field (which is used for degeneration) equals the Grothendieck-Witt ring of its residue field \cite{MPS22}.
With the quadratically enriched correspondence theorem \cite{JPP23}, the first and third author strengthened the connection between tropical geometry and quadratically enriched counting further. 
Indeed, the correspondence theorem already led to new insights into the combinatorial properties of the quadratically enriched (tropical) plane curve count in \cite{JPMPR23}. 
With this project, we generalize \cite{JPP23} and expand the tropical approach to computing $N_\Delta^{\AA^1}$ for a larger class of point conditions.

\subsection{Acknowledgments}
We thank Jesse Pajwani for pointing out that identities in the Gro\-then\-dieck-Witt ring can be checked on multiquadratic finite étale algebras and for his help with computations. We thank Erwan Brugallé and Kirsten Wickelgren for useful discussions and explanations on Witt invariants. We also thank Johannes Rau for corrections of the computations.
We thank Andreas Gross, Marc Levine, Dhruv Ranganathan and Vincenzo Reda for useful discussions.
The first, second and fourth author acknowledge support by DFG-grant MA 4797/9-1. 
The third author acknowledges support by Deutsche Forschungsgemeinschaft (DFG, German Research Foundation) through the Colla\-borative Research Centre TRR 326 \textit{Geometry and Arithmetic of Uniformized Structures}, project number 444845124.
The first author thanks the Universität Duisburg-Essen and the Università degli Studi di Napoli Federico II  for support. 
Part of this work was completed while the authors stayed as Research Fellows at the Mathematisches Forschungsinstitut Oberwolfach in March 2024. We thank the institute for hosting us and for providing ideal working conditions.

\section{Preliminaries}
\label{sec-preliminaries}

\subsection{Embedded plane tropical curves and dual subdivisions}
We start with a brief overview of the theory of tropical plane curves. The exposition in this subsection is based on Chapter 3 in \cite{CMR23}, more details can be found there.

Consider a field $K$. A plane algebraic curve (or, more generally, a curve embedded in a toric surface) over $K$ is given as the zero-set of a polynomial in $K[x, y]$. In tropical geometry we can use the analogous definition, by just adding the word \enquote{tropical} everywhere:
tropical plane curves are tropical vanishing loci of tropical polynomials over the tropical semifield. As this perspective is not of particular relevance for our purpose, we just refer interested readers to surveys on tropical geometry such as \cite{RGST05, BS14}. 

Alternatively, one can also start with an (algebraically closed) field $K$ which is equipped with a non-Archimedean valuation $\val : K^\times \to \RR$ and an algebraic curve $C$ defined over $K$. 

We will work in the field of Puiseux series over the algebraic closure of $k$.

The \emph{tropicalization} of $C$ is the closure of the image of $C \cap (K^\times)^2$ under the \emph{tropicalization map}
\[ \trop : (K^\times)^2 \longrightarrow \RR^2, \qquad (x,y) \longmapsto \big( -\val(x), -\val(y) \big) \]
and this is a tropical plane curve. It is a piecewise linear graph with edges of rational slope in $\mathbb{R}^2$. In fact, the two definitions mentioned here are equivalent by Kapranov's theorem.

Every tropical plane curve is dual to a Newton subdivision of the Newton polygon of its defining equation. The subdivision is given by projecting upper faces of the extended Newton polygon, which takes the valuations of the coefficients into account.

\subsection{Rational tropical stable maps}
 A \emph{rational (abstract) tropical curve} is a metric tree $\Gamma$ with unbounded rays, which are called \emph{ends}. Locally around a point $p \in \Gamma$, the curve $\Gamma$ is homeomorphic to a star with some number of half-rays, and this number is called the \emph{valence} $\val(p)$ of the point $p$.
We require that there are only finitely many points with $\val(p)\neq 2$, these points are called  \emph{vertices}.
By abuse of notation, the underlying graph with this vertex set is also denoted by $\Gamma$.
Correspondingly, we can speak about \emph{edges} and \emph{flags} of $\Gamma$. A flag is a tuple $(v,e)$ of a vertex $v$ and an edge $e$ with $v\in \partial e$. It can be thought of as an element in the tangent space of $\Gamma$ at $v$, i.e.\ as a germ of an edge leaving $v$, or as a half-edge (the half of $e$ that is attached to $v$). Edges which are not ends have a finite length and are called \emph{bounded edges}.

A \emph{marked rational tropical curve} is a rational tropical curve such that some of its ends are labeled. An isomorphism of marked rational tropical curves is a homeomorphism preserving the lengths of the edges and the markings of ends. 
The \emph{combinatorial type} of a tropical curve is obtained by dropping the information on the metric. 

\begin{definition}
    A \emph{rational tropical stable map to $\mathbb{R}^2$} is a tuple $(\Gamma,f)$ where $\Gamma$ is a marked rational abstract tropical curve and $f:\Gamma\to \mathbb{R}^2$ is a piecewise integer-affine map satisfying:
    \begin{itemize}
        \item On each edge $e$ of $\Gamma$, $f$ is of the form 
        $$t\longmapsto a+t\cdot \direction \mbox{ for some } \direction \in \ZZ^2,$$ 
        where we parametrize $e$ as an interval of size the length $l(e)$ of $e$. The vector $\direction$, called the \emph{direction}, arising in this equation is defined up to sign, depending on the starting vertex of the parametrization of the edge. We will sometimes speak of the direction of a flag $\direction(v,e)$. If $e$ is an end we use the notation $\direction(e)$ for the direction of its unique flag.
        \item The \emph{balancing condition} holds at every vertex, i.e.\ 
        $$\sum_{e \text{ with } v \in \partial e} \direction(v,e)=0.$$
    
    \end{itemize}
\end{definition}
For an edge with direction $\direction=(\mathbf{v}_1,\mathbf{v}_2) \in \ZZ^2$, 
we call $m=\gcd(\mathbf{v}_1,\mathbf{v}_2)$ the \emph{weight} or \emph{expansion factor} and $\frac{1}{m}\cdot \direction$ the \emph{primitive direction} of $e$.

An isomorphism of tropical stable maps is an isomorphism of the underlying tropical curves respecting the maps. 

\begin{remark} The image $f(\Gamma)$ in $\mathbb{R}^2$ of a tropical stable map $(\Gamma,f)$ is an embedded plane tropical curve, as considered before, and as such, dual to a Newton subdivision.
\end{remark}

\begin{definition}[Degree of a tropical stable map]\label{def-deg}
    The \emph{degree} $\Delta$ of a tropical stable map is the multiset of the non-zero directions of its ends.   
\end{definition}

\begin{figure}[t]
    \centering

    \tikzset{every picture/.style={line width=0.75pt}} 
    
    \begin{tikzpicture}[x=0.75pt,y=0.75pt,yscale=-1,xscale=1]
        
        \draw    (100,320) -- (160,320) ;
        \draw    (100,260) -- (160,320) ;
        \draw    (100,260) -- (100,320) ;
        \draw    (190,260) -- (190,320) ;
        \draw    (220,260) -- (220,320) ;
        \draw    (190,320) -- (220,320) ;
        \draw    (190,260) -- (220,260) ;
        \draw    (250,260) -- (250,320) ;
        \draw    (250,320) -- (280,320) ;
        \draw    (250,260) -- (280,290) ;
        \draw    (280,290) -- (280,320) ;
        \draw    (349,260) -- (349,300) ;
        \draw    (369,320) -- (379,320) ;
        \draw    (349,260) -- (379,290) ;
        \draw    (379,290) -- (379,320) ;
        \draw    (349,300) -- (369,320) ;
        \draw    (470,270) -- (470,300) ;
        \draw    (490,320) -- (500,320) ;
        \draw    (480,270) -- (500,290) ;
        \draw    (500,290) -- (500,320) ;
        \draw    (470,300) -- (490,320) ;
        \draw    (470,270) -- (480,270) ;
        
        \draw (121,322) node [anchor=north west][inner sep=0.75pt]   [align=left] {$\displaystyle d$};
        \draw (201,322) node [anchor=north west][inner sep=0.75pt]   [align=left] {$\displaystyle a_{2}$};
        \draw (171,280) node [anchor=north west][inner sep=0.75pt]   [align=left] {$\displaystyle a_{1}$};
        \draw (282,293) node [anchor=north west][inner sep=0.75pt]   [align=left] {$\displaystyle a_{1}$};
        \draw (88.33,282.86) node [anchor=north west][inner sep=0.75pt]   [align=left] {$\displaystyle d$};
        \draw (239,282) node [anchor=north west][inner sep=0.75pt]   [align=left] {$\displaystyle d$};
        \draw (246,322) node [anchor=north west][inner sep=0.75pt]   [align=left] {$\displaystyle d-a_{1}$};
        \draw (343.86,305.71) node [anchor=north west][inner sep=0.75pt]   [align=left] {$\displaystyle a_{2}$};
        \draw (381,292) node [anchor=north west][inner sep=0.75pt]   [align=left] {$\displaystyle a_{1}$};
        \draw (491,262) node [anchor=north west][inner sep=0.75pt]   [align=left] {$\displaystyle d-a_{1} -a_{3}$};
        \draw (361,260) node [anchor=north west][inner sep=0.75pt]   [align=left] {$\displaystyle d-a_{1}$};
        \draw (502,293) node [anchor=north west][inner sep=0.75pt]   [align=left] {$\displaystyle a_{1}$};
        \draw (465,304.86) node [anchor=north west][inner sep=0.75pt]   [align=left] {$\displaystyle a_{2}$};
        \draw (462,252) node [anchor=north west][inner sep=0.75pt]   [align=left] {$\displaystyle a_{3}$};
        \draw (306,272) node [anchor=north west][inner sep=0.75pt]   [align=left] {$\displaystyle d-a_{2}$}; 
    
    \end{tikzpicture}

    \caption{The dual polygons of degrees of smooth toric del Pezzo  surfaces.}
    \label{fig-degrees}
\end{figure}
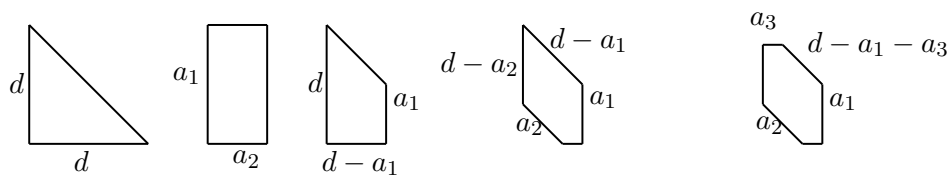

Throughout the paper, we restrict to the case where the degree $\Delta$ consists of normal vectors to a polygon defining a smooth del Pezzo surface, together with a curve class given by the hyperplane section. In this case, we have that $n(\Delta)$ coincides with the number of lattice points on the boundary of the polygon. All possible polygons are depicted in Figure \ref{fig-degrees}. Here, we follow the notation used in \cite[Section 2.2]{IKS09}. For each of these types, we use the following more specific notation. 

\begin{itemize}
    \item If $\Delta$ is the multiset consisting of the vectors $(-1,0)$, $(0,-1)$ and $(1,1)$ each $d$ times, we say the tropical stable map is of \emph{degree $d$}.

    \item If $\Delta$ is the multiset consisting of the vectors $(-1,0)$ and $(1,0)$, each $a_1$ times, and $(0,-1)$ and $(0,1)$ each $a_2$ times, we say that the tropical stable map is of \emph{bidegree $(a_1,a_2)$}. 

    \item For the blowup of the plane in one point, we use the notation $(d;a_1)$ to fix the degree, as depicted in Figure \ref{fig-degrees}. That is, we have $d$ ends of direction $(-1,0)$, $d-a_1$ ends of direction $(0,-1)$ resp.\ $(1,1)$ and $a_1$ ends of direction $(1,0)$.

    \item For the blowup of the plane in two points, we use the notation $(d;a_1,a_2)$ for $d-a_2$ ends in direction $(-1,0)$, $a_2$ ends in direction $(-1,-1)$, $d-a_1-a_2$ ends in direction $(0,-1)$, $a_1$ ends in direction $(1,0)$ and $d-a_1$ ends in direction $(1,1)$.

    \item Finally, for the blowup of the plane in three points, we write $(d;a_1,a_2,a_3)$ for $d-a_2-a_3$ ends of direction $(-1,0)$, $a_2$ ends in direction $(-1,-1)$, $d-a_1-a_2$ ends in direction $(0,-1)$, $a_1$ ends in direction $(1,0)$ and $d-a_1-a_3$ ends in direction $(1,1)$ and $a_3$ ends in direction $(0,1)$.
\end{itemize}

The \emph{combinatorial type} of a tropical stable map is the data obtained when dropping the metric of the underlying graph. More explicitly, it consists of the data of a finite graph $\Gamma$, and for each edge~$e$ of~$\Gamma$, the weight and primitive direction of $e$.

We often fix the convention that the marked ends of a rational tropical stable map have direction~$0$, i.e.\ they are contracted to a point.
Assume we have $n$ marked ends $x_1,\ldots,x_n$.
In our enumerative problem, we can then fix points $p_1,\ldots,p_n\in\mathbb{R}^2$ and consider tropical stable maps $(\Gamma,f)$ which satisfy $f(x_i)=p_i$ for all $i$. 

Fix a degree $\Delta$, and fix $n(\Delta)$ points $p_1,\ldots,p_n$ in tropical general position in $\mathbb{R}^2$ (see Definition~4.7 in \cite{Mi03} resp.\ Definition 5.33 in \cite{Ma06}).

For some piecewise linear graphs in $\RR^2$ (resp.\ embedded tropical plane curves), there are several ways to parametrize them by an abstract tropical curve. So-called simple tropical curves can uniquely be parametrized, once we fix the convention that every point dual to a parallelogram should be viewed as a crossing of two edges and not as a vertex:

\begin{definition}[Simple tropical stable maps]
A tropical stable map $(\Gamma,f)$ is called \emph{simple} if $\Gamma$ is $3$-valent, all ends are of weight $1$ and the Newton subdivision dual to $f(\Gamma)\subset \RR^2$ contains only triangles and parallelograms. 
\end{definition}

As the following proposition shows, restricting to simple tropical stable maps is sufficient for the purpose of counting curves satisfying point conditions.

\begin{proposition}\label{prop-genericsimple}
Fix $n(\Delta)$ points $p_1,\ldots,p_n\in\mathbb{R}^2$ in tropical general position. Then any rational tropical stable map $(\Gamma,f)$ which satisfies $f(x_i)=p_i$ for all $i = 1, \ldots, n$, where the $x_i$ are the marked ends, is simple.
\end{proposition}
For a proof, see e.g.\ Lemma 5.34 in \cite{Ma06}.

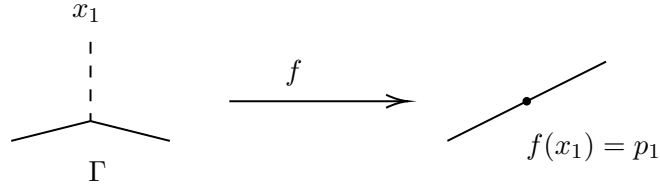
\begin{figure}[t]
    \begin{center}
      
    \tikzset{every picture/.style={line width=0.75pt}} 
    
    \begin{tikzpicture}[x=0.75pt,y=0.75pt,yscale=-1,xscale=1]
    
    \draw  [dash pattern={on 4.5pt off 4.5pt}]  (200,60) -- (200,100) ;
    \draw    (160,110) -- (200,100) ;
    \draw    (200,100) -- (240,110) ;
    \draw    (270,90) -- (358,90) ;
    \draw [shift={(360,90)}, rotate = 180] [color={rgb, 255:red, 0; green, 0; blue, 0 }  ][line width=0.75]    (10.93,-3.29) .. controls (6.95,-1.4) and (3.31,-0.3) .. (0,0) .. controls (3.31,0.3) and (6.95,1.4) .. (10.93,3.29)   ;
    \draw    (380,110) -- (460,70) ;
    \draw  [fill={rgb, 255:red, 0; green, 0; blue, 0 }  ,fill opacity=1 ] (418.35,90) .. controls (418.35,89.09) and (419.09,88.35) .. (420,88.35) .. controls (420.91,88.35) and (421.65,89.09) .. (421.65,90) .. controls (421.65,90.91) and (420.91,91.65) .. (420,91.65) .. controls (419.09,91.65) and (418.35,90.91) .. (418.35,90) -- cycle ;
    
    \draw (200,45.5) node   [align=left] {\begin{minipage}[lt]{13.94pt}\setlength\topsep{0pt}
    $\displaystyle x_{1}$
    \end{minipage}};
    \draw (204.75,125.25) node   [align=left] {\begin{minipage}[lt]{8.67pt}\setlength\topsep{0pt}
    $\displaystyle \Gamma $
    \end{minipage}};
    \draw (304,75.5) node   [align=left] {\begin{minipage}[lt]{8.67pt}\setlength\topsep{0pt}
    $\displaystyle f$
    \end{minipage}};
    \draw (465,112.8) node   [align=left] {\begin{minipage}[lt]{68pt}\setlength\topsep{0pt}
    $\displaystyle f( x_{1}) =p_{1}$
    \end{minipage}};
    \end{tikzpicture}
    \end{center}
    \caption{A local picture of a marked end of a tropical stable map and its image in $\RR^2$. We use the convention to draw marked ends as dashed lines.}\label{fig-pointcondition}
\end{figure}

In particular, a rational tropical stable map satisfying generic point conditions is $3$-valent and locally around a marked end the image looks like a point in the interior of a bounded edge, see Figure \ref{fig-pointcondition}.

\subsection{Moduli spaces of rational tropical stable maps}\label{sec-tropmoduli}

Given a marked abstract tropical curve $\Gamma$, a degree $\Delta$ and the position of one marked end, say $f(x_1)$, the map $f$ can be fully recovered from this data, see \cite{GKM07}, Lemma 4.6. The space of all rational tropical stable maps of a fixed degree can therefore be interpreted as a product of $\mathbb{R}^2$ (for the position of $f(x_1)$) and the space of all abstract rational tropical curves. The latter is also known as the \emph{space of trees} \cite{BHV01, SS04a}. It can be viewed as an abstract cone complex whose cones are indexed by combinatorial types as follows. Given a combinatorial type with $b$ bounded edges, we can vary the lengths of these bounded edges; that is, we obtain an orthant $\mathbb{R}_{>0}^b$ worth of abstract tropical curves of this type. If we shrink an edge to length zero, this gives a new combinatorial type and the corresponding orthant $\RR_{>0}^{b-1}$ of this new combinatorial type can be identified with a face of the orthant of the original combinatorial type. In this way, the orthants can be glued to form an abstract cone complex. This cone complex can also be embedded as a fan into a real vector space, using the distance map which associates the tuple of all distances between pairs of marked ends to an abstract tropical curve, see \cite{SS04a, GKM07}. The embedded fan satisfies the so-called \emph{balancing condition} (i.e.\ the weighted primitive generators of the neighboring top-dimensional cones of each cone $\tau$ of codimension one add up to zero modulo the span of $\tau$). Such \emph{balanced fans} are the local building blocks for tropical varieties.

We call the tropical variety that is 
 the \emph{space of marked rational tropical stable maps} of degree $\Delta$ with~$n$ contracted ends $M_{0,n}(\mathbb{R}^2,\Delta)$. 

For each $i=1,\ldots,n$, there is an evaluation map
$$\ev_i:M_{0,n}(\mathbb{R}^2,\Delta) \longrightarrow \mathbb{R}^2, \qquad (\Gamma,f)\longmapsto f(x_i)$$
which sends a tropical stable map to the evaluation of its $i$-th contracted end.
If $n=n(\Delta)$, the product $\ev:=\ev_1\times\ldots\times\ev_n$ is a morphism of fans of the same dimension (see \cite{GKM07}, Proposition 4.7). 

To count tropical plane curves satisfying point conditions, we can consider the inverse image of a tuple $(p_1,\ldots,p_n)\in(\mathbb{R}^2)^n$ under the evaluation map $\ev$. If we pick the $p_i$ in general position, our preimages in $M_{0,n}(\mathbb{R}^2,\Delta)$ are simple tropical stable maps $(\Gamma,f)$ and we may identify them with their images $f(\Gamma)\subset \mathbb{R}^2$. As is common in tropical geometry, we do not count tropical objects one by one, but with a \emph{multiplicity} which reflects how many algebraic curves are degenerating to one tropical curve. We will discuss multiplicities later (see Section \ref{section:computations}).

As in algebraic geometry, inverse images can be described more accurately in terms of \emph{pull-backs} of morphisms. For an introduction to the theory of tropical divisors and cycles, their pull-backs and push-forwards (i.e.\ to \emph{tropical intersection theory}), see \cite{AR07}.
Using tropical intersection theory, the preimages in the pull-back $\ev^\ast(p_1,\ldots,p_n)$ obtain an intrinsic multiplicity which equals $\mult_{\mathbb{C}}(\Gamma)$, i.e.\ by the correspondence theorem, the number of algebraic curves defined over an algebraically closed field that degenerate to a tropical curve \cite{GM053, Ma06, Mi03}. 

\subsection{Log stable map spaces and their tropicalization}\label{sec-square}

Counts of curves over $K$ in a toric surface~$X$ satisfying point conditions can be  defined in terms of the space of \emph{log stable maps} and the evaluation morphism. For more background, see e.g. \cite{Ran15} and references therein.
We work in this setting because of its strong ties to tropical geometry, in particular the commutative square~\eqref{eq:comm_square} below.
In principle, working with log stable maps allows to fix certain tangency behavior of the curves along the toric boundary of $X$. In our setting, we always fix trivial tangency behavior. We cannot work with higher tangencies since the quadratically enriched algebraic invariants only exist for trivial tangencies. We will always assume that our curves intersect the toric boundary generically. 
By Theorem~2.2.5 in \cite{Ran15} (see also \cite{AC14, GS13}), there is a moduli space $\mathcal{M}_{0,n}(X,\beta)$ of rational log stable maps of fixed curve class $\beta$ to a toric variety $X$. Its boundary consists of degenerate curves, which we  can tropicalize  by taking dual graphs (see Section 3.4 in \cite{Ran15} or Section~2.2.1 in \cite{Tyo09}). Counts of curves can be defined using the evaluation map which evaluates the marked points.
By Theorem~B in \cite{Ran15}, the moduli space of rational log stable maps of fixed curve class to a toric variety with $n$ marked points can be embedded into a toric variety which is given by the corresponding tropical moduli space $M_{0,n}(\mathbb{R}^2,\Delta)$ which we discussed in Section \ref{sec-tropmoduli}. The tropicalization in these coordinates then equals $M_{0,n}(\mathbb{R}^2,\Delta)$, and also coincides with the pointwise tropicalization of a log stable map mentioned before. Furthermore, we obtain the following commutative square (see Theorem C in \cite{Ran15} or Corollary 3.6 in \cite{Gro16}): 

\begin{equation} \label{eq:comm_square}
	  \begin{tikzcd}
		\mathcal{M}_{0,n}(X,\beta) \arrow{d}{\trop} \arrow{r}{\ev} & T^n \arrow{d}{\trop} \\
		M_{0,n}(\mathbb{R}^2,\Delta)  \arrow{r}{\ev} & (\mathbb{R}^2)^n.
	\end{tikzcd}
\end{equation}
Here, $T$ denotes the dense torus in the toric surface $X$.
(Note: \cite{Gro16} makes the global assumption that the characteristics of the ground field is $0$, but this is used only in a later chapter, not for this commutative square.)

\subsection{Definition of Gromov-Witten invariant $N_{S,\beta}$ and Welschinger invariant $W_{S,\beta}$}

Let $S$ be a del Pezzo surface over an algebraically closed field $k$ and $\beta$ an effective Cartier divisor. Set $n$ as~$-K_S\cdot \beta-1$. Then a dimension count shows that there are finitely many degree $\beta$ rational curves through $n$ points in general position on $S$. 
Let $N_{S,\beta}$ be this finite number. 

\begin{example}
    For $k = \CC$, $S=\PP^2_{\mathbb{C}}$, and $\beta=d\cdot H$ where $H$ is the class of a hyperplane, we have that $N_{S,\beta}=N_d$ from the introduction.
\end{example}
\begin{remark}
Instead of viewing the number $N_{S,\beta}$ as a count of curves satisfying certain conditions, one could also introduce these numbers in terms of Gromov-Witten or logarithmic Gromov-Witten theory. Some of the results we rely on (e.g.\ the results concerning the tropicalization of the log stable map space \cite{Ran15}) require the use of logarithmic Gromov-Witten invariants. Remark 5.1.1 of \cite{Ran15} shows that logarithmic Gromov-Witten invariants are enumerative (see also \cite{NS06, MR16,CMR23}), i.e.\ they correspond to the count of curves that we denoted by $N_{S,\beta}$ above. Other results we rely on (e.g.\ the $\GW(k)$-valued count of rational curves we discuss in Section \ref{GWk-counts} following
\cite{KassLevineSolomonWickelgren,KLSWOrientation}) are build on classical Gromov-Witten theory. In our setting of considering curves passing through point conditions in smooth del Pezzo surfaces, the genericity of the points implies that any curve passing through them meets the toric boundary generically and can thus be interpreted as an element in the log stable map space. Gromov-Witten invariants, logarithmic Gromov-Witten invariants and counts of curves satisfying point conditions thus coincide in our setting and we can pass from one description to the other where needed.
\end{remark}

Welschinger proved that one can count real rational curves on $S$ of class $\beta$ with a sign depending on the types of real nodes and get an invariant independent of the chosen point configuration \cite{Wel03,Wel05}. 
More precisely, a real curve in a surface can have two types of nodes defined over the reals called \emph{hyperbolic} and \emph{elliptic}. A real node is hyperbolic if its two branches are defined over~$\RR$ and elliptic if they are not. 
\begin{definition}
    The \emph{Welschinger sign} of a real curve $C$ in a surface is given by
    \[\operatorname{Wel}(C)\coloneqq (-1)^{\#\text{elliptic nodes}}.\]
\end{definition}
Counting real curves with this sign is an invariant: Consider a  generic point configuration of $r$ real points and $s$ pairs of complex conjugate points in $S$ with $r+2s=n$. Define 
\[W_{S,\beta}(r,s)\coloneqq \sum \operatorname{Wel}(C),\]
where the sum goes over all real curves in $S$ of class $\beta$ passing through the point configuration. This is independent of the chosen point configuration. Summing the Welschinger sign over all elements in the fiber of a point of the evaluation map from the moduli space of log stable maps, we can view the Welschinger invariant in the context of logarithmic Gromov-Witten theory.

\subsection{The Grothendieck-Witt ring}
Our enumerative invariants live in the Grothendieck-Witt ring $\GW(M)$ of a finite étale $k$-algebra $M$ where $k$ is a field of characteristic not equal to $2$. 
We recall the definition of $\GW(M)$ here as well as some computational results.
If one endows the set
\begin{equation*}
    \left\{\beta\colon V\times V\rightarrow M \ \middle\vert \ \begin{minipage}{8.7cm}
    $V$ is a finitely generated projective $M$-module and $\beta$ a non-degenerate symmetric bilinear form over $M$
\end{minipage}  \right\} /\text{isometry}
\end{equation*}
with the operations direct sum $\oplus$ and tensor product $\otimes$, one obtains the structure of a semiring.
\begin{definition}\label{df:gw}
    The \emph{Grothendieck-Witt ring} $\GW(M)$ of $M$ is the group completion of the semi-ring of isometry classes of non-degenerate symmetric bilinear forms over $M$ with respect to the operation~$\oplus$. 
\end{definition}

We note that $\chara k \neq 2$ and hence the elements of $\GW(M)$ can be viewed as either classes of bilinear forms or quadratic forms. We will move between these two formulations freely.

\begin{definition}
    Let $\beta\colon V\times V\rightarrow k$ be a non-degenerate symmetric bilinear form over a field $k$. Then the \emph{rank} of $\beta$ is given by $\rk \beta=\dim_kV$. This defines a homomorphism $\rk\colon \GW(k)\rightarrow \ZZ$. 
    Assume $k\subset \RR $ and $\beta\colon V\times V\rightarrow k$ is a non-degenerate symmetric bilinear form over $k$. Then one can diagonalize the Gram matrix of $\beta$, let $d_1,\ldots,d_n$ be the diagonal elements. The \emph{signature} of $\beta$ is given by $\#\{d_i:d_i>0\}-\#\{d_i:d_i<0\}$ and this defines a homomorphism $\sgn\colon \GW(k)\rightarrow \ZZ$.
\end{definition}
\begin{example}
    The Grothendieck-Witt ring of $\CC$ is $\GW(\CC) \cong \ZZ$ and the isomorphism is given by the rank. For $\RR$ we have $\GW(\RR) \cong \ZZ[\mu_2] \cong \gw{1}\ZZ + \gw{-1}\ZZ$.
\end{example}

Recall that a finite étale $k$-algebra $M$ is isomorphic to the product of finitely many finite separable field extensions $L_i$ of $k$
\[M\cong L_1\times\ldots \times L_s.\]
In fact,
\[\GW(M)\cong \GW(L_1)\oplus\ldots \oplus\GW(L_s).\]
For a field $L$, the Grothendieck-Witt ring $\GW(L)$ has the following nice presentation.
For $a\in L^\times$ let $\langle a\rangle\in \GW(L)$ be the class of the symmetric bilinear form $L\times L\rightarrow L$, defined by $(x,y)\mapsto axy$.
Then $\GW(L)$ is generated by symbols $\langle a\rangle$ for varying $a\in L^\times$ subject to the relations
\begin{enumerate}
    \item $\langle a\rangle=\langle ab^2\rangle$ for $a,b\in L^\times$,
    \item $\langle a\rangle\langle b\rangle=\langle ab\rangle$ for $a,b\in L^\times$,
    \item $\langle a\rangle+\langle b\rangle=\langle a+b\rangle+\langle ab(a+b)\rangle$ for $a,b,a+b\in L^\times$.
\end{enumerate}

More generally, for $M$ an \'etale $k$-algebra and $a\in M^\times$ one can also define $\langle a\rangle$ to be the class of the symmetric bilinear form $M\times M\rightarrow M$ defined by $(x,y)\mapsto axy$. In particular, the relation (1), (2) and (3) above also hold for $a,b,a+b$ in~$M^\times$.

\begin{definition}
The element
    $\gw{1}+\gw{-1}\in \GW(M)$ is called the \emph{hyperbolic form} and is denoted $h$.
\end{definition}
It is a nice exercise to show the following relation in $\GW(M)$ :
\begin{lemma}\label{lemma:hyperbolic}
    For $a\in M^\times$ it holds that $\langle a\rangle +\langle -a\rangle=\langle 1\rangle+\langle-1\rangle=h$ in $\GW(M)$. 
\end{lemma}

The following fact is very useful for the tropical approach to $\AA^1$-enumerative geometry, see Theorem 4.7 \cite{MPS22} for a proof. 
\begin{lemma}
    There is a canonical isomorphism $\GW(M)\cong \GW(\Puiseux{M})$.
\end{lemma}

Assume $L$ and $M$ are finite étale $k$-algebras and assume that $M$ is a finitely generated free $L$-module. Then there is a trace map $\tr_{M/L}\colon M\rightarrow L$. This induces a
\emph{trace map}
$\Tr_{M/L}$ from $\GW(M)$ to~$\GW(L)$, which sends the class of a non-degenerate symmetric bilinear form $\beta\colon V\times V\rightarrow M$ over~$M$ to the class of the non-degenerate symmetric bilinear form 
$$V\times V\overset{\beta}{\longrightarrow} M\xrightarrow{\tr_{M/L}}L.$$
We recall some facts about this trace map $\Tr_{M/L}\colon\GW(M)\longrightarrow\GW(L)$.
\begin{proposition}[Properties of the trace map] \;
\label{prop:traces}
    \begin{enumerate}
        \item For $\alpha,\beta\in \GW(M)$ it holds that $\Tr_{M/L}(\alpha+\beta)=\Tr_{M/L}(\alpha)+\Tr_{M/L}(\beta)$.
        \item Let $a\in L^\times$ and $b\in M^\times$. Then 
            \begin{equation*}
                \label{eq:traceofsomethingink}
                \Tr_{M/L}(\gw{a\cdot b})=\gw{a}\cdot \Tr_{M/L}(\gw{b}).
            \end{equation*}
        \item Assume further that $E$ is a finite étale algebra which is a free $M$-module. Then it holds that~$\Tr_{E/L}$ equals the composition of trace maps $\Tr_{M/L}\circ \Tr_{E/M}$.
        \item Assume $\operatorname{rank}_LM=m$, then 
        \begin{equation*}
            \Tr_{M/L}(\h)=m\cdot \h.
        \end{equation*}
        \item For $M=\quotient{L[x]}{(x^m-d)}$ with $d\in L^\times$ one gets
        \begin{equation*}
            \Tr_{M/L}\langle 1\rangle=\begin{cases}
                \langle m\rangle +\frac{m-1}{2}\h& m\text{ odd}\\
                \langle m\rangle+\langle md\rangle+\frac{m-2}{2}\h&m\text{ even}
            \end{cases}
        \end{equation*}
        and 
        \begin{equation*}
            \Tr_{M/L}\langle x\rangle=\begin{cases}
                \langle md\rangle +\frac{m-1}{2}\h& m\text{ odd}\\
                \frac{m}{2}\h&m\text{ even}.
            \end{cases}
        \end{equation*}
    \item Assume $M_1$ and $M_2$ are both finite étale $k$-algebras and both free $L$-modules. Set $M=M_1\otimes_LM_2$. Further let $a=a_1\otimes a_2\in M^\times$. Then
    \begin{equation*}
    \label{eq:traceoftensorproduct}
        \Tr_{M/L}(\gw{a})=\Tr_{M_1/L}(\gw{a_1})\cdot \Tr_{M_2/L}(\gw{a_2}).
    \end{equation*}

    \end{enumerate}
\end{proposition} 
\begin{proof}
    Statement (3) follows directly from the definition of the trace, statements (1), (2), and (6) follow from the properties of the algebraic trace. Statements (4) and (5) are shown in \cite[Proposition 2.13 and Lemma 2.12]{puentes2022quadratically}. 
\end{proof}

\subsection{$\GW(k)$-valued count of rational curves}\label{GWk-counts}

We start by recalling some results from \cite{LevineWelschinger,KassLevineSolomonWickelgren,KLSWOrientation} about the algebraic curve count $N_{S, \beta}^{\AA^1}(\sigma)$.
Let $S$ be a smooth del Pezzo surface over a field $k$ and $\beta$ an effective Cartier divisor on $S$. A dimension count yields that there are finitely many degree $\beta$ rational curves in $S$ through $n=-K_S\cdot \beta-1$ points in general position.

\begin{definition} \label{def:quadweightofcurve}
    Let $u\colon C\rightarrow S$ where $C$ is a nodal curve defined over the extension field $\kappa(u) \supset k$. 
     \begin{enumerate}
        \item Let $p$ be a node of the image $u(C)$ defined over the field $\kappa(p) \supset \kappa(u)$. The two branches of~$p$ are defined over and extension $\kappa(p)\subset \kappa(p)[\sqrt{a}]$ for some $a\in \kappa(p)^\times/(\kappa(p)^\times)^2$. The \emph{mass} of the node $p$ is
        \[\mass(p)\coloneqq \big\langle N_{\kappa(p)/\kappa(u)}  (a) \big\rangle\in \GW(\kappa(u)),\]
        where $N_{L/k}$ is the field norm.
        \item The \emph{quadratic weight} of $u$ is
        \[\Wel(u)\coloneqq \prod_{\text{nodes } p}\mass(p)\in \GW(\kappa(u)).\]
    \end{enumerate} 
\end{definition}
\begin{remark}
    \cref{def:quadweightofcurve} generalizes Welschinger's definition of the mass of a node and the Welschinger sign in the real case from \cite{Wel03,Wel05}:
    Recall that Welschinger defined the mass $\mass(p)$ of a real node to be $0$ when the node was hyperbolic and $1$ when it was elliptic. The Welschinger sign of a real curve is $\operatorname{Wel}(u)=\prod_{\text{real nodes }p}(-1)^{\mass(p)}$. Now it is easy to check that for a curve $u$ defined over $\RR$ we have
    \[\Wel(u)=\langle \operatorname{Wel}(u)\rangle\in \GW(\mathbb{R}).\]
\end{remark}

\begin{remark}
\label{remark:computingWel}
One can write down a formula for $\Wel(u)$ as follows:
    If the curve $u(C)$ is defined by the vanishing of $f\in \kappa(u)[x,y]$ and $p$ is a node, then 
    \[\Wel(u)=\gw{\prod_{\text{nodes }p}N_{\kappa(p)/\kappa(u)}\left(-\det \operatorname{Hessian}f(p)\right)}\in \GW \big(\kappa(u(C)) \big)\]
    where $\kappa(u)$ is the field of definition of the curve $u(C)$ and $\kappa(p)$ is the residue field of $p$.
\end{remark}

Now let $\sigma=(L_1,\ldots,L_q)$ be a $q$-tuple of field extensions $k\subset L_i\subset \overline{k}$ such that $\sum_{i=1}^q[L_i:k]=n$.
The following theorem was originally stated in \cite{KassLevineSolomonWickelgren} in a more general setting. Here, we specialize it to the case of a toric del Pezzo surface $S$ defined over the base field $k$.

\begin{theorem}[Levine, Kass-Levine-Solomon-Wickelgren]
\label{thm:KLSW}
    Assume $k$ is a perfect field of characteristic not equal to $2$ or $3$. Let $S$ be a toric del Pezzo surface defined over $k$ and let $p_1,\ldots,p_q$ in $S$ be any generically chosen points with $\kappa(p_i)=L_i$. Define
    \[N^{\A^1}_{S,\beta}(\sigma)= \sum \Tr_{\kappa(u)/k}\Wel(u)\in \GW(k),\]
    where the sum runs over all rational curves of degree $\beta$ through $p_1,\ldots,p_q$. Then $N^{\A^1}_{S, \beta}$ is well-defined, that is independent of the chosen point configuration.
\end{theorem}

\begin{definition}\label{def-quadraticallyenrichedcount}
    We work in \cref{setting}, so for us $S$ is a toric surface with Newton polygon $\Delta$ as in \cref{fig-degrees}, $n=n(\Delta)$, and $\sigma= \big(k,\ldots,k,k(\sqrt{d_1}),\ldots,k(\sqrt{d_s}) \big)$. In our situation we use the notation 
    $N^{\A^1}_\Delta(r,(d_1,\ldots,d_s)):=N^{\A^1}_{S,\beta}(\sigma)$.
\end{definition}

\begin{example}
    For $k=\mathbb{R}$ when we take the signature of $N^{\A^1}_{S,\beta}$ we obtain Welschinger's invariant. Note that in this case $L_i=\mathbb{R}$ or $L_i=\mathbb{C}$, that is we can assume that $\sigma=(\mathbb{R},\ldots,\mathbb{R},\mathbb{C},\ldots,\mathbb{C})$. Then 
    \begin{align*} 
        \sgn\left(N^{\A^1}_{S,\beta}(\sigma)\right)&=\sgn\left(\sum \Tr_{\kappa(u)/k}\Wel(u)\right)\\
        &=\sum_{\kappa(u)=\RR} \sgn \big(\Wel(u) \big)+\sum_{\kappa(u)=\CC} \sgn \big(\Tr_{\CC/\mathbb{R}}\Wel(u) \big)\\
        &=\sum_{\kappa(u)=\RR}\operatorname{Wel}(u)+0 \\
        &= W_{S,\beta}(r,s).
    \end{align*}
\end{example}

The following follows directly from \cite[Proposition 5.26]{BrugalleWickelgren}.
\begin{lemma}
\label{lm:GWPuiseux}
Let $S$ be a del Pezzo surface over $k$ and let $S_{\Puiseux{k}}$ be its base change. Further, let $\sigma$ be a tuple~$(L_1,\ldots,L_q)$ of field extension of $k$ and denote $\sigma_{\Puiseux{k}}= \big(\Puiseux{L_1},\ldots,\Puiseux{L_q} \big)$.
    Under the natural isomorphism $\GW \big(\Puiseux{k} \big) \cong \GW(k)$ the invariant $N^{\A^1}_{S_{\Puiseux{k}},\beta}(\sigma_{\Puiseux{k}})$ is mapped to $N^{\A^1}_{S,\beta}(\sigma)$.
\end{lemma}


\section{Tropical curves with double point conditions}\label{sec-tropdoublepoints}

In this section, we study the combinatorial structure of the tropical stable maps which appear in the tropical version of our enumerative problem. To this end recall \cref{setting}: we are trying to count rational plane curves of some degree $\Delta$ defined over an extension of $K$ which pass through~$r$ points defined over $K$ and $s$ pairs of conjugate points defined over $K(\sqrt{d_i})$ for $i=1\ldots,s$, with $r+2s=n(\Delta)$.
When we tropicalize pairs of conjugate points, we obtain the same point in $\RR^2$. 
Therefore, in the tropical approach to our counting problem, we cannot work with $n=n(\Delta)$ points in tropical general position. Rather, we have to consider $r+s$ points, but each of the $s$ points obtained by tropicalizing two conjugate points imposes a \enquote{double condition}. 
Following the convention in \cite{Shu06b}, we draw such points as \emph{fat points} in $\mathbb{R}^2$, whereas the $r$ points coming from $K$-rational points are drawn as \emph{thin points}.

\begin{definition}[Simple and double point conditions]\label{def-thinfat}
Fix a degree $\Delta$ and let $r+2s=n(\Delta)$. Assume the contracted ends of a rational tropical stable map $(\Gamma,f)$ of degree $\Delta$ are marked with $x_1,\ldots,x_r$, $y_1,\ldots,y_s$ and $z_1,\ldots,z_s$. Fix $r$ (\enquote{simple}) points $p_1,\ldots,p_r\in \mathbb{R}^2$ and $s$ (\enquote{double}) points $q_1,\ldots,q_s$ in $\RR^2$. We say that $(\Gamma,f)$ \emph{satisfies the simple and double point conditions} if $f(x_i)=p_i$ for all $i=1,\ldots,r$ and $f(y_i)=f(z_i)=q_i$ for all $i=1,\ldots,s$.
\end{definition}

\begin{remark}\label{rem-tropicalizationofstablemaps}
By commutative square \eqref{eq:comm_square}, counting tropical curves through $r$ simple and $s$ double points is the correct thing to do. Indeed, fix $r$ points over $K$ and $s$ pairs of conjugate points over $K(\sqrt{d_i})$ (with $n=r+2s$) in $T^n$. To solve the algebraic enumerative problem, one has to study the preimage of this configuration under $\ev$ in the moduli space of log stable maps. If we tropicalize the log stable maps in this preimage we obtain precisely the set of tropical stable maps passing through $r$ simple and $s$ double points. 
\end{remark}

\subsection{The moduli space of rational tropical stable maps satisfying double point conditions}
\begin{definition}
    Let $MC_{0,s}(\mathbb{R}^2,\Delta)$ denote the pull-back of zero under $$(\ev_{y_1}-\ev_{z_1})\times\ldots\times(\ev_{y_s}-\ev_{z_s} ) : M_{0, r+2s}(\RR^2, \Delta) \longrightarrow (\RR^2)^s.$$ 
    This is the space of tropical stable maps which can meet $s$ double points.
\end{definition}

The following statement follows from basics in tropical intersection theory (see Proposition 4.7 and Definition 3.8 in \cite{AR07}):

\begin{lemma}\label{lem-MC}
The moduli space $MC_{0,s}(\mathbb{R}^2,\Delta)$ of rational tropical stable maps satisfying $s$ double point conditions is a polyhedral subcomplex of $M_{0,n}(\mathbb{R}^2,\Delta)$ of dimension $2\cdot(n(\Delta))-2s$. 
\end{lemma}

\begin{lemma}\label{lem-evinverse}
In the setting of Definition \ref{def-thinfat}, the set of tropical stable maps satisfying $r$ simple and $s$ double point conditions $(p_i,q_j)$ equals the inverse image of $(p_i,q_j)$ under $\ev_{x_1}\times\ldots\times\ev_{x_r}\times\ev_{y_1}\times\ldots\times\ev_{y_s}$ restricted to $MC_{0,s}(\mathbb{R}^2,\Delta)$.
\end{lemma}
\begin{proof}
Notice first that $\ev_{x_1}\times\ldots\times\ev_{x_r}\times\ev_{y_1}\times\ldots\times\ev_{y_s}$ is a morphism starting from a polyhedral complex of dimension $2\cdot(n(\Delta))-2s$ to a polyhedral complex of dimension $2(r+s)$. But since $r+2s=n(\Delta)$, these are polyhedral complexes of the same dimension. Thus, generically we have a finite fiber. The tropical stable maps in the inverse image of $(p_i,q_j)$ satisfy $\ev_{x_i}(\Gamma,f)=f(x_i)=p_i$ and $\ev_{y_i}(\Gamma,f)=f(y_i)=q_i$. Furthermore, since we restrict the map to $MC_{0,s}(\mathbb{R}^2,\Delta)$, $\ev_{y_i}(\Gamma,f)=ev_{z_i}(\Gamma,f)=f(z_i)=q_i$. Hence any tropical stable map in the inverse image satisfies the simple and double point conditions  $(p_i,q_j)$. The converse is clear.
\end{proof}

\begin{definition}
A set of $r+s$ simple and double points in $\mathbb{R}^2$ is in \emph{general position}, if all tropical stable maps that pass through them are contained in the interior of top-dimensional cells of $MC_{0,s}(\mathbb{R}^2,\Delta)$.
\end{definition}

\begin{lemma}\label{lem-no-deform}
A tropical stable map $(\Gamma,f)$ of degree $\Delta$ passing through $r+s$ simple and double points in general position with $r+2s=n(\Delta)$ is locally fixed by the point conditions. That is, there is no deformation of $(\Gamma,f)$ within $MC_{0,s}(\mathbb{R}^2,\Delta)$ (i.e.\ a tropical stable map of the same combinatorial type which is close in the corresponding orthant) that still passes through the points. 
\end{lemma}
\begin{proof}
This follows from a dimension count: by Lemma \ref{lem-evinverse}, $(\Gamma,f)$ is in the inverse image of
 $(p_i,q_j)$ under $\ev_{x_1}\times\ldots\times\ev_{x_r}\times\ev_{y_1}\times\ldots\times\ev_{y_s}$, which is a morphism of polyhedral complexes of the same dimension. In particular, it is locally a linear map, and  for point conditions in general position, we can have at most one preimage in a cell. 
\end{proof}

\begin{lemma}\label{lem-notoverfixed}
Let $(\Gamma,f)$ be a tropical stable map of degree $\Delta$ passing through  $r+s$ simple and double points in general position with $r+2s=n(\Delta)$. Then $\Gamma$ minus the (closures of the) marked ends cannot contain a bounded component.
\end{lemma}
\begin{proof}
Assume it did, and the bounded component was adjacent to $m$ marked points, i.e.\ it has $m$ one-valent vertices (adjacent to leaves) and all other vertices are of valence at least three. We first argue that it has at least one vertex of valence $l$ which is adjacent to $l-1$ leaves. Let $t$ denote the number of vertices of valence at least $3$. If there was no such vertex, we had
$ m\leq \sum_{v} (\val(v)-2)$, where the sum goes over all vertices $v$ of valence at least $3$. But $t=m-2-\sum_v (\val(v)-3)= m-2 - \sum_v \val(v) +3k$, hence $2k= \sum_v \val(v) -m+2 $. Inserting this in the inequality above, we obtain
$$ m\leq \sum_{v} \big(\val(v)-2 \big) =  \sum_v \val(v) -  \sum_v \val(v) +m-2 = m-2$$
which is a contradiction.

Thus there is a vertex $v_0$ of valence $l$ which is adjacent to $l-1$ leaves. 
For $v_0$, the lines into which its leaf edges are mapped via $f$ are fixed by the point conditions for the marked point adjacent in $\Gamma$. Accordingly, also the position of $v_0$ and with it, of its adjacent bounded edge is fixed. Recursively, we find a leaf whose position is already fixed by the other point conditions; its position is overdetermined. We can thus drop the marked point adjacent to this leaf and its corresponding point condition.
Dropping the marked point reduces the dimension of the moduli space by one, but the target space of the evaluation map (evaluating only the other points) is reduced by dimension $2$. Accordingly, there is a one-dimensional fiber of the remaining point conditions, and every tropical stable map in this fiber automatically also passes through the dropped point condition, as it passes through the points required  for the other marked points adjacent to the bounded component. This is a contradiction to Lemma \ref{lem-no-deform}.
\end{proof}

\subsection{Local building blocks of tropical stable maps satisfying double point conditions}

We can deform a double point into two simple points. As generic points are dense, we can achieve such a deformation in such a way that we obtain $n$ points in tropical general position. If $(\Gamma,f)$ is a rational tropical stable map satisfying the simple and double point conditions, we can deform it to obtain a simple tropical stable map through $n$ points. To understand the local building blocks of tropical stable maps that satisfy simple and double point conditions, we can reverse the process: We start with a simple tropical stable map through $n$ points and then move pairs of points close together until they get identified.

With the following lemma, we classify all combinatorial possibilities how a tropical stable map around a double point condition can look like.

\begin{figure}[t]
\begin{center}
    \tikzset{every picture/.style={line width=0.75pt}} 
\setlength\tabcolsep{0.5cm}
\begin{tabular}{c c c c c}
    \begin{tikzpicture}
        \path[draw] (0,0) --+ (90:1);
        \path[draw] (0,0) --+ (210:1);
        \path[draw] (0,0) --+ (-30:1);
        \fatpoint{0, 0}
        \draw (0, -1) node {\crtcrossreflabel{(A)}[fatPointOnVertex]};
    \end{tikzpicture}
    &
    \begin{tikzpicture}
        \path[draw] (30:-1) -- (30:1);
        \path[draw] (60:-1) -- (60:1);
        \fatpoint{0, 0}
        \draw (0, -1) node {\crtcrossreflabel{(B)}[fatPointOnParallelogram]};
    \end{tikzpicture}
    &
    \begin{tikzpicture}
        \path[draw] (0,0) -- (135:0.5);
        \path[draw] (0,0) -- (-135:0.5);
        \path[draw, double] (0,0) -- (1,0);
        \fatpoint{0.5, 0}
        \draw (0.5, -1) node {\crtcrossreflabel{(C)}[fourValentVertex]};
    \end{tikzpicture}
    &
    \begin{tikzpicture}
        \path[draw, double] (30:-1) -- (30:0.7);
        \path[draw] (60:-1) -- (60:0.7);
        \fatpoint{30:-0.6}
        \thinpoint{60:-0.6}
        \draw (0, -1) node {\crtcrossreflabel{(D)}[parallelogramWithOneDoubleEdge]};
    \end{tikzpicture} 
    &
    \begin{tikzpicture}
        \path[draw, double] (30:-1) -- (30:0.7);
        \path[draw, double] (60:-1) -- (60:0.7);
        \fatpoint{30:-0.6}
        \fatpoint{60:-0.6}
        \draw (0, -1) node {\crtcrossreflabel{(E)}[parallelogramWithTwoDoubleEdges]};
    \end{tikzpicture}
    \\[1cm]
    \begin{tikzpicture}
        \path[draw] (0,\doubleedgesep) --+ (135:0.5);
        \path[draw] (0,\doubleedgesep) --+ (-90:0.5);
        \path[draw] (0,\doubleedgesep) -- (1,\doubleedgesep);
        \path[draw] (-0.5, -\doubleedgesep) -- (1, -\doubleedgesep);
        \fatpoint{0.5, 0}
        \draw (0.5, -1) node {\crtcrossreflabel{(F)}[triangleWithMergedEdge]};
    \end{tikzpicture}
    &
    \begin{tikzpicture}
        \path[draw] (0,\doubleedgesep) --+ (135:0.5);
        \path[draw, double] (0,\doubleedgesep) --+ (-90:0.7);
        \path[draw] (0,\doubleedgesep) -- (1,\doubleedgesep);
        \path[draw] (-0.5, -\doubleedgesep) -- (1, -\doubleedgesep);
        \fatpoint{0.5, 0}
        \fatpoint{0, -0.4}
        \draw (0.5, -1) node {\crtcrossreflabel{(G)}[fourValentVertexWithMergedEdge]};
    \end{tikzpicture}
    &
    \begin{tikzpicture}
        \path[draw] (\doubleedgesep,\doubleedgesep) --+ (135:0.5);
        \path[draw] (\doubleedgesep,\doubleedgesep) -- (\doubleedgesep, -0.7);
        \path[draw] (\doubleedgesep,\doubleedgesep) -- (1,\doubleedgesep);
        \path[draw] (-0.5, -\doubleedgesep) -- (1, -\doubleedgesep);
        \path[draw] (-\doubleedgesep, -0.7) -- (-\doubleedgesep, 0.5);
        \fatpoint{0.5, 0}
        \fatpoint{0, -0.4}
        \draw (0.5, -1) node {\crtcrossreflabel{(H)}[triangleWithTwoMergedEdges]};
    \end{tikzpicture}
    &
    \begin{tikzpicture}
        \path[draw, double] (0,0) -- (1, 0);
        \path[draw, double] (0,-0.7) -- (0, 0) --+ (135:0.5);
        \fatpoint{0.5, 0}
        \fatpoint{0, -0.4}
        \draw (0.5, -1) node {\crtcrossreflabel{(I)}[allDoubleVertex]};
    \end{tikzpicture}
    &
    \begin{tikzpicture}
        \path[draw, double] (0,0) -- (1, 0);
        \path[draw, double] (0,-0.7) -- (0, 0);
        \path[draw] (0,0) --+ (115:0.5);
        \path[draw] (0,0) --+ (165:0.5);
        \fatpoint{0.5, 0}
        \fatpoint{0, -0.4}
        \draw (0.5, -1) node {\crtcrossreflabel{(J)}[mergedTriangles]};
    \end{tikzpicture}
\end{tabular}

 
\end{center}
\caption{Possibilities for double point conditions of a tropical stable map. We depict the image $f(\Gamma)$ but indicate the parametrization by the abstract graph. We call these local building blocks the \emph{vertex types}.}\label{fig-casesfatpoint}
\end{figure}
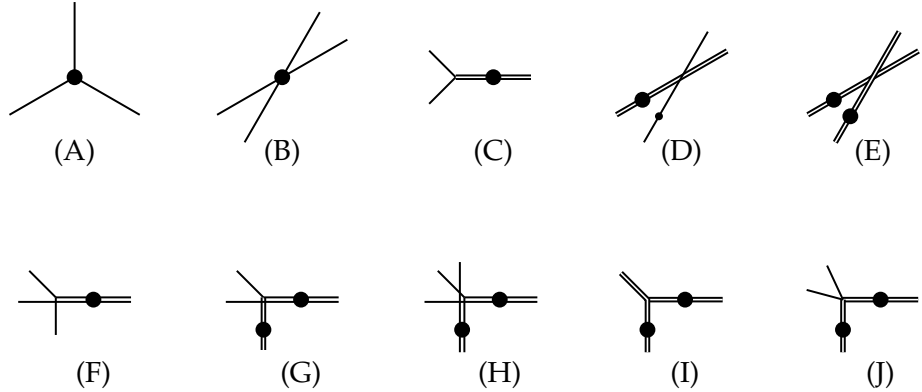

\begin{lemma}\label{lem-casesfatpoint}
Let $(\Gamma,f)$ be a rational tropical stable map satisfying $r$ simple and $s$ double point conditions. Let $q$ be a double point. Then, locally around $q$,  the image $f(\Gamma)$ (and its parametrization) can look like one of the cases described below and depicted in Figure \ref{fig-casesfatpoint}.
\begin{enumerate}
\item The double point can sit on a $3$-valent vertex of the image $f(\Gamma)$ (Case~\fatPointOnVertex{}). 
\item The double point can sit on a crossing of two edges of $f(\Gamma)$ (Case~\fatPointOnParallelogram{}). 
\item The double point can sit on a ''double edge'', i.e.\ an edge of the image $f(\Gamma)$ to which two edges of $\Gamma$ of the same primitive direction are mapped.
\end{enumerate}
The last case encompasses several subcases (Cases~\fourValentVertex{} -~\mergedTriangles{}) according to the possible end vertices of the two parallel edges as depicted in Figure \ref{fig-casesfatpoint}. The picture shows the image $f(\Gamma)$, but indicates the parametrization by the abstract tropical curve. Double edges as in cases~\fourValentVertex{} -~\mergedTriangles{} can be ends (if their direction is in $\Delta$) or be  adjacent to another of these cases. Then,  two pieces of local building blocks are glued at a double point of each. Case~\fourValentVertex{} cannot be combined with itself nor with the vertical double edge in case~\triangleWithMergedEdge{}.

In cases~\fourValentVertexWithMergedEdge{} -~\mergedTriangles{}, the vertical double edge does  not need to pass through a double point, but its position has to be fixed by double point conditions.
\end{lemma}

\begin{proof}
The proof works by resolving the double point via deformation into two points, then study the possibilities for tropical stable maps meeting these simple point conditions, and finally specializing back to a double point.
Case~\fatPointOnVertex{} arises if two adjacent marked points (i.e.\ two points on edges sharing a $3$-valent vertex) merge into a $5$-valent vertex in $\Gamma$, shrinking the bounded edges that separated them before to length $0$.
Case~\fatPointOnParallelogram{} arises if two marked points adjacent to edges of different direction are forced to map to the same point.
The remaining cases arise if edges of the same primitive direction are forced to map to overlapping segments of the same image line. We call two such edges double edges. Case~\fourValentVertex{} arises if the double edges are adjacent to a shrinking edge.
Notice that \fourValentVertex{} implies the existence of a $4$-valent vertex.
Cases \parallelogramWithOneDoubleEdge{}, \parallelogramWithTwoDoubleEdges{} arise if two such double edges cross each other, resp.\ if a double edge crosses a single edge.
Case~\triangleWithMergedEdge{} arises if one of the  double edges passes through the adjacent vertex of the other. In the dual subdivision, we obtain a trapezoid. As one of the two edges could also pass through an adjacent vertex which is of the type in \fourValentVertex{}, we also obtain case \fourValentVertexWithMergedEdge{} as another possibility.
One of the two edges could also pass through a vertex through which already another pair of double edges passes, leading to case \triangleWithTwoMergedEdges{}. 
Finally, the two edges could end at a double vertex, i.e.\ a vertex adjacent to another double edge forced by a double point. If all pairs of double edges through double points are  of the same direction (up to multiple), also the third adjacent edges for each of the vertices, respectively, are of the same direction by balancing. This is case \allDoubleVertex{}. Case \mergedTriangles{} appears if the edges are of the same primitive directions but of different weights.
In cases \fourValentVertexWithMergedEdge{} - \mergedTriangles{}, the vertical  double edges do not need to pass through a double point directly, it could also be that this double edge is fixed by other double points like the third double edge of case \allDoubleVertex{}.

All double edges can be ends if their direction is in $\Delta$. Else, any of the pictures in cases \fourValentVertex{} - \mergedTriangles{} can be combined with any suitable other such case.
Case \fourValentVertex{} cannot be combined with itself nor with the vertical double edge in case \fourValentVertexWithMergedEdge{}, as this would produce a cycle in the abstract curve~$\Gamma$ which is not possible since we consider rational tropical curves, i.e.\ $\Gamma$ is a tree. Because of \cref{lem-notoverfixed}, we deduce that two pieces of local building blocks are glued at a double point of each, resp.\ at an edge which is fixed by other double points as described above.
\end{proof}

In the following, when we speak about edges and vertices of a tropical stable map $(\Gamma,f)$, by abuse of notation we refer to edges and vertices of the image $f(\Gamma)$. The graph structure of the underlying abstract tropical curve is controlled by the local building blocks. We call these the \emph{vertex types}. 

\begin{remark}
Notice that the list of local building blocks from Lemma \ref{lem-casesfatpoint} and Figure \ref{fig-casesfatpoint} confirms our dimension statement from Lemma \ref{lem-MC}, meaning that each double point lowers the dimension of the moduli space by 2. 
\begin{itemize}
    \item In case \fatPointOnVertex{}, we shrink two bounded edges, producing a $5$-valent vertex in $\Gamma$ with $f$ contracting two ends and producing a $3$-valent image, which reduces the dimension by two.

    \item In case \fatPointOnParallelogram{}, the fact that two marked ends have to be precisely at the point of intersection of two edges also reduces the dimension by two. 
    
    \item In case \fourValentVertex{} and for the lower point of case \fourValentVertexWithMergedEdge{}, we shrink a bounded edge, producing a $4$-valent vertex. Furthermore, we require two marked points to have the same distance from the $4$-valent vertex. Altogether, we lower the dimension by two as well.
    
    \item The double points in cases \parallelogramWithOneDoubleEdge{} and \parallelogramWithTwoDoubleEdges{} must be adjacent to a vertex in another case, which yields the dimension reduction.
    
    \item In case \triangleWithMergedEdge{}, for the other double point in case \fourValentVertexWithMergedEdge{} and for both double points in cases \triangleWithTwoMergedEdges{}, \allDoubleVertex{} and \mergedTriangles{}, the fact that two edges of the same primitive direction have to map into the same line yields an equation on the lengths of the bounded edges that form the path from one of the edges to the other. Furthermore, the equality of the position of two marked ends yields another condition. Altogether, we reduce the dimension by two again.    
\end{itemize}
As we consider rational tropical stable maps, multiple such conditions will be independent.
\end{remark}

\subsection{Structure of connected components of double edges}

Let $(\Gamma,f)$ be a rational tropical stable map satisfying simple and double point conditions. In the image $f(\Gamma)$, consider a connected component of double edges, i.e.\ edges of $f(\Gamma)$ to which two edges of $\Gamma$ of the same primitive direction are mapped. An end point of such a connected component can arise from any of the cases \fourValentVertex{} - \mergedTriangles{} depicted in Figure \ref{fig-casesfatpoint}; cases \fourValentVertexWithMergedEdge{}, \triangleWithTwoMergedEdges{} and \mergedTriangles{} yield two connected components of double edges in~$f(\Gamma)$ as we treat these end vertices as separating. 

\begin{lemma}\label{lem-doubleedges3valent}
Every vertex in a connected component of double edges is $3$-valent.
\end{lemma}

\begin{proof}
This follows as case \allDoubleVertex{} in Figure \ref{fig-casesfatpoint} is the only case which yields vertices in a connected component of~$\Gamma$ mapping to double edges, and all these vertices are then $3$-valent.
\end{proof}

We treat each end point of such a connected component as an end, i.e.\ we parametrize each connected component of double edges in $f(\Gamma)$ by a $3$-valent abstract tropical curve, possibly prolonging end points arising from cases \fourValentVertex{} - \triangleWithTwoMergedEdges{}, \mergedTriangles{} in Figure \ref{fig-casesfatpoint} to ends.

We orient the double edges of $f(\Gamma)$ such that they point away from the double points, and such that if two inward-oriented edges meet at a vertex, the third adjacent edge is oriented outwards. In this way we can inductively orient all edges of a connected component of double edges. 

\begin{figure}[t]
    \begin{center}
    	\tikzset{every picture/.style={line width=0.75pt}} 

\begin{tikzpicture}[x=0.75pt,y=0.75pt,yscale=-1,xscale=1]

\draw    (349.1,332.13) -- (338.5,321.6) ;
\draw    (349.39,353.76) -- (349.25,331.23) ;
\draw    (351.07,354.43) -- (351.13,333.28) ;
\draw    (387.1,331.6) -- (348.5,331.6) ;
\draw  [fill={rgb, 255:red, 0; green, 0; blue, 0 }  ,fill opacity=1 ] (369.42,332.73) .. controls (369.42,330.99) and (370.82,329.58) .. (372.56,329.58) .. controls (374.29,329.58) and (375.7,330.99) .. (375.7,332.73) .. controls (375.7,334.46) and (374.29,335.87) .. (372.56,335.87) .. controls (370.82,335.87) and (369.42,334.46) .. (369.42,332.73) -- cycle ;
\draw    (385.95,334.2) -- (350.98,334.18) ;
\draw  [fill={rgb, 255:red, 0; green, 0; blue, 0 }  ,fill opacity=1 ] (347.07,345.05) .. controls (347.07,343.31) and (348.47,341.91) .. (350.21,341.91) .. controls (351.94,341.91) and (353.35,343.31) .. (353.35,345.05) .. controls (353.35,346.78) and (351.94,348.19) .. (350.21,348.19) .. controls (348.47,348.19) and (347.07,346.78) .. (347.07,345.05) -- cycle ;
\draw    (351.94,334.07) -- (350.41,332.7) ;
\draw    (349.73,330.59) -- (339.78,320.7) ;
\draw    (280,300) -- (300,320) ;
\draw    (300,320) -- (280,340) ;
\draw  [fill={rgb, 255:red, 0; green, 0; blue, 0 }  ,fill opacity=1 ] (288.7,310) .. controls (288.7,309.28) and (289.28,308.7) .. (290,308.7) .. controls (290.72,308.7) and (291.3,309.28) .. (291.3,310) .. controls (291.3,310.72) and (290.72,311.3) .. (290,311.3) .. controls (289.28,311.3) and (288.7,310.72) .. (288.7,310) -- cycle ;
\draw  [fill={rgb, 255:red, 0; green, 0; blue, 0 }  ,fill opacity=1 ] (288.7,330) .. controls (288.7,329.28) and (289.28,328.7) .. (290,328.7) .. controls (290.72,328.7) and (291.3,329.28) .. (291.3,330) .. controls (291.3,330.72) and (290.72,331.3) .. (290,331.3) .. controls (289.28,331.3) and (288.7,330.72) .. (288.7,330) -- cycle ;
\draw    (298.64,321.47) -- (339.36,322.01) ;
\draw    (300,320) -- (336.64,320.01) ;
\draw    (340.33,292.17) -- (340.33,321.08) ;
\draw    (338.64,292.01) -- (338.5,321.6) ;
\draw    (269.6,339.87) -- (280,340) ;
\draw    (280,340) -- (280.07,350.27) ;
\draw  [fill={rgb, 255:red, 0; green, 0; blue, 0 }  ,fill opacity=1 ] (272.3,339.73) .. controls (272.3,339.02) and (272.88,338.43) .. (273.6,338.43) .. controls (274.32,338.43) and (274.9,339.02) .. (274.9,339.73) .. controls (274.9,340.45) and (274.32,341.03) .. (273.6,341.03) .. controls (272.88,341.03) and (272.3,340.45) .. (272.3,339.73) -- cycle ;
\draw  [fill={rgb, 255:red, 0; green, 0; blue, 0 }  ,fill opacity=1 ] (278.7,346.17) .. controls (278.7,345.45) and (279.28,344.87) .. (280,344.87) .. controls (280.72,344.87) and (281.3,345.45) .. (281.3,346.17) .. controls (281.3,346.88) and (280.72,347.47) .. (280,347.47) .. controls (279.28,347.47) and (278.7,346.88) .. (278.7,346.17) -- cycle ;
\draw    (540.44,330.8) -- (531.73,322.74) ;
\draw    (540.72,352.42) -- (540.59,329.89) ;
\draw    (542.4,353.09) -- (542.47,331.95) ;
\draw    (578.44,330.27) -- (539.84,330.27) ;
\draw  [fill={rgb, 255:red, 0; green, 0; blue, 0 }  ,fill opacity=1 ] (560.75,331.39) .. controls (560.75,329.66) and (562.16,328.25) .. (563.89,328.25) .. controls (565.63,328.25) and (567.03,329.66) .. (567.03,331.39) .. controls (567.03,333.13) and (565.63,334.53) .. (563.89,334.53) .. controls (562.16,334.53) and (560.75,333.13) .. (560.75,331.39) -- cycle ;
\draw    (577.28,332.87) -- (542.32,332.85) ;
\draw  [fill={rgb, 255:red, 0; green, 0; blue, 0 }  ,fill opacity=1 ] (538.4,343.71) .. controls (538.4,341.98) and (539.81,340.57) .. (541.54,340.57) .. controls (543.28,340.57) and (544.68,341.98) .. (544.68,343.71) .. controls (544.68,345.45) and (543.28,346.86) .. (541.54,346.86) .. controls (539.81,346.86) and (538.4,345.45) .. (538.4,343.71) -- cycle ;
\draw    (543.27,332.74) -- (541.74,331.36) ;
\draw    (540.79,328.76) -- (532.96,321.68) ;
\draw    (471.33,298.67) -- (487.18,314.56) ;
\draw    (487.18,322.56) -- (471.33,338.67) ;
\draw  [fill={rgb, 255:red, 0; green, 0; blue, 0 }  ,fill opacity=1 ] (480.03,308.67) .. controls (480.03,307.95) and (480.62,307.37) .. (481.33,307.37) .. controls (482.05,307.37) and (482.63,307.95) .. (482.63,308.67) .. controls (482.63,309.38) and (482.05,309.97) .. (481.33,309.97) .. controls (480.62,309.97) and (480.03,309.38) .. (480.03,308.67) -- cycle ;
\draw  [fill={rgb, 255:red, 0; green, 0; blue, 0 }  ,fill opacity=1 ] (480.03,328.67) .. controls (480.03,327.95) and (480.62,327.37) .. (481.33,327.37) .. controls (482.05,327.37) and (482.63,327.95) .. (482.63,328.67) .. controls (482.63,329.38) and (482.05,329.97) .. (481.33,329.97) .. controls (480.62,329.97) and (480.03,329.38) .. (480.03,328.67) -- cycle ;
\draw    (487.18,322.56) -- (531.73,322.74) ;
\draw    (525.96,291.01) -- (526.04,314.93) ;
\draw    (531.79,292.13) -- (531.73,322.74) ;
\draw    (460.93,338.53) -- (471.33,338.67) ;
\draw    (471.33,338.67) -- (471.4,348.93) ;
\draw  [fill={rgb, 255:red, 0; green, 0; blue, 0 }  ,fill opacity=1 ] (463.63,338.4) .. controls (463.63,337.68) and (464.22,337.1) .. (464.93,337.1) .. controls (465.65,337.1) and (466.23,337.68) .. (466.23,338.4) .. controls (466.23,339.12) and (465.65,339.7) .. (464.93,339.7) .. controls (464.22,339.7) and (463.63,339.12) .. (463.63,338.4) -- cycle ;
\draw  [fill={rgb, 255:red, 0; green, 0; blue, 0 }  ,fill opacity=1 ] (470.03,344.83) .. controls (470.03,344.12) and (470.62,343.53) .. (471.33,343.53) .. controls (472.05,343.53) and (472.63,344.12) .. (472.63,344.83) .. controls (472.63,345.55) and (472.05,346.13) .. (471.33,346.13) .. controls (470.62,346.13) and (470.03,345.55) .. (470.03,344.83) -- cycle ;
\draw    (487.18,314.56) -- (487.18,322.56) ;
\draw    (487.18,314.56) -- (526.04,314.93) ;
\draw    (530.46,319.43) -- (526.04,314.93) ;

\end{tikzpicture}
    \end{center}
    
    \caption{If an end of a component of double edges was not fixed by the double point conditions, we could deform the double edges into two edges.}\label{fig-deformnonfixeddoubleedge}
\end{figure}
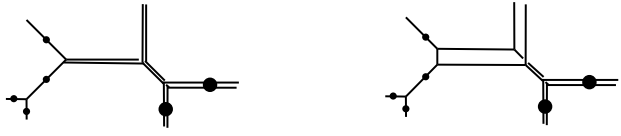

\begin{lemma}\label{lem-pointsondoubleparts}
Each connected component of double edges as above with $e$ ends passes through $e-1$ double points. Its $e$ ends are all outward oriented.
\end{lemma}
\begin{proof}
This proof follows the ideas of \cite{Mi03}, Lemma 4.20.
Assume an end of a component of double edges was not oriented outwards. Then its position is not fixed by the double point conditions. We can accordingly deform $(\Gamma,f)$ such that the two edges mapping to the same line do not map to the same line anymore, see Figure \ref{fig-deformnonfixeddoubleedge}. This deformation would still pass through the point conditions, which contradicts \cref{lem-no-deform}. Thus every end of a connected component has to be oriented outwards. We now decompose the connected component further into connected components by removing the double points. Removing one double point increases the number of components by one, so that we obtain $l+1$ components if we have $l$ double points. By  \cref{lem-notoverfixed}, no such component can be bounded. For that reason, each must contain at least one end. If one component had more than one end, it could be deformed again. Thus, every component must contain precisely one end and so we have $e-1$ double points.
\end{proof}

\begin{lemma} \label{lem-structuredoublepart}
    Any connected component of double edges must end at a vertex of type \fourValentVertex{} or \triangleWithMergedEdge{} --- \mergedTriangles{}. 
\end{lemma}
\begin{proof}
    The only way how a connected component of double edges could not end is when all of $f(\Gamma)$ consisted of double edges. But then $\Gamma$ would have two connected components, which were both mapped identically. This would only be possible if both components had weight $1$, as we require ends to have weight $1$. This required all end directions to come in pairs, i.e.\ all lattice lengths of the sides of the polygon defining the degree (see Figure \ref{fig-degrees}) would be even. By a dimension count, we obtained an odd number of point conditions if all points were thin. This implies that we cannot have only double points (which arise from two simple points), and so we cannot have only double edges, yielding a contradiction. From the classification of vertex types, we can conclude that a connected component of double edges must end at a vertex of type \fourValentVertex{} or \triangleWithMergedEdge{} --- \mergedTriangles{}.
\end{proof}

The following particular type of connected component of double edges plays an important role in the later sections of this article. 

\begin{figure}[t]
    \centering
    \begin{tikzpicture}
		\path[draw, double] (-2, 4) --++ (1, 0) --++ (1, 1);
		\fatpoint{-1.5, 4}
        \draw[thick, gray] (-1, 4) circle (4pt);
		
		\path[draw, double] (-2, 2) --++ (2, 0) --++ (1, 1);
		\fatpoint{-0.5, 2}
        \draw[thick, gray] (0, 2) circle (4pt);
		
		\path[draw, double] (-1, 4) --++ (0, -3);
		\fatpoint{-1, 3.5}
		
		\path[draw, double] (0, 2) -- ++ (0, -2);
		\fatpoint{0, 1}
		
		\path[draw, double] (-2, 1) --++ (1, 0) --++ (1, -1) --++ (0.5, -1) --++ (1, 0) -- ++ (1,1);
		\fatpoint{0.25, -0.5}
		
		\path[draw, double] (0.5, -1) --++ (0, -1); 
		\fatpoint{0.5, -1.5}
		\draw (0.3, -1.5) node[anchor = east] {\edgeweight{$(2, 2)$, \quadextension{$d_1$}}};
		
		\path[draw, double] (1.5, -1) --++ (0, -2.5);
		\fatpoint{1.5, -2.5}
		
		\path[draw] (0, -2) --++ (0.5, 0) --++ (0.5, -1);

		\draw (0, -3.5) node[anchor = north] {$\cT_3$};
        \draw (0, -4) node[anchor = north] {$t = 7$, \quad $m_\circ = 2 + 1$};
        \draw (0, -4.4) node[anchor = north, align = center] {$\mult^{\AA^1}(\cT_3) = \Big( 2\big(\gw{1} + \gw{-d_1} \big) + 6\h \Big) \cdot \gw{2^6}$ \\ $\cdot \Big( \sum_{\substack{I \subseteq \{1, \ldots 7\}  \\ |I| \text{ odd}}} \gw{\prod_{i \in I} d_i}\Big)$};
		
		\begin{scope}[shift={(-8, 3)} = ]
			\path[draw, double] (0, 1) -- ++ (0, -2);
			\fatpoint{0, -0}
			\path[draw] (-0.5, -1) --++ (0.5,0) --++ (0.25, -0.5);
			
			\draw (0, -1.5) node[anchor = north] {$\cT_1$};
            \draw (0, -2) node[anchor = north] {$t = 1$, \quad $m_\circ = 1 + 1$};
            \draw (0, -2.5) node[anchor = north] {$\mult^{\AA^1}(\cT_1) = \gw{1}$};
		\end{scope}
	
	\begin{scope}[shift={(-8, -2)}]
		\path[draw,double] (-1, 0) --++ (1, 0) --++ (1,1);
		\fatpoint{-0.5, 0}
		\path[draw, double] (0, 0) -- ++ (0, -1);
		\fatpoint{0, -0.5}
		\path[draw] (-0.5, -1) --++ (0.5,0) --++ (0.25, -0.5);
		
		\draw (0, -1.5) node[anchor = north] {$\cT_2$};
        \draw (0, -2) node[anchor = north] {$t = 2$, \quad $m_\circ = 1 + 0$};
        \draw (0, -2.5) node[anchor = north] {$\mult^{\AA^1}(\cT_2) = \gw{2d_1} + \gw{2 d_2}$};
	\end{scope}
		
	\end{tikzpicture}
    \caption{Examples for twin trees and their multiplicities as in \cref{def-multtwintree}. The color coding from \cref{not:colors} applies. In the picture on the right, the encircled vertices are possible choices for $v_0$ in the proof of \cref{lem-algebratwintree}.} 
    \label{fig-twintree}
\end{figure}
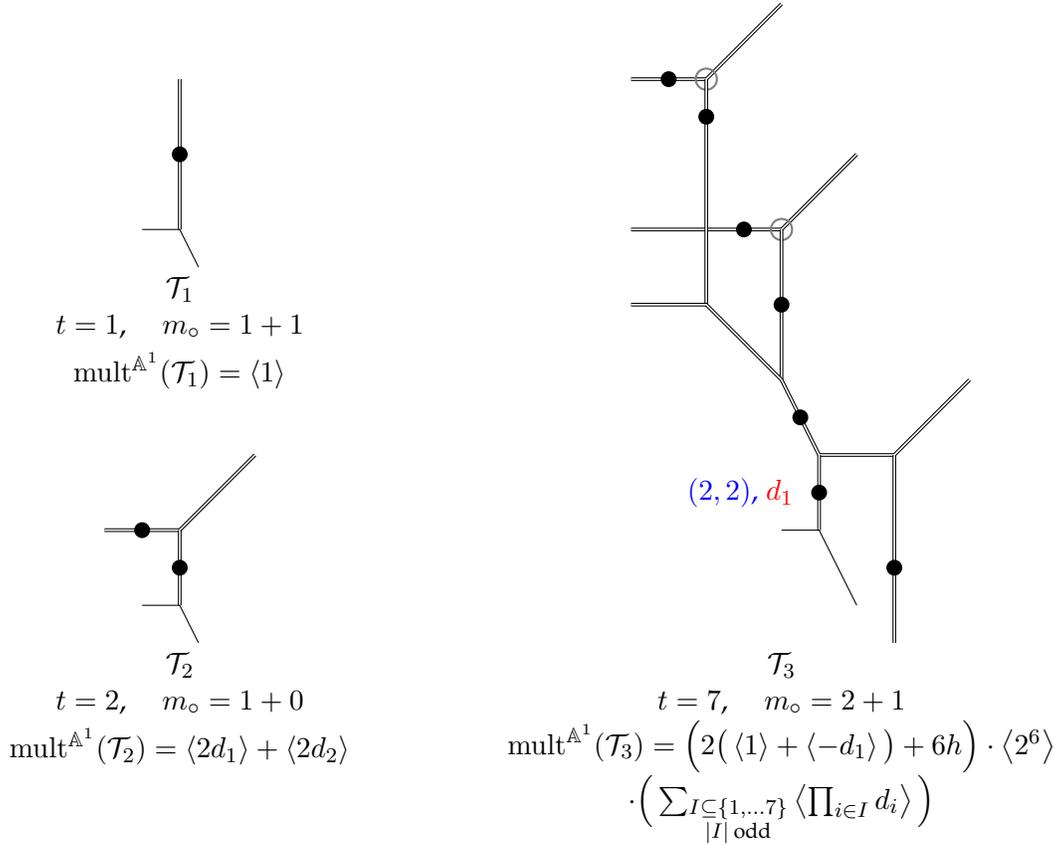

\begin{definition}\label{def-twintree}
    A connected component of double edges which consists of two identical parts, necessarily ending at a vertex of type \fourValentVertex{} (resp.\ the \fourValentVertex{}-part of type \fourValentVertexWithMergedEdge{}), is called a \emph{twin tree}. See Figure \ref{fig-twintree} for examples. A double edge which is part of a twin tree is called \emph{twin edge}.
\end{definition}

\subsection{Vertically stretched point configurations}

In tropical enumerative geometry we count tropical curves passing through a collection of points in general position. As the number of curves does not depend on the point configuration, we may choose the point conditions to be in a specific configuration (which is still generic) with the goal of reducing the combinatorial complexity of the tropical objects we are counting. This can be achieved by restricting to \emph{vertically stretched} point conditions as was introduced by \cite{BM08}.

\begin{definition}
    A set of $r$ simple and $s$ double points is in \emph{vertically stretched position} if their $x$-coordinates are contained in a small interval, while the distances between their $y$-coordinates is large. 
\end{definition}

\begin{definition}
    Let $(\Gamma,f)$ be a tropical stable map satisfying vertically stretched point conditions. 
    An \emph{elevator edge} of $(\Gamma,f)$ is an edge of $\Gamma$ of primitive direction $(0,\pm 1)$.
    If an elevator edge is adjacent to a marked point followed by an adjacent elevator edge, we call the union of these two edges and the marked point an \emph{elevator}.
\end{definition}

\begin{lemma}\label{lem-floordecomposed}
     Let $(\Gamma,f)$ be a rational tropical stable map to a smooth del Pezzo surface with degree $\Delta$ as in Definition \ref{def-deg}
     satisfying $r$ simple and $s$ double vertically stretched point positions such that $r+2s=n(\Delta)$.
     Then every connected component of $\Gamma$ minus its elevator edges which is not just a marked point (called a \emph{floor}) contains 2 non-vertical ends (either horizontal or diagonal). 
     Each edge which belongs to a floor has direction $(\pm 1,y)$ for some $y$. In particular, all edges in floors have weight 1.
     
     Each floor contains a marked point. A marked point meeting a double point condition can fix two floors for which the images of two edges overlap (i.e.\ double edges, as in types \fourValentVertex{}-\mergedTriangles{}), or for which two edges cross (type \fatPointOnParallelogram{}). 
    A marked point meeting a simple point condition fixes one floor.
     
     Each elevator edge is adjacent to a marked point. A double point can fix an elevator edge and its adjacent floor (type \fatPointOnVertex{}), or a crossing of an elevator with a floor (type \fatPointOnParallelogram{}), or a double elevator.
\end{lemma}
\begin{proof}
    This follows from the fact that a tropical stable map satisfying $r+2s$ vertically stretched point conditions is \emph{floor decomposed}, i.e.\ every connected component of $\Gamma$ minus its elevator edges which is not just a marked point is a floor and contains precisely one marked point. Furthermore, every elevator edge is adjacent to a marked point followed by another elevator edge such that they form an elevator. This is shown e.g.\ in Section 5, \cite{BM08}. If we now merge $s$ pairs of simple points to form double points, we can obtain the cases listed in the claim.  
\end{proof}

\begin{proposition}\label{prop-vertextypesinfloordecomposedcurve}
    Let $(\Gamma,f)$ be a tropical stable map satisfying $r$ simple and $s$ double vertically stretched point conditions.
    \begin{enumerate}
        \item Assume $(\Gamma,f)$ contains a vertex of type \fatPointOnVertex{}. Then two adjacent edges are of direction $(\pm 1,y)$ for some $y$.

        \item Assume $(\Gamma,f)$ contains a vertex of type \fourValentVertex{}. Then the double edges are elevator edges, and the two other edges are of direction $(\pm 1,y)$ for some $y$. 
        
        \item Assume $(\Gamma,f)$ contains a double vertex of type \allDoubleVertex{}. Then for each of the two components, two adjacent edges are of direction $(\pm 1,y)$ for some $y$.

        \item Assume $(\Gamma,f)$ contains a vertex of type \fourValentVertexWithMergedEdge{}. Then the double edges of the part containing the $4$-valent vertex are elevator edges, and the three other edges are of direction $(\pm 1,y)$ for some $y$.
    \end{enumerate}
\end{proposition}

\begin{proof}
    This follows from Lemma \ref{lem-floordecomposed}: for a $3$-valent vertex, the balancing condition implies that two adjacent edges must be floor edges and one must be an elevator edge. This settles the cases \fatPointOnVertex{} and \allDoubleVertex{}. For case \fourValentVertex{}, assume the double edges were floor edges. As in the proof of Lemma \ref{lem-floordecomposed}, they occur by merging two points to become a double point. Before merging, the double edges must have belonged to two adjacent edges, which must be connected with a shrinking elevator edge. But this is a contradiction, as every elevator edge is adjacent to a marked point itself. Thus the double edges are elevator edges. For the case \fourValentVertexWithMergedEdge{}, notice that it consists of a part of type \fourValentVertex{} crossed by an edge. Using the previous argument, the double edges of the part of type \fourValentVertex{} containing the $4$-valent vertex are elevator edges.
    The claim follows.
\end{proof}

\begin{example}\label{ex-floordecomposedstablemaps}
    Figure \ref{fig-exfloordecomposed} shows a tropical stable map of degree $4$ passing through $5$ double and one simple vertically stretched points, and two tropical stable maps of bidegree $(2,3)$ passing through~$4$ double and one simple vertically stretched point.

    \begin{figure}[t]
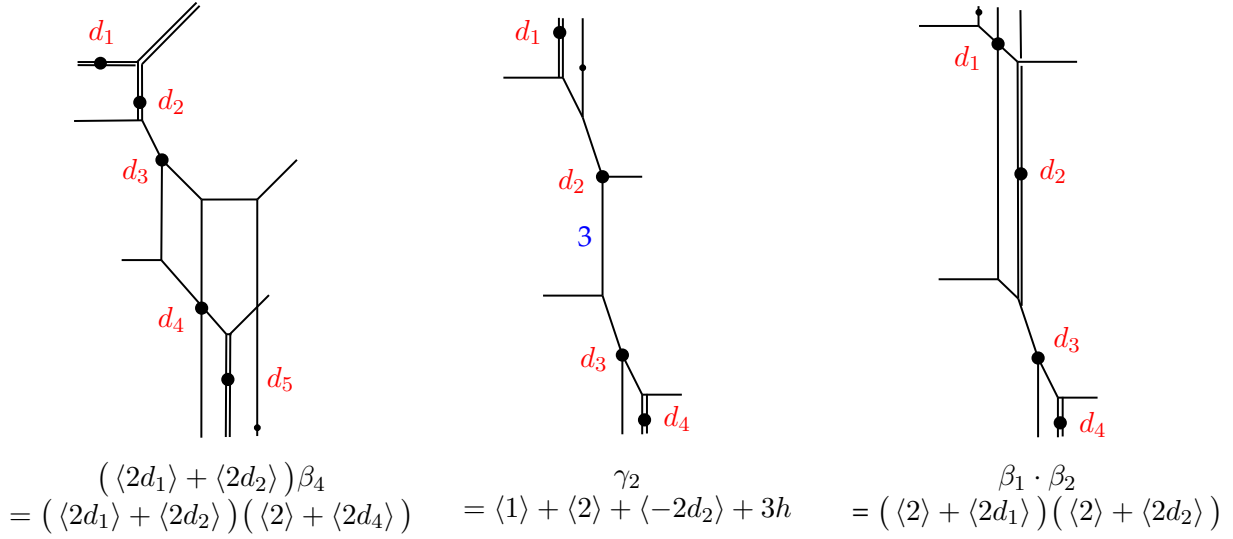

        \centering
        \include{figures/tropCurves_floorDecomposed}        
        \caption{Examples of tropical stable maps passing through vertically stretched simple and double point conditions. The color convention of \cref{not:colors} applies. The curve on the left is the tropicalization of a stable map to a plane quartic curve in $\PP^2$, the middle and right curves are coming from curves in $\PP^1 \times \PP^1$ of bidegree $(2, 3)$. Below each curve we give its quadratically enriched multiplicity according to \cref{def-A1mult-vertices}.} 
        \label{fig-exfloordecomposed}
    \end{figure}
\end{example}

\section{The quadratically enriched multiplicity of a floor decomposed tropical stable map}
\label{section:computations}

In \cref{sec-tropdoublepoints} we determined the structure and local features that tropical stable maps appearing in the context of our tropical enumerative problem can have. In this section we define the quadratically enriched multiplicity $\mult^{\AA^1}(\Gamma, f)$ of a tropical stable map $(\Gamma, f)$ passing through vertically stretched point conditions. For this purpose, we introduce the following notation.

\begin{notation}\label{not-vertexedgeweight}
    Let $(\Gamma, f)$ be a tropical stable map passing through vertically stretched point conditions.
    Let $v$ be a $3$-valent vertex of the underlying graph $\Gamma$ (which is not adjacent to a marked point).    
    By the vertically stretched position of the point conditions, $v$ is adjacent to two edges of weight~$1$ (see Proposition \ref{prop-vertextypesinfloordecomposedcurve}). Denote the weight of the third edge by~$m_v$.

    Let $v$ be a $4$-valent vertex of $\Gamma$ with a double edge, i.e.\ a vertex mapping to type \fourValentVertex{} or a part of \fourValentVertexWithMergedEdge{}. By \cref{prop-vertextypesinfloordecomposedcurve}, the double edges are vertical (and the other edges are of weight $1$). We denote the weights of the double edges by $m_{v}$ and $m_{v}'$. 
    
    Let $m$ be a positive integer. We will write $m^{\AA^1}$ to mean the following quantity in $\GW(k)$:
    \[ m^{\AA^1} = \begin{cases}
        \qinv{m}+
        \frac{m-1}{2}\cdot \h  & \text{if $m$ is odd,} \\
        \frac{m}{2}\cdot \h                      & \text{if $m$ is even.}   
    \end{cases} \]
    Let $m$ be a positive integer and $d \in k^\times$ and define the following quantities: 
    \[
        (m, d)^{\AA^1} = \begin{cases}
        \qinv{1} + \frac{m-1}{2} \big( \qinv{2} + \qinv{-2d} \big) + \frac{m(m-1)}{2} \cdot\h & \text{if $m$ is odd,} \\
        \frac{m}{2} \big(\gw{1} + \gw{-d} \big) + \frac{m(m-1)}{2} \cdot \h    & \text{if $m$ is even}   
    \end{cases}
    \]
    and
    \begin{equation*}
        \gamma(m, d) = \begin{cases}
                \gw{m} + \frac{m-1}{2} \big( \gw{2m} + \gw{-2dm} \big) & \text{if $m$ odd,} \\
                \frac{m}{2} \cdot \h & \text{if $m$ even.}
            \end{cases}
    \end{equation*}
\end{notation}

\begin{definition}\label{def-A1mult-vertices}
    Let $(\Gamma,f)$ be a tropical stable map passing through $r$ simple and $s$ double point conditions which are vertically stretched. Assume the double point $q_i$ is the tropicalization of a pair of conjugate points in the field extension $ K(\sqrt{d_i})$ of the field of Puiseux series~$K$.
    Let $v$ be a vertex of $\Gamma$. Depending on the type of its image under $f$ we define its multiplicity $\mult^{\AA^1}(v)$ as follows.
    \begin{itemize}
        \item Type \fatPointOnVertex{}, then $\mult^{\AA^1}(v) = \gamma(m_v, d_v)$, where $d_v$ defines the quadratic field extension of the double point condition on $v$.
        
        \item Type \fatPointOnParallelogram{}, then $\mult^{\AA^1}(v) = \gw{1}$.
        
        \item Type \fourValentVertex{} or \fourValentVertexWithMergedEdge{} and a twin tree starts at $v$, then $\mult^{\AA^1}(v) = (m_v, d_\elev)^{\AA^1}$, where $m_v$ is the weight of either of the parts of the twin elevator edge $(e, e')$ starting at $v$ and $d_\elev$ defines the field extension of the point condition on that twin elevator.
        
        \item Type \allDoubleVertex{} and part of a twin tree, then let $(e, e')$ be the twin elevator edge incident to $v$, let $m_v$ be the weight of $e$ (and $e'$), let $d_{\elev}$ define the field extension of the double point on $(e, e')$ and write $d_{\floor}$ similarly for the double point on the twin floor incident to $v$. Then 
        \[\mult^{\AA^1}(v) = (m_v, d_{\elev})^{\AA^1} \gw{d_{\elev}^{m_v}} \big(\gw{2} + \gw{2d_{\floor}d_{\elev}} \big) . \]

        \item In all other cases, if the there is a single elevator edge of weight $m_v$ incident to $v$, then $\mult^{\AA^1}(v) = m_v^{\AA^1}$ and if there are two elevator edges (the two parts of a double but not twin edge) of weights $m_v$ and $m'_v$ incident to $v$, then $\mult^{\AA^1}(v) = m_v^{\AA^1} \cdot {m'_v}^{\AA^1}$.
    \end{itemize}
    With this notation we define the quadratically enriched multiplicity (short: quadratic multiplicity) of $(\Gamma, f)$ as
    \[ \mult^{\AA^1}(\Gamma, f) = \prod_{v} \mult^{\AA^1}(v) \cdot \prod_{i \in R} \beta_i, \]
    where the first product runs over all vertices and the second product runs over $R \subseteq \{1, \ldots, s\}$, the subset of indices of double points which are neither of type \fatPointOnVertex{} nor on edges of twin trees and we write $\beta_i=\langle 2\rangle +\langle 2d_i\rangle $.
\end{definition}

The multiplicity in \cref{def-A1mult-vertices} is (almost) a product over vertex contributions -- a shape very familiar to tropical geometers. Nonetheless, we will now present a rearrangement of the multiplicity which treats twin trees a separate parts of the tropical curve. This perspective appears more natural from the point of view of the proof of our main theorem.

\begin{definition}\label{def:twinedgemult}

    The multiplicity of a non-twin edge $e$ of $\Gamma$ of weight $m$ is given by
    \begin{equation*}
        \operatorname{mult}^{\AA^1}(e) = \big( m^{\AA^1} \big)^2 =\begin{cases}
            \qinv{1}+\frac{m^2-1}{2}\h & \text{$m$ odd}\\
            \frac{m^2}{2}\h & \text{$m$ even.}
        \end{cases}
    \end{equation*}
    The multiplicity of a twin edge $e,e'$ with weight $m$ (i.e. both $e$ and $e'$ have weight $m$) carrying a double point condition corresponding to the field extension $K(\sqrt{d}) / K$ is
    \begin{equation*}
        \mult^{\AA^1}(e,e')  = \big( (m, d)^{\AA^1}  \big)^2 = \begin{cases}
           \gw{1}+\frac{m^2-1}{2} \big(\gw{1}+\gw{-d} \big)+\frac{m^4-m^2}{2}\h & \text{if $m$ is odd}\\
           \frac{m^2}{2} \big(\gw{1}+\gw{-d} \big)+\frac{m^4-m^2}{2}\h & \text{if $m$ is even.}
        \end{cases}
    \end{equation*}
\end{definition}

\begin{definition}\label{def-multtwintree}
     Let $(\Gamma,f)$ be a tropical stable map as in \cref{def-A1mult-vertices}.
     Let $\mathcal{T}$ be a twin tree of $(\Gamma,f)$ and let $q_1,\ldots,q_t$ be the double points on the twin edges of $\mathcal{T}$. $\mathcal{T}$ is rooted at $\Gamma$ at a twin edge of weight $m_{\mathrm{root}}$, adjacent to a vertex mapping to type \fourValentVertex{} or a part of \fourValentVertexWithMergedEdge{} (see Definition \ref{def-twintree}). Define $m_\circ$ to be $m_{\mathrm{root}}$ plus the number of unbounded twin elevators in $\cT$.
     We define the quadratically enriched multiplicity of the twin tree $\mathcal{T}$ to be
     \begin{equation*}
        \mult^{\AA^1}(\mathcal{T})=\prod_{(e, e') \text{ elevator of $\cT$}} \mult^{\AA^1}(e,e') \cdot 
        \qinv{2^{t-1}}\cdot
        \sum_{\substack{I\subset \{1,\dots,t\}\\ \vert I\vert\equiv m_{\circ}\mod 2}} \qinv{\prod_{i\in I} d_i}  
        \in \GW(k).
    \end{equation*}
\end{definition}

\begin{remark}
    Some remarks about \cref{def-multtwintree} are in order.
    \begin{enumerate}
        \item Every twin tree has a unique root by \cref{lem-structuredoublepart}, hence $m_\mathrm{root}$ is well defined.
        
        \item We are working with vertically stretched point conditions. Hence every elevator twin edge $(e, e')$ contains a double point and hence the twin edge multiplicity of $(e, e')$ which occurs in the definition of $\mult^{\AA^1}(\cT)$ is indeed defined. 

        \item The simplest twin tree $\cT$ consists of only one unbounded end which is double and this is necessarily an elevator edge by \cref{lem-floordecomposed}. In this case $m_\circ = 1 +1$ and hence the twin tree multiplicity of such a $\cT$ is always $\gw{1}$.
    \end{enumerate}
\end{remark}

\begin{example}
    Consider the twin trees in \cref{fig-twintree}. In each of them we label the double points $q_1, \ldots, q_t$ arbitrarily (except in $\cT_3$, where we specify $q_1$ to lie on the root elevator). Then the multiplicities are as indicated in the figure.
\end{example}

\begin{lemma} \label{lem-A1mult-twintree}
    Let $\cT$ be a twin tree. Then $\mult^{\AA^1}(\cT)$ in the sense of \cref{def-multtwintree} equals the product of vertex multiplicities involved in $\cT$.
\end{lemma}

\begin{proof} Writing the product over multiplicities of vertices in $\cT$, it is easy to see that for every bounded elevator edge $(e, e')$ there are two factors $(m, d)^{\AA^1}$ coming from the multiplicities of the vertices at either end. These combine into one factor of $\mult^{\AA^1}(e, e')$. For unbounded twin elevators there is only a single such factor, however, it is equal to $\gw{1}$ which already happens to be the multiplicity of a twin elevator of weight 1. In formulas, using $V_{\allDoubleVertex{}}$ for the set of all vertices in $\cT$ of type \allDoubleVertex{}:
    \begin{align*}
        &\underbrace{(m_{\mathrm{root}}, d_{\mathrm{root}})^{\AA^1}}_{\text{mult of unique root vertex of type \fourValentVertex{}}} \cdot \prod_{v \in V_{\allDoubleVertex{}}} \mult^{\AA^1}(e, e')\\
        ={} &(m_{\mathrm{root}}, d_{\mathrm{root}})^{\AA^1} \underbrace{\gw{d_{\mathrm{root}}^{m_{\mathrm{root}}}}^2}_{= \gw{1}} \cdot \prod_{v \in V_{\allDoubleVertex{}}} (m_v, d_{\elev})^{\AA^1} \gw{d_{\elev}^{m_v}} \big(\gw{2} + \gw{2d_{\floor}d_{\elev}} \big) \\
        ={} & \gw{d_{\mathrm{root}}^{m_{\mathrm{root}}}} \cdot \prod_{\substack{(e, e') \text{ bounded} \\\text{elevator of } \cT}} \mult^{\AA^1}(e, e') \underbrace{\gw{d_{\elev}^{m_v}}^2}_{= \gw{1}} 
        \cdot \prod_{\substack{(e, e') \text{ unbounded} \\\text{elevator of } \cT}} \underbrace{(1, d_{\elev})^{\AA^1}}_{= \gw{1} = \mult^{\AA^1}(e, e')} \gw{d_{\elev}^1} \\
        & \qquad \cdot \gw{2}^{\# V_{\allDoubleVertex{}}} \prod_{v \in V_{\allDoubleVertex{}}} \big(\gw{1} + \gw{d_{\floor}d_{\elev}} \big) \\
        ={} &\prod_{(e, e') \text{ elevator of } \cT} \mult^{\AA^1}(e, e') \cdot \underbrace{\gw{d_{\mathrm{root}}^{m_{\mathrm{root}}}} \cdot \prod_{\substack{(e, e') \text{ unbounded} \\\text{elevator of } \cT}} \gw{d_{\elev}} }_{= (\ast)} \cdot \gw{2}^{\# V_{\allDoubleVertex{}}} \prod_{v \in V_{\allDoubleVertex{}}} \big(\gw{1} + \gw{d_{\floor}d_{\elev}} \big)
    \end{align*}
    
    Next it is easy to see that the number of type \allDoubleVertex{} vertices in $\cT$ is one less than the number $t$ of double point conditions on $\cT$, and hence
    \[ \gw{2}^{\# V_{\allDoubleVertex{}}} = \gw{2}^{t-1}. \]
    
    Next we let $d_1, \ldots, d_t \in K^\times$ be the elements defining the field extensions of the point conditions in $\cT$, listed in any arbitrary order, and we show that 
    \begin{equation} \label{eq-product_expanded_4}
        \prod_{v \in V_{\allDoubleVertex{}}} \big( \gw{1} + \gw{d_{\floor}d_{\elev}} \big) = \sum_{\substack{I\subset \{1,\dots,t\}\\ \vert I\vert\equiv 0\mod 2}} \qinv{\prod_{i\in I} d_i}.
    \end{equation}
    For $t = 1$ the product on the left is empty and the only summand on the right is the empty product, hence $\gw{1} = \gw{1}$. Now let $t \geq 2$ and pick a vertex $v_0 \in V_{\allDoubleVertex{}}$ which is incident to an unbounded elevator or a floor dual to a simplex of size 1. This is always possible due to the restrictive class of degrees $\Delta$ which we consider in this article. In either case, there is going to be one $d$ among $\{d_{\floor}, d_{\elev} \}$ of $v_0$ which appears in the left hand side of \cref{eq-product_expanded_4} exclusively in the factor corresponding to $v_0$. Without loss of generality order the $d$'s such that $d_1$ is this unique one and $d_2$ the other one. By induction we compute:
    \begin{align*}
        &\prod_{v \in V_{\allDoubleVertex{}}} \big( \gw{1} + \gw{d_{\floor}d_{\elev}} \big) \\
        ={} &\big(\gw{1} + \gw{d_1d_2} \big) 
        \cdot \sum_{\substack{I\subset \{2,\dots,t\}\\ \vert I\vert\equiv 0\mod 2}} \qinv{\prod_{i\in I} d_i} \\
        ={} &\sum_{\substack{I\subset \{2,\dots,t\}\\ \vert I\vert\equiv 0\mod 2}} \qinv{\prod_{i\in I} d_i} 
        + \sum_{\substack{I\subset \{2,\dots,t\}\\ \vert I\vert\equiv 0\mod 2 \\ 2 \in I}} \qinv{\prod_{i\in I \setminus \{2\} \cup \{1\} } d_i} 
        + \sum_{\substack{I\subset \{2,\dots,t\}\\ \vert I\vert\equiv 0\mod 2\\ 2 \not\in I}} \qinv{\prod_{i\in I \cup \{1, 2\}} d_i}  \\
        ={} &\sum_{\substack{I\subset \{1,\dots,t\}\\ \vert I\vert\equiv 0\mod 2 \\ 1 \not\in I}} \qinv{\prod_{i\in I} d_i} 
        + \sum_{\substack{I\subset \{1,\dots,t\}\\ \vert I\vert\equiv 0\mod 2 \\ 1 \in I \text{ and } 2 \not\in I}} \qinv{\prod_{i\in I } d_i} 
        + \sum_{\substack{I\subset \{1,\dots,t\}\\ \vert I\vert\equiv 0\mod 2\\ 1,2 \in I}} \qinv{\prod_{i\in I} d_i} \\
        ={} &\sum_{\substack{I\subset \{1,\dots,t\}\\ \vert I\vert\equiv 0\mod 2}} \qinv{\prod_{i\in I} d_i}.
    \end{align*} 
    This shows \cref{eq-product_expanded_4}. To complete the proof of the lemma it suffices to show that for any $d_i$ we have
    \[ \gw{d_i} \cdot \sum_{\substack{I\subset \{1,\dots,t\}\\ \vert I\vert\equiv 0\mod 2}} \qinv{\prod_{i\in I} d_i} = \sum_{\substack{I\subset \{1,\dots,t\}\\ \vert I\vert\equiv 1\mod 2}} \qinv{\prod_{i\in I} d_i} \]
    and hence all the remaining factors in $(\ast)$ are absorbed into the right hand side of \cref{eq-product_expanded_4}, inducing a total of $m_\circ$ many parity changes. With this final observation, we have successfully transformed the product over vertex multiplicities into the $\mult^{\AA^1}(\cT)$ from \cref{def-multtwintree}.
\end{proof}

The following immediately follows from \cref{lem-A1mult-twintree}.

\begin{corollary}\label{cor-A1mult-twintrees}
    Let $(\Gamma, f)$ be a tropical curve as in \cref{def-A1mult-vertices}. 
    Assume furthermore that there are $T$ twin trees $\mathcal{T}_1,\ldots,\mathcal{T}_T$.
    Then the quadratic multiplicity may also be written as
    \begin{equation} \label{eq-A1mult-twintrees}
        \mult^{\AA^1}(\Gamma,f)=\prod_{i=1}^T
    \mult^{\AA^1}(\mathcal{T}_i)\cdot \prod_{v} \mult^{\AA^1}(v) \cdot \prod_{i\in R} \beta_i,
    \end{equation}
    where 
    \begin{enumerate}
        \item $\mult^{\AA^1}(\mathcal{T}_i)$ denotes the quadratically enriched multiplicity of the twin tree $\mathcal{T}_i$ as defined in \cref{def-multtwintree},

        \item the second product runs over all vertices which are not part of a twin tree and $\mult^{\AA^1}(v)$ is the vertex multiplicity from \cref{def-A1mult-vertices}, and
    
        \item $R \subseteq \{1, \ldots, s\}$ and $\beta_i$ are as in \cref{def-A1mult-vertices}.  
    \end{enumerate}
\end{corollary}

\begin{example}
    If we label the double points and their corresponding quadratic extensions from top to bottom, the quadratically enriched multiplicity of the three tropical stable maps from Example \ref{ex-floordecomposedstablemaps} are given in \cref{fig-exfloordecomposed}.
\end{example}

\begin{remark} \label{rem-A1mult-twintree}
    Another equivalent presentation of the quadratic multiplicity can be obtained by combining the multiplicities of two vertices which are connected by a bounded elevator into a single edge multiplicity. This point of view yields
    \[ \mult^{\AA^1}(\Gamma, f) = \prod_{i=1}^T \mult^{\AA^1}(\mathcal{T}_i) \cdot \prod_{v \text{ of type \fatPointOnVertex{}}} \gamma(m_v, d_v) m_v^{\AA^1} \cdot \prod_{\substack{e \text{ bounded elevator} \\ \text{not incident to type \fatPointOnVertex{}}}} \mult^{\AA^1}(e) \cdot \prod_{i\in R} \beta_i. \]
\end{remark}

\begin{definition}[Tropical quadratically enriched count of curves satisfying rational point conditions and point conditions in quadratic extensions]\label{def-ourcount}
    Consider a configuration of $r$ simple and $s$ double point conditions which are vertically stretched. Assume the double point $q_i$ is the tropicalization of a pair of conjugate points in the field extension $K\subset K(\sqrt{d_i})$. 
    We define the tropical quadratically enriched count of curves satisfying these point conditions to be

    $$ N_\Delta^{\AA^1,\trop}(r,(d_1,\ldots,d_s))\coloneqq \sum_{(\Gamma,f)} \operatorname{mult}^{\A^1}(\Gamma,f),$$
    where the sum goes over all tropical stable maps $(\Gamma,f)$ satisfying the point conditions, and the multiplicity is as defined in \cref{def-A1mult-vertices}.
\end{definition}

The goal of this paper is to prove the quadratically enriched correspondence theorem, see \cref{thm-main} in the introduction, or \cref{thm:correspondence} below.
For this we will use the following theorem to show that the sum of quadratic weights of the log stable maps satisfying the conditions that tropicalize to a given tropical stable map $(\Gamma,f)$ satisfying the conditions equals $\operatorname{mult}^{\AA^1}(\Gamma,f)$ defined in \cref{def-A1mult-vertices}.
A log stable map that tropicalizes to $(\Gamma,f)$ can be described in terms of the combinatorics: the vertices $v$ of $\Gamma$ correspond to local pieces $C_v$ defined over a finite field extension of $k$ (with a map to the toric surface given by the polygon $\Delta_v$ dual to the image of $v$ in $\mathbb{R}^2$) and the edges of $\Gamma$ govern how these local pieces are glued together. Following \cite{Shu06b}, we distinguish \emph{deformation patterns} for the gluing of edges and \emph{refined point conditions} for edges containing a marked point (also defined over a finite field extension of $k$). In \cite{Shu06b}, the methods to compute preimages under tropicalization using local pieces $C_v$, deformation patterns and refined point conditions, are provided. We enrich these techniques by computing the quadratic weight of each preimage. This yields the \emph{quadratic weight} of the log stable map that arises that way, that is the quadratic weight of the underlying stable map, as stated in the following theorem.

\begin{theorem}
\label{theorem:multofatropstablemap}
The different choices for the local pieces together with the deformation patterns and the refined point conditions determine all log stable maps satisfying the point conditions that tropicalize to a given tropical stable map $(\Gamma,f)$ satisfying the conditions.
 Furthermore, the quadratic weight of a log stable map can be computed as the product of quadratic weights of the local pieces times the quadratic weights of the deformation patterns and the refined point conditions.
\end{theorem}

\cref{theorem:multofatropstablemap} essentially follows from the general techniques for correspondence theorems for counts of tropical curves that have been established in \cite{Mi03, Shu04, Shu06b, NS06, AB22, Tyo09, Ran15, Gro16}. 
    The tropical stable map hands us the data of a toric degeneration of the surrounding surface and pieces of algebraic curves satisfying tangency conditions with the toric boundary in the little toric surfaces. For each piece, we will compute the algebra parametrizing the algebraic curves tropicalizing to this local piece. This is a refinement: a mere count is not sufficient for us here, since we need to know precisely over which field extension the algebraic curves live. Moreover, we need to determine their nodes and their quadratic weights. We will do this in \cref{section:computations local pieces}.
    The tangency conditions allow us to glue the pieces to obtain a log stable map. In \cite{Shu04}, this gluing procedure is described in terms of patchworking (see \cite{Viro, IMS09}) via the so-called \emph{deformation patterns} and \emph{refined point conditions}.
    When gluing a higher weight edge, further choices of coefficients are involved and these might live over larger field extensions. We also want to record these larger field extensions. We will compute these field extensions in Proposition \ref{prop-deformpatternandrefinedpoint} for the different choices of deformation patterns and in Proposition \ref{prop:algebrarefinedpointconditions} for the different choices of refined point conditions.
    
    Shustin's deformation patterns tell us how to glue the pieces to obtain log stable maps. Furthermore, with the deformation patterns, we also discover hidden nodes for which we will compute the quadratic weight. Since the quadratic weight of a curve is given as the product of the quadratic weights of the nodes, as in \cite{JPP23} the quadratic weight of a preimage of a tropical stable map under tropicalization is a product of the contributions of the local pieces and the deformation patterns.

\section{Compatibility with existing tropical curve counts}\label{sec-compatible}

Before starting the proof of the main theorem, we use this section to demonstrate that our quadratically enriched count of tropical curves from Definition \ref{def-ourcount} generalizes and unifies many existing tropical curve counts:
\begin{enumerate}
    \item the quadratically enriched count of tropical curves satisfying $k$-rational point conditions from \cite{JPP23}, 
    \item the complex count of tropical curves which provides the plane Gromov-Witten invariants \cite{Mi03}, 
    \item the real count of tropical curves (satisfying real point conditions) which provides the Welschinger invariant \cite{Mi03}, 
    \item and the real count of tropical curves satisfying real point conditions and pairs of complex conjugate point conditions \cite{Shu06b}, providing the more general Welschinger invariant.
\end{enumerate}

We illustrate this in Table \ref{tab:Felixfantasticpicture} by providing our quadratically enriched count and the four counts mentioned above simultaneously for the case of plane cubics. 
This shows how universally usable our quadratically enriched count from \cref{def-ourcount} is.  By the end of this section, \cref{thm-generalizing} from the introduction will be proven.

\begin{table}[t]
    \centering
    \input{figures/table_cubics}
    \caption{Plane cubic curves through a fixed configuration of 3 double points and 2 simple points.}
    \label{tab:Felixfantasticpicture}
\end{table}
Furthermore, setting $s = \lfloor \frac{n}{2} \rfloor$ provides a universal count for all multiquadratic sigma. In fact, as we show in Proposition \ref{prop-specializetokrational}, the quadratic count $N_\Delta^{\AA^1,\trop}(r+2, (d_1, \ldots, d_{s-1}))$ can be computed from $N_\Delta^{\AA^1,\trop}(r, (d_1, \ldots, d_s))$ by setting $d_s = 1$ in the formula for $N_\Delta^{\AA^1,\trop}(r, (d_1, \ldots, d_s))$.

\subsection{Compatibility with the quadratically enriched count of tropical curves}\label{sec-compatibilityquadraticallyenriched}

Recall that the quadratically enriched count of tropical curves with $k$-rational point conditions was introduced in \cite{JPP23} and we recall the definition:

\begin{definition}
    Assume $(\Gamma,f)$ is a floor decomposed rational tropical stable map of degree $\Delta$ passing through $n(\Delta)$ vertically stretched points, which we view as tropicalizations of $k$-rational points.
    Following \cite{JPP23}, the quadratically enriched multiplicity of $(\Gamma,f)$ is $$\mult^{\mathbb{A}^1}(\Gamma,f) = \prod_v \ (m_v)^{\mathbb{A}^1},$$
    where the product goes over all (necessarily $3$-valent) vertices $v$ of the image $f(\Gamma)$ and $m_v$ denotes the weight of their adjacent vertical edge as in Notation \ref{not-vertexedgeweight}.
    
    The \emph{quadratically enriched tropical count} is defined as
    $$ N_\Delta^{\AA^1,\trop}\coloneqq \sum_{(\Gamma,f)} \operatorname{mult}^{\A^1}(\Gamma,f)$$
    where the sum goes over all tropical stable maps passing through the points.
    For the algebraic counterpart of the count of rational curves with $k$-rational point conditions we use the notation~$N_\Delta^{\AA^1}$.
\end{definition}

\begin{theorem}[Quadratically enriched correspondence theorem, \cite{JPP23}]
    \label{thm:generalize_JPP}

    The quadratically enriched count of curves satisfying $k$-rational point conditions (see Definition \ref{def-quadraticallyenrichedcount}) equals its tropical counterpart,
    $$ N_\Delta^{\AA^1,\trop}= N_\Delta^{\AA^1}.
    $$
\end{theorem}

Our quadratically enriched tropical count from Definition \ref{def-ourcount} extends this count: 

\begin{proposition}\label{prop-specializetokrational_s_0}
    If $s = 0$, then our count recovers that of \cite{JPP23}, i.e.
    \[ N_\Delta^{\AA^1, \trop}(r, \emptyset) = N_\Delta^{\AA^1, \trop}. \]
\end{proposition}

\begin{proof}
    Let $(\Gamma, f)$ be a tropical stable map passing through $r = n$ simple points. And let $v$ be a vertex. Note that $\mult^{\AA^1}(v)$ in the sense of \cref{def-A1mult-vertices} is just $(m_v)^{\AA^1}$, which is the multiplicity in the sense of \cite{JPP23}. Moreover, the second product in $\mult^{\AA^1}(\Gamma, f)$ in the sense of \cref{def-A1mult-vertices} is empty and thus the claim follows.
\end{proof}

\begin{proposition}\label{prop-specializetokrational}
    If $d_s \in k^\times$ is a square, then 
    \[ N_\Delta^{\AA^1, \trop}(r, (d_1, \ldots, d_s)) = N_\Delta^{\AA^1, \trop}(r+2, (d_1, \ldots, d_{s-1})). \]
    In particular, if all $d_i$ are squares, then $N_\Delta^{\AA^1, \trop}(r, (d_1, \ldots, d_s)) = N_\Delta^{\AA^1, \trop}$.
\end{proposition}

\begin{proof}
    Given a tropical stable map $(\Gamma,f)$ passing through vertically stretched $r$ simple and $s$ double points, we consider a deformation of the double point condition $q_s$ where we replace the double point by two close-by simple points $p_{r+1}$ and $p_{r+2}$, and we deform $(\Gamma,f)$ along. 

    We define the \emph{deformed multiplicity} $\mult^D(\Gamma, f)$ of $(\Gamma,f)$  
    to be the sum of the quadratically enriched multiplicities of all tropical stable maps close to $(\Gamma,f)$ that pass through the deformed points. 
    We claim that the deformed multiplicity is equal to $\mult^{\AA^1}(\Gamma, f)$ with $d_s$ set to be a square. The statement on the tropical counts follows.
    To show the claim, we work with the formula for $\mult^{\AA^1}(\Gamma,f)$ from \cref{cor-A1mult-twintrees} and distinguish cases depending on which factor of \cref{eq-A1mult-twintrees} the double point $q_s$ was contributing to. 

    \textbf{Case 1:} If $q_s$ was part of a twin tree $\cT$, then the deformation of $\cT$ is not a twin tree anymore. Moreover, the deformation of the stable map is unique since the two branches of $\cT$ were indistinguishable and hence $\mult^D(\Gamma, f)$ is the quadratic multiplicity of the unique deformation. In this case the remaining double points $q_1, \ldots, q_{t-1}$ which used to be on $\cT$ now contribute to the third product in \cref{eq-A1mult-twintrees} of $\mult^D(\Gamma, f)$ and their contribution is
    \[ \prod_{i = 1}^{t-1} \beta_i = \gw{2^{t-1}} \prod_{i = 1}^{t-1} \big( \gw{1} + \gw{d_i} \big) = \gw{2^{t-1}} \sum_{I \subseteq \{1, \ldots, t-1\}} \gw{\prod_{i \in I} d_i }. \]
    The last equality can easily be verified by induction. Moreover, every twin elevator $(e, e') \in \cT$ with weight $(m, m)$ will contribute a total factor of $(m^{\AA^1})^4$ to the second product of $\mult^D(\Gamma, f)$ (see \cref{eq-A1mult-twintrees}).
    
    We now investigate what happens when we set $d_s$ to be square in the multiplicity of the twin tree and show that it is the same as the deformed multiplicity. Without loss of generality we assume $q_s$ occurs with index $t$ in $\cT$, i.e. we set $d_s = d_t$ to be a square. Recall the definition of the quadratic multiplicity of $\cT$ 
    \[\mult^{\AA^1}(\cT) =\prod_{(e, e') \text{ elevator of $\cT$}} \mult^{\AA^1}(e,e') \cdot \qinv{2^{t-1}}\cdot \sum_{\substack{I\subset \{1,\dots,t\}\\ \vert I\vert\equiv m_{\circ}\mod 2}} \qinv{\prod_{i\in I} d_i}.  \] 
    With $\gw{d_t} = \gw{1}$ the sum on the right hand side becomes
    \begin{align*}
        \sum_{\substack{I \subseteq \{1, \ldots, t\} \\ |I| \equiv m_{\circ}\mod 2} } \gw{\prod_{i \in I} d_i} 
        =& 
        \sum_{\substack{ I \subseteq \{1, \ldots, t\} \\ |I| \equiv m_{\circ}\mod 2 \\ t \in I }} \gw{\prod_{i \in I} d_i} 
        + \sum_{\substack{ I \subseteq \{1, \ldots, t\} \\ |I| \equiv m_{\circ}\mod 2 \\ t \not\in I }} \gw{\prod_{i \in I} d_i} \\
        =& 
        \gw{d_t} \sum_{\substack{I \subseteq \{1, \ldots, t-1\} \\ |I| \not\equiv m_{\circ}\mod 2}} \gw{\prod_{i \in I} d_i} 
        + \sum_{\substack{I \subseteq \{1, \ldots, t-1\} \\ |I| \equiv m_{\circ}\mod 2}} \gw{\prod_{i \in I} d_i} \\
        =&
        \sum_{I \subseteq \{1, \ldots, t-1\}} \gw{\prod_{i \in I} d_i }. 
    \end{align*}
For ease of notation we set $D = \sum_{I \subseteq \{1, \ldots, t-1\}} \gw{\prod_{i \in I} d_i }$.
    We need to show that 
    \[D\cdot \prod_{(e,e')} \mult^{\AA^1}(e,e')=D\cdot \prod_{(e,e')} (m_e^{\A^1})^4.\]
    To do this, let $(e, e') \in \cT$ be the twin elevator which contains the double point $q_j$ and let $(m, m)$ be the weight of $(e, e')$.  For $j = 1, \ldots, t-1$ we compute 
    \begin{align*}
        \big( \gw{1} + \gw{-d_j} \big) \cdot D
        &= \sum_{I \subseteq \{1, \ldots, t-1\}} \gw{\prod_{i \in I} d_i } + \sum_{I \subseteq \{1, \ldots, t-1\}} \gw{-\prod_{i \in I} d_i } 
        = 2^{t-1} \h  
        = \h \cdot D
    \end{align*}
    From this and the definition of $\mult^{\AA^1}(e, e')$
    it follows that $\mult^{\AA^1}(e, e') \cdot D = (m^{\AA^1})^4 \cdot D$ in this case. 
    Finally, for $j = t$, we obtain directly $\mult^{\AA^1}(e, e') = (m^{\AA^1})^4$ and this finishes the proof of the first case.

    \textbf{Case 2:} The double point $q_s$ was on a vertex of type~\fatPointOnVertex{}. If we separate the double point into two close-by simple points, they will be on two of the edges adjacent to the vertex, but we have just one local deformation for $(\Gamma,f)$.
    For Definition~\ref{def-A1mult-vertices} applied before the deformation,
    the vertex of type~\fatPointOnVertex{} contributes to $\mult^{\AA^1}(\Gamma, f)$ with a factor of $\gamma(m, 1) = m^{\AA^1}$. But this is the same as the contribution of the vertex after the deformation to $\mult^D(\Gamma, f)$ and we can conclude that again $\mult^{\AA^1}$ and $\mult^D$ are equal. 

    \textbf{Case 3:} Consider $q_s$ to be a double point like type \fatPointOnParallelogram{} or on any double edge which is not part of a twin tree.
    Deforming the double point into two close-by simple points leads to two options to locally deform the tropical stable map accordingly, which  have the same quadratically enriched multiplicity, as they only differ combinatorially by which edge passes through which of the two simple points in the deformation.
    
    The specialization of the factor of $\beta_i$ arising for any such double point in $\mult^{\AA^1}(\Gamma, f)$ 
    yields a factor of $\langle 2 \rangle + \langle 2 d_i\rangle = \langle 2 \rangle + \langle 2 \rangle = \langle 1 \rangle + \langle 1\rangle$ (see \cref{eq:gwrealtion2}). We thus have a bijection between local deformations and summands, and the contribution to the quadratically enriched multiplicity of each deformation equals the contribution of a summand. 

    In every case, a possible $4$-valent vertex of type \fourValentVertex{} (or part of \fourValentVertexWithMergedEdge{}) is deformed into two $3$-valent vertices. However, the $4$-valent vertex' contribution to the quadratic multiplicity is $\big((m_v)^{\AA^1}\cdot(m_v')^{\AA^1} \big)$, which is already the product of the contributions of the two deformed $3$-valent vertices. 
\end{proof}

\subsection{Compatibility with the complex count of tropical curves}
It was already shown in \cite{JPP23} that the quadratically enriched count of tropical curves satisfying only $k$-rational point conditions extends both the complex count of tropical curves and the real count of tropical curves satisfying real point conditions. It thus follows easily from \cref{prop-specializetokrational_s_0} that also our count from Definition \ref{def-ourcount} extends these counts as well. For the sake of completeness, we repeat these statements in the following two subsections.

\begin{definition}
    Assume $(\Gamma,f)$ is a floor decomposed rational tropical stable map of degree $\Delta$ passing through $n(\Delta)$ vertically stretched points, which we view as tropicalizations of complex points.
 Following \cite{Mi03}, the complex multiplicity of $(\Gamma,f)$ is $$\mult_{\mathbb{C}}(\Gamma,f) = \prod_v \ m_v,$$
 where the product goes over all (necessarily $3$-valent) vertices $v$ of the image $f(\Gamma)$ and $m_v$ denotes the weight of their adjacent vertical edge as in Notation \ref{not-vertexedgeweight}.

We define $$ N_\Delta^{\trop}\coloneqq \sum_{(\Gamma,f)} \operatorname{mult}_{\CC}(\Gamma,f),$$ 
where the sum goes over all tropical stable maps passing through the points.

\end{definition}
\begin{theorem}[Correspondence Theorem, \cite{Mi03}]
The complex count of rational tropical curves yields the Gromov-Witten invariant $N_\Delta$ counting complex curves in the toric surface defined by $\Delta$ of curve class defined by $\Delta$ passing through $n(\Delta)$ points in general position:
$$N_\Delta^{\trop}= N_\Delta.$$    
\end{theorem}

The following is a corollary of Proposition \ref{prop-specializetokrational_s_0}:

\begin{corollary}
    \label{cor-Mikhalkincount}
    Assume $k=\mathbb{C}$ (or more generally, an algebraically closed field) and necessarily all point conditions are defined over $k$.
    Then the rank of our quadratically enriched count from Definition \ref{def-ourcount} equals the complex count of tropical curves:
   $$ \operatorname{rk}\Big(N_\Delta^{\AA^1,\trop}(r,(d_1,\ldots,d_s)) \Big)=N_\Delta^{\trop}. $$ 
\end{corollary}

\subsection{Compatibility with the tropical Welschinger count for real point conditions}
Again, it was already shown in \cite{JPP23} that the quadratically enriched count of tropical curves satisfying only $k$-rational point conditions extends both the complex count of tropical curves and the real count of tropical curves satisfying real point conditions. It thus follows easily from the Section \ref{sec-compatibilityquadraticallyenriched} that also our count from Definition \ref{def-ourcount} extends these counts. For the sake of completeness, we repeat this statement.

\begin{definition}
    Assume $(\Gamma,f)$ is a floor decomposed rational tropical stable map of degree $\Delta$ passing through $n(\Delta)$ vertically stretched points, which we view as tropicalizations of real points.
    Following \cite{Mi03}, the real multiplicity $\mult_{\mathbb{R}}(\Gamma,f)$ is $0$ if there is an edge of even weight and $1$ else.
    We define 
    $$ W_\Delta^{\trop}\coloneqq \sum_{(\Gamma,f)} \operatorname{mult}_{\RR}(\Gamma,f)$$
    where the sum goes over all tropical stable maps passing through the points.
\end{definition}

\begin{theorem}[Correspondence Theorem, \cite{Mi03}]
    The real count of rational tropical curves yields the Welschinger invariant $W_\Delta$ which is the signed count of real curves in the toric surface defined by $\Delta$ of curve class defined by $\Delta$ passing through $n(\Delta)$ real points in general position:
    $$W_\Delta^{\trop}= W_\Delta.$$    
\end{theorem}

The following is a corollary of Proposition \ref{prop-specializetokrational}:

\begin{corollary}
    \label{cor-MikhalkincountR}
    Assume $k=\mathbb{R}$ and all point conditions are defined over $k$, that is all $d_i$ are positive.
    Then the signature of our count from Definition \ref{def-ourcount} equals the signed real count of tropical curves:
   $$ \operatorname{sgn} \Big(N_\Delta^{\AA^1,\trop}(r,(d_1,\ldots,d_s)) \Big)=W_\Delta^{\trop}. $$
\end{corollary}

\subsection{Compatibility with the tropical count of real curves satisfying real point conditions and pairs of complex point conditions}

\begin{definition}[Real  multiplicity]
    Let $(\Gamma,f)$ be a tropical stable map passing through $r$ simple and $s$ double point conditions which are vertically stretched.
    Assume $(\Gamma,f)$ has no double edges which are not part of a twin tree, no vertices of type \fatPointOnParallelogram{}, and no edges of even weight.
    Then, following \cite{Shu06b}, the real multiplicity of $(\Gamma,f)$ is expressed in terms of Newton subdivision dual to the image $f(\Gamma)$, it equals 
    $$\mult_{\mathbb{R}}(\Gamma,f)= (-1)^{a+b}\cdot 2^{-b-T}\cdot \prod_\Delta \mbox{Area}(\Delta),$$
    where 
    \begin{enumerate}
        \item $a$ is the number of interior lattice points in triangles dual to vertices in the image $f(\Gamma)$,
        \item $b$ equals the number of interior vertices of twin trees,
        \item $T$ is the number of twin trees,
    \end{enumerate}
    and the product goes over all triangles of even lattice area or dual to a vertex of type \fatPointOnVertex{} (see (2.12) in \cite{Shu06b}).

    Note that double edges which are not part of a twin tree, vertices of type \fatPointOnParallelogram{} and edges of even weight do not appear in the real count in \cite{Shu06b}, so we may view their real multiplicity as $0$.

    We define $$ W_\Delta^{\trop}(r,s)\coloneqq \sum_{(\Gamma,f)} \operatorname{mult}_{\RR}(\Gamma,f)$$
    where the sum goes over all tropical stable maps passing through the $r$ simple and $s$ double points.
\end{definition}

\begin{theorem}[Correspondence Theorem, \cite{Shu06b}]
    The tropical count of real curves satisfying real point conditions and pairs of complex conjugate point conditions equals the Welschinger invariant $W_\Delta(r,s)$ which is the signed count of real rational curves in the toric surface defined by $\Delta$ of curve class defined by~$\Delta$ passing through $r$ real points and $s$ pairs of complex conjugate points:
    $$ W_\Delta^{\trop}(r,s)=W_\Delta(r,s).$$
\end{theorem}

\begin{proposition}
    \label{prop-Shustincount}
    Assume $k=\mathbb{R}\subset \mathbb{C}$ are the field extensions for all of the double points (i.e.\ $d_i=-1$ for all $i$).
    Then the signature of our count from Definition \ref{def-ourcount} equals the tropical count of real curves satisfying real point and pairs of complex conjugate point conditions:
    $$ \operatorname{sgn} \Big(N_\Delta^{\AA^1,\trop}(r,(d_1,\ldots,d_s)) \Big)= W_\Delta^{\trop}(r,s).$$
\end{proposition}

\begin{proof}

Again we work with the version of the quadratic multiplicities from \cref{cor-A1mult-twintrees}.
Let $\mathcal{T}$ be a twin tree of $(\Gamma,f)$. Assume that there are $t$ double points on $\mathcal{T}$. 
 Assume furthermore that the $t-1$ bounded pairs of twin edges $e_1,e_1',\ldots,e_{t-1},e_{t-1}'$ (see Lemma \ref{lem-pointsondoubleparts}) of $\mathcal{T}$ are of weights $m_1,\ldots,m_{t-1}$.
Every vertex in the interior of a twin tree is dual to the Minkowski sum of two equal triangles of lattice height $1$ and base length $m$, where $m$ equals the weight of the dual double elevator edges. The area of this triangle is thus $4m$. It contains $m-1$ interior lattice points. 
The end vertex of the twin tree, which is necessarily of type \fourValentVertex{}, is dual to a triangle of area $2m_{\mathrm{root}}$, if $m_{\mathrm{root}}$ is the weight of its adjacent double edges. 
The double edge adjacent to the root vertex is adjacent to two vertices, one dual to a triangle of height $1$ which accordingly has no interior lattice points, and one dual to a triangle with $m_{\mathrm{root}}-1$ many interior lattice points.
Every other interior edge of the twin tree is adjacent to two vertices dual to triangles with the same number of interior lattice points. Thus, the sign contribution we obtain from the lattice points inside triangles is $-1$ if 
$m_{\mathrm{root}}$ is even and $+1$ if $m_{\mathrm{root}}$ is odd. 
Altogether, we obtain a contribution of
$$(-1)^{t-1+\delta} \cdot \left(\prod_{j=1}^{t-1}m_j^2 \right)\cdot  2^{-t} \cdot 2 \cdot 4^{t-1}
= (-2)^{t-1} \cdot (-1)^{\delta}\cdot \left(\prod_{j=1}^{t-1}m_j^2 \right),$$ 
where $\delta=1$ if  $m_{\mathrm{root}}$ is even and $0$ if $m_{\mathrm{root}}$ is odd. 
The number $t$ of points on the twin tree equals the number of ends. Each floor has two ends, thus $t$ is congruent to the number of unbounded elevators mod 2. The parity of $t-1+\delta$ equals the parity of $t+m_{\mathrm{root}}$, which by the above mod 2 equals $m_{\mathrm{root}}$ plus the number of unbounded elevators, i.e.\ $m_{\circ}$. (Recall that by Definition \ref{def-multtwintree}, $m_{\circ}$ exactly equals $m_{\mathrm{root}}$ plus the number of unbounded elevators of the twin tree.) 
Thus the sign contribution for the real multiplicity is $1$ if $m_{\circ}$ is even and $-1 $ if $m_{\circ}$ is odd. 

We observe that the factors of $m_j^2$ also occur in the definition of the quadratically enriched multiplicity of a twin tree in \cref{def-multtwintree}, when specialized to the reals.
Furthermore, this definition also contains a factor of $\langle 2^{t-1}\rangle$ and a factor of $ \sum \qinv{\prod_{i\in I} d_i}$, where the sum runs over all sets $I\subset \{1,\dots,t\}$, such that $\vert I\vert\equiv m_{\circ}\mod 2$. For the latter, notice that each summand equals $\langle 1 \rangle$ if $m_{\circ}$  is even and $\langle -1 \rangle$ if $m_{\circ}$ is odd. Thus, as in the proof of Proposition \ref{prop-specializetokrational}, we only have to count how many odd or even subsets there are of a set of size $t$. Using binomial coefficients, it is easy to see that there are always $2^{t-1}$ summands. If $t=1$, this is just one summand and the previous factor of $\langle 2^{t-1}\rangle$ is just $\langle 1 \rangle$. If $t>1$, there is an even number of summands, and so even after multiplying with $\langle 2^{t-1} \rangle$, we end up having a sum with $2^{t-1}$ summands of $\langle  1 \rangle$ if $m_{\circ}$ is even and $\langle -1 \rangle$ if $m_{\circ}$ is odd.
Taking the signature, we thus obtain $\pm (2)^{t-1} \cdot (\prod_{j=1}^{t-1}m_j^2)$, where the sign is $ +1 $ if $m_{\circ}$  is even and $-1 $ if $m_{\circ}$ is odd.
This equals to the above
 expression for the real multiplicity of the twin tree.

For a vertex $v$ of type \fatPointOnVertex{} (which is adjacent to edges of odd weight only), the specialization to the reals of the contribution $\gamma(m_v, -1)$ equals
$\langle m_v \rangle + \frac{m_v-1}{2} \big(\langle 2m_v \rangle + \langle 2m_v\rangle \big)$, which has signature $m_v$. This in turn is equal to the area of its dual triangle, which is also its contribution to the real multiplicity.

If $(\Gamma,f)$ has a vertex of type \fatPointOnParallelogram{} or double edges which are not part of a twin tree, then there is a factor of $\beta_i$ in its quadratic multiplicity and specializing to $d_i = -1$ we see $\beta_i = \langle 2\rangle +\langle -2\rangle=h$, which has signature $0$.
An edge of even weight $m$ which is not part of a twin tree produces a factor of $\frac{m}{2}\cdot h$ in the quadratic multiplicity. This has signature $0$ over $\mathbb{R}$ again.
This is in accordance with the real multiplicity.

Other $3$-valent unmarked vertices of $\Gamma$ contribute a factor of $1$ to the quadratic multiplicity to the reals, as well as to the real multiplicity. 

A $4$-valent vertex of $\Gamma$ with a double edge which is not the root of a twin tree occurs only in connection with factors $\beta_i$, so the quadratic multiplicity of a tropical stable map with such a vertex has signature $0$, which equals its real multiplicity.

It follows that if $k\subset\mathbb{R}$, then the signature of $\mult^{\mathbb{A}^1}(\Gamma,f)$ equals $\mult_{\mathbb{R}}(\Gamma,f)$.
\end{proof}

\section{Quadratic enrichments of local pieces}
\label{section:computations local pieces}

We are in \cref{setting}. Let $(\Gamma, f)$ be a tropical stable map subject to vertically stretched simple and double point conditions. In Lemma \ref{lem-casesfatpoint} we gave a complete list of possible combinatorial types for the local picture around a vertex of $f(\Gamma)$.  
In this section we go over these combinatorial types and compute the $k$-algebra of algebraic lifts and their quadratic weight.
The overview over the computational results of this section is given in Tables~\ref{tab:possibilities} and~\ref{tab:possibilities2}.

The strategy in this section is as follows. Let $v$ be a vertex of the image plane tropical curve and let $v$ be fixed by some point conditions. We choose coordinates for the algebraic lifts of these point conditions. Let $F_v$ denote the finite étale $k$-algebra determined by the algebraic lift of the local piece of $f(\Gamma)$ around $v$, and $C_v$ the unique algebraic lift defined over $F_v$. To compute the quadratic weight $\Wel(C_v) \in \GW(F_v)$, we use the formula in Remark \ref{remark:computingWel}. For each node of $C_v$, this requires computing the norm (in the sense of Galois theory) of the negative determinant of the Hessian evaluated at the node. Since the norm is the product of Galois conjugates, we will express it directly as such in the computations.

Throughout this section we use \cref{not:colors} when visualizing the local situation at a vertex. 

\begin{notation} \label{not:colors}
    We will draw (pieces of) tropical curves with the following conventions:
    \begin{itemize}
        \item Simple point conditions are drawn as small dots.
        \item Double point conditions are drawn as big dots, with $d$ defining the field extension of the algebraic point conditions indicated next to them in red (if relevant).
        \item Edge weights are given in blue. Double edges have two weights. If no weight is indicated then it is assumed to be equal to 1.
        \item Coordinates of algebraic lifts of point conditions are given in green (if relevant).
    \end{itemize}
\end{notation}

\def\extraheight{8mm}
\settowidth\rotheadsize{in $\GW(C_v)$}
\begin{table}
    \caption{Contributions from vertices in tropical curves subject to simple and double point conditions.}
    \label{tab:possibilities}    
	\tikzset{every picture/.style={line width=0.75pt}} 
	\begin{tabular}{c|c|c|c|c}
        type & pictogram & $k$-algebra $F_v$ of algebraic lifts $C_v$ & \rotcell{$\Wel(C_v)$ \\ in $\GW(F_v)$} & reference \\
		\hline
        & 
		\begin{tikzpicture}[baseline={([yshift=-.5em]current bounding box.center)}]
			\path[draw] (0,0) --+ (180:0.8);
			\path[draw] (0,0) --+ (270:0.8);
			\path[draw] (0,0) --+ (45:0.8);
			
			\thinpoint{-0.5, 0}
			\thinpoint{0, -0.5}
            \draw (0.1, -0.7) node[anchor = west] {\edgeweight{$m$}};
		\end{tikzpicture}  
        & $k$ & $\gw{1} $ & \label{dummy} \makecell{\cref{lem:thin_triangle} \\ or \cite{JPP23} \\ Lem.~5.5 and~5.6} \\[\extraheight]

        &
        \begin{tikzpicture}[baseline={([yshift=-.5em]current bounding box.center)}]
			\path[draw] (0, -0.7) -- (0, 0.7);
			\path[draw] (-0.5, -0.5) -- (0.5, 0.5);
			\thinpoint{0, -0.5}
			\thinpoint{-0.3, -0.3}
            \draw (0.1, -0.6) node[anchor = west] {\edgeweight{$m$}};
		\end{tikzpicture}
		& $k$ & $\gw{1} $ & \cref{lem:par1} \\[\extraheight]
        
        \fatPointOnVertex{} & 
        \begin{tikzpicture}[baseline={([yshift=-.5em]current bounding box.center)}]
            \path[draw] (0,0) --+ (180:0.8);
            \path[draw] (0,0) --+ (270:0.8);
            \draw (0.1, -0.7) node[anchor = west] {\edgeweight{$m$}};
            \path[draw] (0,0) --+ (45:0.8);
            
            \fatpoint{0, 0}
            \draw (0.1, 0) node [anchor = west] {\quadextension{$d$}};
        \end{tikzpicture}
        & \makecell{Fixed locus $L^\phi$ where \\
        $L\coloneq \quotient{k[\iota, u]}{\big(\iota^2 - d, u^m-\eta / \conj{\eta} \big)}$ \\
        and \\
        $\begin{array}{cccc}
        \phi :  & L         &\longrightarrow    & L \\
                &\iota   &\longmapsto        & -\iota \\
                &u &\longmapsto   &\frac{\conj\eta}{\eta}u^{m-1}\\
        \end{array}$} & $\gw{1}$ & \cref{lem:fatpoint_thintriangle}\\[\extraheight] 

        \fatPointOnParallelogram{} &
        \begin{tikzpicture}[baseline={([yshift=-.5em]current bounding box.center)}]
			\path[draw] (0, -0.7) -- (0, 0.7);
			\path[draw] (-0.5, -0.5) -- (0.5, 0.5);
			\fatpoint{0, 0}
            \draw (0, 0.1) node[anchor = south east] {\quadextension{$d$}};
            \draw (0.1, -0.6) node[anchor = west] {\edgeweight{$m$}};
		\end{tikzpicture} 
        & $\quotient{k[\iota]}{(\iota^2 - d)}$ & $\gw{1}$ & \cref{lem:par5} \\[\extraheight]

        \fourValentVertex{} &
        \begin{tikzpicture}[baseline={([yshift=-.5em]current bounding box.center)}]
			\path[draw] (0,0) --+ (180:0.8);
			\path[draw, double] (0,0) --+ (270:1.2);
            \draw (-0.1, -1.1) node[anchor = east] {\edgeweight{$m_1$}};
            \draw (0.1, -1.1) node[anchor = west] {\edgeweight{$m_2$}};
			\path[draw] (0,0) --+ (45:0.8);
			
			\thinpoint{-0.5, 0}
			\fatpoint{0, -0.5}
            \draw (0.1, -0.5) node[anchor = west] {\quadextension{$d$}};
		\end{tikzpicture} 
        & $\begin{cases}
            \quotient{k[\iota]}{(\iota^2 - d)} & \text{if } m_1 \neq m_2, \\
            k & \text{if } m_1 = m_2.
        \end{cases}$ & $\gw{1} $ & \cref{lem:type_C} \\[\extraheight]

        \parallelogramWithOneDoubleEdge{} &
        \begin{tikzpicture}[baseline={([yshift=-.5em]current bounding box.center)}]
			\path[draw] (0, -1.2) -- (0, 0.7);
			\path[draw, double] (-0.5, -0.5) -- (0.5, 0.5);
			\thinpoint{0, -0.5}
			\fatpoint{-0.3, -0.3}
            \draw (0.1, -1.1) node[anchor = west] {\edgeweight{$m$}};
		\end{tikzpicture} or
        \begin{tikzpicture}[baseline={([yshift=-.5em]current bounding box.center)}]
			\path[draw, double] (0, -1.2) -- (0, 0.7);
			\path[draw] (-0.5, -0.5) -- (0.5, 0.5);
			\fatpoint{0, -0.5}
			\thinpoint{-0.3, -0.3}
            \draw (-0.1, -1.1) node[anchor = east] {\edgeweight{$m_1$}};
            \draw (0.1, -1.1) node[anchor = west] {\edgeweight{$m_2$}};
		\end{tikzpicture}
		& $\begin{cases}
            \quotient{k[\iota]}{(\iota^2 - d)} & \text{if } m_1 \neq m_2, \\
            k & \text{if } m_1=m_2.
        \end{cases}$ & $\gw{1} $ & \cref{lem:par3} \\[\extraheight]

        \parallelogramWithTwoDoubleEdges{} &
        \begin{tikzpicture}[baseline={([yshift=-.5em]current bounding box.center)}]
			\path[draw, double] (0, -1.2) -- (0, 0.7);
			\path[draw, double] (-0.5, -0.5) -- (0.5, 0.5);
			\fatpoint{0, -0.5}
			\fatpoint{-0.3, -0.3}
            \draw (-0.1, -1.1) node[anchor = east] {\edgeweight{$m_1$}};
            \draw (0.1, -1.1) node[anchor = west] {\edgeweight{$m_2$}};
		\end{tikzpicture} 
        & $k$ & $\gw{1} $ & \cref{lem:par4} \\[\extraheight]

        \triangleWithMergedEdge{} &
        \begin{tikzpicture}[baseline={([yshift=-.5em]current bounding box.center)}]
            \path[draw] (0,0) --+ (180:0.8);
			\path[draw] (0,0) --+ (270:1.2);
			\path[draw] (0,0) --+ (45:0.8);
            \path[draw] (-0.8, -0.05) -- (0.8, -0.05);

            \fatpoint{-0.5, -0.025}
            \draw (-0.5, 0.1) node[anchor = south east] {\quadextension{$d$}};
            \thinpoint{0, -0.5}
            \draw (0, -1.1) node[anchor = west] {\edgeweight{$m$}};
        \end{tikzpicture} or
        \begin{tikzpicture}[baseline={([yshift=-.5em]current bounding box.center)}]
            \path[draw] (0,0) --+ (180:0.8);
			\path[draw] (0,0) --+ (270:1.2);
			\path[draw] (0,0) --+ (45:0.8);
            \path[draw] (-0.05, -1.2) -- (-0.05, 0.8);

            \thinpoint{-0.5, 0}
            \fatpoint{-0.025, -0.5}
            \draw (0.1, -0.5) node[anchor = west] {\quadextension{$d$}};
            \draw (-0.15, -1.1) node[anchor = east] {\edgeweight{$m_1$}};
            \draw (0.1, -1.1) node[anchor = west] {\edgeweight{$m_2$}};
        \end{tikzpicture}
        & $\quotient{k[\iota]}{(\iota^2 - d)}$ & $\gw{1}$ & \cref{lem:trapezoid}   
	\end{tabular}
\end{table} 

\begin{table}[t]
    \tikzset{every picture/.style={line width=0.75pt}} 
	\begin{tabular}{c|c|c|c|c}
        type & pictogram & $k$-algebra $F_v$ of algebraic lifts $C_v$ & $\Wel \in \GW(F_v) $ & reference \\
		\hline
        
        \fourValentVertexWithMergedEdge{}{} &
        \begin{tikzpicture}[baseline={([yshift=-.5em]current bounding box.center)}]
            \path[draw] (0,0) --+ (180:0.8);
			\path[draw] (0,0) --+ (270:1.2);
			\path[draw] (0,0) --+ (45:0.8);
            \path[draw] (-0.8, -0.05) -- (0.8, -0.05);
            \path[draw] (-0.05, -1.2) -- (-0.05, 0);
            
            \fatpoint{-0.5, -0.025}
            \draw (-0.5, 0.1) node[anchor = south east] {\quadextension{$d'$}};
            \fatpoint{-0.025, -0.5}
            \draw (0.1, -0.5) node[anchor = west] {\quadextension{$d$}};
            \draw (-0.15, -1.1) node[anchor = east] {\edgeweight{$m_1$}};
            \draw (0.1, -1.1) node[anchor = west] {\edgeweight{$m_2$}};
        \end{tikzpicture} 
        %
        & $\begin{cases}
            \quotient{k[\iota, \iota']}{\big(\iota^2 - d, (\iota')^2 - d'\big)} & \text{if } m_1 \neq m_2, \\
            \quotient{k[\iota']}{\big((\iota')^2 - d'\big)} & \text{if } m_1 = m_2.
        \end{cases}$ & $\gw{1}$ & \cref{lem:trapezoid}   \\[\extraheight]  
        
        \triangleWithTwoMergedEdges{}{}{} &
        \begin{tikzpicture}[baseline={([yshift=-.5em]current bounding box.center)}]
            \path[draw] (0,0) --+ (180:0.8);
			\path[draw] (0,0) --+ (270:1.2);
			\path[draw] (0,0) --+ (45:0.8);
            \path[draw] (-0.8, -0.05) -- (0.8, -0.05);
            \path[draw] (-0.05, -1.2) -- (-0.05, 0.8);
            
            \fatpoint{-0.5, -0.025}
            \draw (-0.5, 0.1) node[anchor = south east] {\quadextension{$d_1$}};
            \fatpoint{-0.025, -0.5}
            \draw (0.1, -0.5) node[anchor = west] {\quadextension{$d_2$}};
            \draw (-0.15, -1.1) node[anchor = east] {\edgeweight{$m_1$}};
            \draw (0.1, -1.1) node[anchor = west] {\edgeweight{$m_2$}};
        \end{tikzpicture} 
        & $\quotient{k[\iota_1, \iota_2]}{\big(\iota_1^2 - d_1, \iota_2^2 - d_2 \big)}$ & $\gw{1}$ & \cref{lem:trapezoid} \\[\extraheight]

        \allDoubleVertex{} &
        \begin{tikzpicture}[baseline={([yshift=-.5em]current bounding box.center)}]
			\path[draw, double] (0,0) --+ (180:0.8);
			\path[draw, double] (0,0) --+ (270:1.2);
			\path[draw, double] (0,0) --+ (45:0.8);
            \draw (-0.1, -1.1) node[anchor = east] {\edgeweight{$m$}};
            \draw (0.1, -1.1) node[anchor = west] {\edgeweight{$m$}};
   
			\fatpoint{-0.5, 0}
            \draw (-0.5, 0.1) node[anchor = south] {\quadextension{$d_2$}};
			\fatpoint{0, -0.5}
            \draw (0.1, -0.5) node[anchor = west] {\quadextension{$d_1$}};
		\end{tikzpicture}
        & $\quotient{k[\iota]}{(\iota^2 - d_1d_2)}$ & $\begin{cases}
    		    \gw{d_1} & m \text{ odd} \\
                \gw{1} & m \text{ even}
    		\end{cases}$ & \cref{lem:fatvertex} \\[\extraheight]

        \mergedTriangles{} &
        \begin{tikzpicture}[baseline={([yshift=-.5em]current bounding box.center)}]
            \path[draw] (0,0) --+ (180:0.8);
			\path[draw] (0,0) --+ (270:1.2);
			\path[draw] (0,0) --+ (25:0.8);
            \path[draw] (0,0) --+ (65:0.8);
            \path[draw] (-0.8, -0.05) -- (0, -0.05);
            \path[draw] (-0.05, -1.2) -- (-0.05, 0);
            
            \fatpoint{-0.5, -0.025}
            \draw (-0.5, 0.1) node[anchor = south east] {\quadextension{$d_1$}};
            \fatpoint{-0.025, -0.5}
            \draw (0.1, -0.5) node[anchor = west] {\quadextension{$d_2$}};
            \draw (-0.15, -1.1) node[anchor = east] {\edgeweight{$m_1$}};
            \draw (0.1, -1.1) node[anchor = west] {\edgeweight{$m_2$}};
        \end{tikzpicture} 
        & $\quotient{k[\iota_1, \iota_2]}{\big(\iota_1^2 - d_1, \iota_2^2 - d_2 \big)}$ & $\gw{1}$ & \cref{lem:trapezoid}
    \end{tabular}
    \caption{\cref{tab:possibilities} continued.}
    \label{tab:possibilities2}
\end{table}

\subsection{Situations dual to a triangle}
\label{sec:triangle}
We work on the cases where the vertex $v$ is dual to a triangle $\Delta_v$ in the Newton subdivision of $\Delta$. Throughout we make use of \cite[Lemma~5.3]{JPP23}, which allows us to assume that $\Delta_v$ has its vertices at $(0, 0)$, $(0, m)$ (or $(0, m_1 + m_2)$), and $(1, 0)$ (or $(2, 0)$) depending on which of the edges incident to $v$ are double. 

We start with the base case of two simple point conditions adjacent to a trivalent vertex. This is in fact a special case of \cite[Lemmas~5.5 and~5.6]{JPP23}. However, due to the assumptions on the tropical point configuration being vertically stretched, the proof simplifies dramatically and we give it here for the convenience of the reader.

\begin{lemma} \label{lem:thin_triangle}
    The $k$-algebra $F_v$ of algebraic lifts of 
    
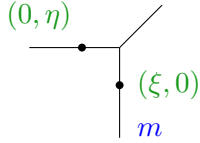
\begin{center}
        \begin{tikzpicture}[baseline={([yshift=-.5em]current bounding box.center)}]
            \path[draw] (0,0) --+ (180:1.2);
            \path[draw] (0,0) --+ (270:1.2);
            \draw (0.1, -1.1) node[anchor = west] {\edgeweight{$m$}};
            \path[draw] (0,0) --+ (45:0.8);
            
            \thinpoint{-0.5, 0}
            \draw (-0.5, 0.1) node[anchor = south east] {\coordinates{$(0, \eta)$}};
            \thinpoint{0, -0.5}
            \draw (0.1, -0.5) node[anchor = west] {\coordinates{$(\xi, 0)$}};
        \end{tikzpicture} 
    \end{center}
    is $k$ and the quadratic weight of the unique curve is $\gw{1} \in \GW(k)$. 
\end{lemma}

\begin{proof}
    By the assumption on the point conditions being vertically stretched, $\Delta_v$ has height 1 and after applying a unimodular transformation we may assume without loss of generality that 
    \[ 
        \Delta_v = \conv \big\{ (0,0), (0, 1), (m, 0) \big\}.
    \]
    In order to determine the algebraic lifts $V(f)$ for some polynomial $f$ with Newton polytope $\Delta_v$ we work with the $k$-rational point conditions on the toric boundary. Since $f$ has to vanish with order $m$ at $(\xi, 0)$ and $f$ has to vanish at $(0, \eta)$ we see, that up to scaling, $f$ is uniquely determined to be
    \[ 
        f(x, y) = y - \frac{\eta}{(-\xi)^m} (x - \xi)^m.
    \]
    In particular, we see that $f$ is defined over $k$. So the $k$-algebra of possible algebraic lifts is just $k$. Moreover, since $\Delta_v$ has no interior lattice points, $V(f)$ has no nodes and hence the quadratic weight of $V(f)$ is $\gw{1} \in \GW(k)$.
\end{proof}

\begin{lemma} \label{lem:fatpoint_thintriangle}
    Consider a vertex $v$ of type \fatPointOnVertex{}:
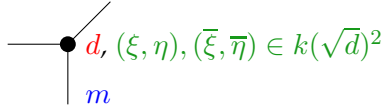
    \begin{center}
        \begin{tikzpicture}[baseline={([yshift=-.5em]current bounding box.center)}]
            \path[draw] (0,0) --+ (180:0.8);
            \path[draw] (0,0) --+ (270:0.8);
            \draw (0.1, -0.7) node[anchor = west] {\edgeweight{$m$}};
            \path[draw] (0,0) --+ (45:0.8);
            
            \fatpoint{0, 0}
            \draw (0.1, 0) node [anchor = west] {\quadextension{$d$}, \coordinates{$(\xi, \eta), (\conj \xi, \conj \eta) \in k(\sqrt{d})^2$}};
        \end{tikzpicture}
    \end{center}
    Let $L\coloneq k[\iota, \mu]/\left(\iota^2 - d, \mu^m-\eta / \conj{\eta}\right)$ and define an involution of $k$-algebras $\phi : L \to L$ by 
    \[
    \begin{array}{cccc}
    \phi :  & L         &\longrightarrow    & L \\
            &\iota   &\longmapsto        & -\iota \\
            &\mu &\longmapsto   &\frac{\conj\eta}{\eta} \mu^{m-1}.\\
    \end{array}
    \] 
    Then the algebra $F_v$ defined by the defining equations of $C_{v}$ is the algebra $L^\phi$ consisting of the fixed points of $\phi$.     
    Furthermore, any algebraic lift has quadratic weight $\gw{1} \in \GW(F_v)$.
\end{lemma}

\begin{proof}
    Again we use the assumption on the point configuration in $\RR^2$ being vertically stretched to conclude the that dual triangle to $v$ is
    \[ \Delta_v = \conv \big\{ (0,0), (0, 1), (m, 0) \big\}. \]
    Hence, the defining polynomial of any algebraic lift is of the form
    \[ f(x, y) = y - a(x - b)^m \]
    for some parameters $a$ and $b$ which we determine using the point condition $f(\xi, \eta) = f(\conj \xi, \conj \eta) = 0$. 
    Divide the equations stemming from the point conditions to obtain
    \[
        \frac{\eta}{\conj \eta} = \left( \frac{\xi - b}{\conj \xi - b} \right)^m .    
    \]
    We get that the solutions for $a$ and $b$ are given by
    \[  
        b = \frac{\xi-\mu \conj \xi }{ 1-\mu}, \qquad
        a=\conj\eta\left(\frac{1-\mu }{\conj \xi-\xi}\right)^m,
    \]
    where $\mu$ is an $m$-th root of $\eta / \conj \eta$. 
    Observe that $\phi(\mu)=\mu^{-1}$, and thus $\phi(a)=a, \phi(b)=b$.
    Hence, the algebra $k[a,b]$ defined by $a$ and $b$ is a subset of the fix point algebra $L^{\phi}$ consisting of fixed points of the involution $\phi$. Since both of these algebras have dimension $m$ over $k$, we have that $k[a,b]=L^{\phi}$. 

    Finally, any algebraic lift $C_v = V(f)$ has quadratic weight $\gw{1}$ because the Newton polygon of $f$ has no interior lattice points and hence $C_v$ has no nodes.
\end{proof}

\begin{lemma}
\label{lemma:algebratypeA}
    Let $L^{\phi}$ be the algebra consisting of fixed points of the involution $\phi$ as in Lemma~\ref{lem:fatpoint_thintriangle}. If~$m$ is odd, then we have that
    \[
    \Tr_{L^\phi/k}\big(\Wel(C_v)\big) = \Tr_{L^{\phi}/k}\qqinv{1}=\qinv{m}+\frac{m-1}{2}\qinv{2m}+\frac{m-1}{2}\qinv{-2dm}.
    \]
\end{lemma}

\begin{proof}
    Put $e_0\coloneq2, e_i\coloneq \mu^i+\mu^{-i}, f_i\coloneq\sqrt{d}\left(\mu^i-\mu^{-i}\right)$, for $i\in\ZZ$, and $m'\coloneq\frac{m-1}{2}$.
     These elements of~$L^{\phi}$ satisfy the following relations:
    \begin{enumerate}
        \item $e_i e_j=e_{i+j}+e_{i-j}$
        \item $f_i f_j=de_{i+j}-de_{i-j}$
        \item $e_i f_j=f_{i+j}+f_{j-i}$
    \end{enumerate}
    The set $\mathcal{B}\coloneq\{e_0,e_i,f_i\}_{i=1,\ldots, m'}$ is a $k$-vector space basis of $L^{\phi}$. From relations (1), (2) and (3) we compute the algebraic traces of the basis elements to be $\tr_{L^{\phi}/k}(e_0)=2m$ and $\tr_{L^{\phi}/k}(e_i)=\tr_{L^{\phi}/k}(f_j) = 0$, for $i=1,\ldots,m-1$ and $j=1,\ldots,m-1$.

    This yields the Gram matrix of the representative of the quadratic form $\Tr_{L^{\phi}/k} \qinv{1}$ determined by the chosen basis. Its $(i,j)$-entry is given by  
\[
\tr_{L^{\phi}/k}(b_i \cdot b_j),
\]  
where the elements $b_l$ of $\mathcal{B}$ are given by $b_l = e_i$ and $b_l = f_j$ for indices  
$i = 0, \ldots, m' \quad \text{and} \quad j = 1, \ldots, m'$.
One can use relations (1), (2), and (3) above to compute the entries of the Gram matrix.

    \begin{equation*}
        \tr(e_ie_j)=\begin{cases}
            4m&i=j=0,\\2m&i=j=1,\ldots,m',\\0&i\neq j,
        \end{cases} \qquad
        \tr(e_if_j)=0, \qquad
        \tr(f_if_j)=\begin{cases}
            -2dm&i=j=1,\ldots,m',\\0&i\neq j.
        \end{cases}
    \end{equation*}
    Now it is easy to see that the quadratic form of the Gram matrix with these entries represents the claimed class in the Grothendieck-Witt ring.
\end{proof}

\begin{lemma} \label{lem:type_C}
    Consider a vertex $v$ of type \fourValentVertex{}:
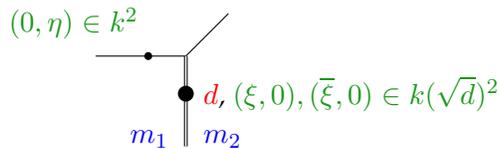
    \begin{center}
        \begin{tikzpicture}[baseline={([yshift=-.5em]current bounding box.center)}]
			\path[draw] (0,0) --+ (180:1.2);
			\path[draw, double] (0,0) --+ (270:1.2);
            \draw (-0.1, -1.1) node[anchor = east] {\edgeweight{$m_1$}};
            \draw (0.1, -1.1) node[anchor = west] {\edgeweight{$m_2$}};
			\path[draw] (0,0) --+ (45:0.8);
			
			\thinpoint{-0.5, 0}
            \draw (-0.5, 0.1) node[anchor = south east] {\coordinates{$(0, \eta) \in k^2$}};
			\fatpoint{0, -0.5}
            \draw (0.1, -0.5) node[anchor = west] {\quadextension{$d$}, \coordinates{$(\xi, 0), (\conj{\xi}, 0) \in k(\sqrt{d})^2$}};
		\end{tikzpicture}
    \end{center}
    Then there is a unique algebraic lift $C_v$ defining the $k$-algebra
    \[
        F_v = \begin{cases}
            \quotient{k[\iota]}{(\iota^2 - d)} & \text{if } m_1 \neq m_2 \\
            k & \text{if } m_1 = m_2
        \end{cases}
    \]
    and the quadratic weight of $C_v$ is $\gw{1} \in \GW(F_v)$.
\end{lemma}

\begin{proof}
    Assume the coordinates of the simple point are $(0, \eta)$ with $\eta \in k$ and the coordinates of the two $k(\sqrt{d})$-rational points making up the double point are $(\xi, 0)$ and $(\conj \xi, 0)$ with $\xi \in k(\sqrt{d}) \setminus k$. If $m_1=m_2$, then the unique algebraic lift $C_v$ of this tropical situation is given by the defining equation
    \[
        f(x,y) = y - a(x - \xi)^{m_1} (x - \conj \xi)^{m_2} \qquad \text{with} \qquad a = \frac{\eta}{(-\xi)^{m_1} (-\conj \xi)^{m_2}}.
    \]
    This curve is defined over $k$. Otherwise, if $m_1\neq m_2$, then $f$ and its conjugate are the unique algebraic lifts, and they are defined over $\quotient{k[\iota]}{(\iota^2 - d)}$.
    Finally, the quadratic weight of $C_v$ is $\Wel (C_v) = \gw{1} \in \GW(F_v)$ since the Newton polygon of $f$ has no interior lattice points and hence $C_v$ has no nodes. Note that here we are using crucially that the double edge is vertical, which is always the case, see \cref{prop-vertextypesinfloordecomposedcurve}.
\end{proof}

The following lemma will be helpful in future computations. 
\begin{lemma}\label{lm:produnits}\label{lm:fullproduct}
    Let $\zeta_m$ be a primitive $m$-th root of unity. Then, 
    \[\displaystyle\prod_{i=1}^m \zeta_m^i=(-1)^{m+1},\qquad \text{ and }\qquad
    \displaystyle\prod_{i=1}^m (1-a\zeta_m^i) =1-a^m.\]
\end{lemma}

\begin{proof}
    We have that \[p_a(x)\coloneqq \prod_{i=1}^m (x-a\zeta_m^i)=x^m-a^m.\]
    Hence, $\displaystyle\prod_{i=1}^m \zeta_m^i=(-1)^{m}p_1(0)=(-1)^{m+1}$, and $\displaystyle\prod_{i=1}^m (1-a\zeta_m^i)=p_a(1)=1-a^m$.
\end{proof}

For the next lemma, recall that a vertex of type \allDoubleVertex{} only occurs if the two merged vertical edges have the same weight. This is because all other edges have weight 1. The merging of two trivalent vertices with different weights on the vertical edges yields a vertex of type \mergedTriangles{}.

\begin{lemma} \label{lem:fatvertex}
    Consider a vertex $v$ of type \allDoubleVertex{} fixed by two double points:
    \begin{center}
        \begin{tikzpicture}[baseline={([yshift=-.5em]current bounding box.center)}]
			\path[draw, double] (0,0) --+ (180:1.2);
			\path[draw, double] (0,0) --+ (270:1.2);
			\path[draw, double] (0,0) --+ (45:0.8);
            \draw (-0.1, -1.1) node[anchor = east] {\edgeweight{$m$}};
            \draw (0.1, -1.1) node[anchor = west] {\edgeweight{$m$}};
   
			\fatpoint{-0.5, 0}
            \draw (-0.5, 0.1) node[anchor = south east] {\quadextension{$d_2$}, \coordinates{$(0, \eta), (0, \conj{\eta}) \in k(\sqrt{d_2})^2$}};
			\fatpoint{0, -0.5}
            \draw (0.1, -0.5) node[anchor = west] {\quadextension{$d_1$}, \coordinates{$(\xi, 0), (\conj{\xi}, 0) \in k(\sqrt{d_1})^2$}};
		\end{tikzpicture}
    \end{center}
    Then the $k$-algebra of algebraic lifts is given by
    \[
        F_v = \quotient{k[\iota]}{(\iota^2 - d_1d_2)}
    \]
    and the quadratic weight of $C_v$ is $\gw{d_1} \in \GW(F_v)$ if $m$ is odd and $\gw{1} \in \GW(F_v)$ if $m$ is even.
\end{lemma}

\begin{proof}
    Since we are working with parametrized tropical curves, it is clear that an algebraic lift of this tropical vertex will consist of two conjugate curves $C_1$, $C_2$, defined by equations $f_1$, $f_2$, respectively such that the Newton polygon of each $f_i$ is $\conv \big\{ (0,0), (0, 1), (m, 0) \big\}$. Without loss of generality $C_1$ contains $(0, \eta)$ and $(\xi, 0)$ and $C_2$ contains the conjugate points. There is only one option:
    \begin{equation} \label{eq:options_fat_vertex}
        \begin{aligned}
            f_1(x,y) &= y - a (x - \xi)^m  & &\text{with } a = \frac{\eta}{(-\xi)^m} \\
            f_2(x,y) &= y - b (x - \conj\xi)^m & &\text{with } b = \frac{\conj\eta}{(-\conj\xi)^m}.
        \end{aligned}
    \end{equation}
    Here, we are abusing notation and denote by $\conj{(\cdot)}$ conjugation with respect to the field extensions $k(\sqrt{d_1}) / k$ and $k(\sqrt{d_2}) / k$, depending on context. To make this clearer, we define the following $k$-algebra automorphism of $L = k[\iota_1, \iota_2]/(\iota_1^2 - d_1, \iota_2^2 - d_2)$ extending both conjugations
    \[
        \phi \colon L \longrightarrow L \colon
        \iota_i \longmapsto -\iota_i, \qquad \text{for } i = 1,2.
    \]
    With this notation we see that $b = \phi(a)$ and $f_2 = \phi(f_1)$ and hence the combined curve $C_1 \cup C_2$ defined by the product $f_1f_2$ is in fact defined over 
    \[F_v = \{\alpha \in L \mid \phi(\alpha) = \alpha \} = \quotient{k[\iota]}{(\iota^2 - d_1d_2)}. \] 
    
    Now we compute the quadratic weight. To this end note that the irreducible components $C_i$ have no nodes themselves, simply because the Newton polygons of their defining equations have no interior lattice points. Hence, all nodes of $C_1 \cup C_2$ arise from the intersection of the two branches. 
    These points have $x$-coordinate determined by the equation
    \begin{equation} \label{eq:fatvertex_findnodes}
        a (x - \xi)^m = (\phi a) (x - \phi\xi)^m
    \end{equation}
    and any solution of this equation determines a unique corresponding $y$-coordinate via $f_1(x,y) = f_2(x, y) = 0$. In order to solve \cref{eq:fatvertex_findnodes}, we choose $\mu_i$ to be an $m$-th root of $(\phi a)/a$ and solve for $x$:
    \[
        \mu = \frac{x- \xi}{x- \phi\xi} \qquad \Longleftrightarrow \qquad x = \frac{\xi - \mu \cdot \phi\xi}{1 - \mu}.
    \]
    For later reference we determine the corresponding $y$-coordinate:
    \begin{equation} \label{eq:fatvertex_y}
        y = (\phi a) \Big(\frac{\xi - \mu \cdot \phi\xi}{1 - \mu} - \phi\xi\Big)^m = (\phi a) \Big(\frac{\xi - \phi \xi}{1 - \mu} \Big)^m.
    \end{equation}
    This determines all nodes of $C_1 \cup C_2$ and we now compute the Hessian of the defining equation
    \[ \Hess_{f_1f_2} = \begin{pmatrix}
        (f_1f_2)_{xx} & (f_1f_2)_{xy} \\
        (f_1f_2)_{xy} & (f_1f_2)_{yy}
    \end{pmatrix} = \begin{pmatrix}
        2(f_1)_x(f_2)_x + f_1(f_2)_{xx} + f_2(f_1)_{xx} & (f_1)_x + (f_2)_x \\
        (f_1)_x + (f_2)_x & 2
    \end{pmatrix}. \]
    Note that if we plug in the coordinates $(x_0, y_0)$ of a node, then $f_1(x_0, y_0) = f_2(x_0, y_0) = 0$ and hence 
    \begin{equation} \label{eq:fatvertex_detHess}
        \begin{aligned}
            -\det \Hess_{f_1f_2}(x_0, y_0) &= \big( (f_1)_x(x_0, y_0) - (f_2)_x(x_0, y_0) \big)^2 \\
            &= \Big(\frac{my_0}{x_0 - \xi} - \frac{my_0}{x_0 - \phi\xi}\Big)^2 \\
            &= m^2 y_0^2 \left(\frac{\xi - \phi\xi}{ \mu \frac{\xi - \phi \xi}{ 1 - \mu} \frac{\xi - \phi\xi}{ 1- \mu}} \right)^2 \\
            &= m^2 y_0^2 \frac{1}{(\xi - \phi\xi)^2} \frac{(1- \mu)^4}{\mu^2}.
        \end{aligned}
    \end{equation}
    By definition, the quadratic weight of $C_1 \cup C_2$ is now given as the product of \cref{eq:fatvertex_detHess} taken over all nodes and then considered as an element in $\GW(F_v)$. 
    This can be easily computed using the representation of $y_0$ which we provided in \cref{eq:fatvertex_y} and using the formulas established in Lemma~\ref{lm:produnits}. 
    To see this, recall that the different choices for the nodes are given by
    \[
    x_i = \frac{\xi - \mu \cdot \phi\xi}{1 - \mu_i},\qquad y_i = (\phi a) \Big(\frac{\xi - \phi \xi}{1 - \mu_i} \Big)^m.
    \]
    The result is
    \begin{equation} \label{eq:fatvertex_Wel}
        \begin{aligned}
            \prod_{(x_0, y_0) \text{ node}} \Big( (f_1)_x(x_0, y_0) - (f_2)_x(x_0, y_0) \Big)^2 
            &= \prod_{i = 0}^{m-1} m^2 y_0^2 \frac{1}{(\xi - \phi\ xi)^2} \frac{(1- \mu\zeta_m^i)^4}{\mu^2 (\zeta_m^i)^2 } \\
            &= m^{2m} \frac{ (\phi a)^{2m} (\xi - \phi\xi)^{2m^2}}{\big(1 - (\phi a)/a \big)^{2m}} \frac{1}{(\xi - \phi\xi)^{2m}} \frac{\big(1- (\phi a)/a \big)^4}{\big( (\phi a)/a \big)^2} .
        \end{aligned}
    \end{equation}
    We claim that in $\GW(F_v)$ we have
    \[ 
        \Big \langle\text{right hand side of \cref{eq:fatvertex_Wel}}\Big\rangle = \begin{cases}
            \gw{d_1} & \text{if $m$ is odd} \\
            \gw{1} & \text{if $m$ is even.}
        \end{cases} 
    \]
    To show this, we exhibit most of the terms appearing in the right hand side of \cref{eq:fatvertex_Wel} as squares of elements in $F_v$. 
    \begin{itemize}
        \item Clearly $m \in F$ and hence $m^{2m}$ is a square in $F_v$.

        \item  From the computation
        \[ 
            \frac{(\phi a)^{2m}}{\big( 1 - (\phi a)/a \big)^{2m}} = \Big( \frac{a (\phi a)}{a - \phi a} \Big)^{2m}
            = \big(\underbrace{a (\phi a)}_{\in F} \big)^{2m} \cdot \Big(\underbrace{\frac{1}{(a - \phi a)\sqrt{d_1}}}_{\in F}\Big)^{2m} d_1^m
        \]
        we see that up to multiplication by squares in $F_v$ this is equal to $d_1^m$, which in turn is equal to $d_1$ if $m$ is odd and 1 if $m$ is even.
        
        \item $\xi - \phi\xi \in \sqrt{d_1} \cdot k$ and hence $(\xi-\phi\xi)^2 \in k \subseteq F_v$. The total power of $\xi - \phi\xi$ in \cref{eq:fatvertex_Wel} is $2m(m-1)$, which is divisible by 4. Therefore, all occurrences of $\xi - \phi\xi$ combined are a square in $F_v$.

        \item $(\phi a)/a$ has the property that $\phi \big( (\phi a)/a \big) = a / (\phi a)$. Therefore 
        \[ \frac{\big(1 - (\phi a)/a \big)^4}{\big( (\phi a)/a \big)^2} = \left(\frac{1 - (\phi a)/a}{(\phi a)/a} \right)^2 \big(1 - (\phi a)/a \big)^2 = \Big[ \underbrace{\phi \big(1-  (\phi a)/a \big) \cdot \big(1 - (\phi a)/a\big)}_{\in F} \Big]^2 \]
        is a square in $F_v$.
    \end{itemize}
    Combined, these observations prove the claim. 
\end{proof}

\begin{lemma}
    If \[F\coloneqq\frac{k[\iota]}{(\iota^2-d_1d_2)},\]
    then \[ \Tr_{F/k}{\qinv{d_1}}=\qinv{2d_1}+\qinv{2d_2}.\]
\end{lemma}

\begin{proof}
    Since the morphism $\Tr_{F/k}$ is $\GW(k)$-linear, we have that
    \begin{align*}
        \Tr_{F/k}{\qinv{d_1}}
        &=\qinv{d_1}\cdot \Tr_{F/k}{\qinv{1}}\\
        &=\qinv{d_1}\cdot \big(\qinv{2}+\qinv{2d_1d_2} \big)\\
        &=\qinv{2d_1}+\qinv{2d_2},
    \end{align*}
    since $\qinv{2d_1^2d_2}=\qinv{2d_2}\in\GW(k)$.
\end{proof}

\subsection{Parallelogram cases}
We deal with vertices $v$ which are dual to a parallelogram in the dual subdivision. More precisely, the relevant cases are the $k$-rational situation and the types \fatPointOnParallelogram{}, \parallelogramWithOneDoubleEdge{}, and \parallelogramWithTwoDoubleEdges{}. 
Throughout this subsection we can assume without loss of generality that the lattice lengths of the parallelogram~$\Delta_v$ are $1$ or $2$. 
Indeed, suppose the lattice length was larger, i.e.\ the defining equation of the curve contains a factor of the form $(x-a)^n$. This local equation corresponds to a singularity that is not nodal. However, its local contribution is given by the contribution of a perturbation. A perturbation can be used to split this high multiplicity branch into $n$ branches of multiplicity $1$ each.
We will see that in all parallelogram cases, the nodes in the perturbations are always \emph{split}, meaning that their mass is $\gw{1}$, over the field of definition of the algebraic curves that lift the combinatorial configuration around the parallelogram.  
To prove this, we can assume that all powers are equal to $1$. It then follows that for curves where some powers are greater than $1$, the quadratic weight is $\gw{1}$.

The product $(x-a)^{m_1}(x-\bar{a})^{m_2}$, for $a\in k(\sqrt{d})\setminus k$, is defined over $k$ if and only if $m_1=m_2$. Hence, we have that for $m_1\neq m_2$, our deformation can treat the branches independently and so we can restrict ourselves to study the case of a single branch $x-a$, or double conjugated branches~$(x-a)(x-\bar{a})$.
Consequently, in \cref{lem:par3} and \cref{lem:par4} we assume $m_1=m_2$. For $m_1\neq m_2$ we use \cref{lem:par1} twice.

Using a lattice preserving transformation we take $\Delta_v$ to have three of its corners given by $(0,0), (1, 0), (p, q)$ with $\gcd(p, q) =1$, or $(0,0), (2, 0), (2p, 2q)$, depending on the combinatorial type of $v$.

The following lemma is in fact the same as \cite[Proposition~5.7]{JPP23}. Nevertheless, we repeat the proof to introduce the reader to the computational strategy used throughout this subsection. 

\begin{lemma} \label{lem:par1}
    Consider the $k$-rational parallelogram:
    
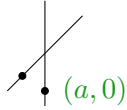
\begin{center}
        \begin{tikzpicture}[baseline={([yshift=-.5em]current bounding box.center)}]
			\path[draw] (0, -0.7) -- (0, 0.7);
			\path[draw] (-0.5, -0.5) -- (0.5, 0.5);
			\thinpoint{0, -0.5}
            \draw (0.1, -0.5) node[anchor = west] {\coordinates{$(a, 0)$}};
			\thinpoint{-0.3, -0.3}
		\end{tikzpicture}
    \end{center}
	Without loss of generality this is dual to $\Delta_v = \conv \big\{ (0,0), (1, 0), (p, q), (p+1, q) \big\}$ with $\gcd(p, q) =1$. This situation has a unique algebraic lift and its quadratic weight is the factor $\gw{1} \in \GW(k)$.
\end{lemma}

\begin{proof}
	In this case the algebraic curve $C$ is given by the equation $f(x, y) = (x-a)(x^py^q-b)$ for some $a, b \in k$, which are uniquely determined by the simple point conditions. We compute the internal nodes of $C$, i.e. the nodes with coordinates $x \neq 0$ and $y \neq 0$. Consider the derivatives of $f$
	\begin{align*}
		f_x &= (x^py^q - b) + (x-a)px^{p-1} y^q \\
		f_y &= (x-a) qx^py^{q-1}
	\end{align*}
	and set $f_y= 0$. Since $x,y\neq 0$ we obtain $x = a$ which we plug into $f_x = 0$ to obtain $y^q = b/a^p$. In particular, this shows that there are $q$ nodes $\eta_1,\ldots,\eta_q$ in total. We spell out the $y$-coordinates of $\eta_i$ as follows. Let $y_0$ be a $q$-th root of $b/a^p$ and let $\zeta_q$ be a primitive $q$-th root of unity. Then $\eta_i = (a, \zeta_q^i y_0)$ for $i = 0, \ldots, q-1$.
	
	We compute the quadratic weight of $C$. For this note that 
	\[ f_{yy} = (x-a) q(q-1)x^py^{q-2} = 0, \]
    if we plug in $x=a$.
	So in order to evaluate $-\det \Hess_f$ at a node $\eta_i$, is suffices to compute $-\big(-f_{xy}^2(\eta_i)\big)$. Evaluating the mixed second 
    derivative
    \[ f_{xy} = qx^py^{q-1} + (x-a)pqx^{p-1}y^{q-1}\]
    at $\eta_i=(a, \zeta_q^i y_0)$ we obtain
	\[f_{xy}(\eta_i)=qa^p\zeta_q^{i(q-1)}y_0^{q-1}=\frac{qb}{y_0\zeta_q^i} \]
	and hence 
	\[\Wel(C) = \gw{ \prod_{i = 0}^{q-1} \frac{q^2b^2}{ y_0^2 \zeta_q^{2i}} } = \gw{ \left(\frac{a^p}{b} \right)^2 \frac{1}{\zeta_q^{2\frac{q(q-1)}{2}}} } = \gw{1}.\]
	We remark that this computes the quadratic weight irrespective of whether $\zeta_q \in k$ or not. If $\zeta_q \in k$, then this is the product over $q$ many $k$-rational nodes of $C$. If $\zeta_q \not\in k$, then this is the product of norms as in Remark \ref{remark:computingWel}. 
\end{proof}

\begin{lemma} \label{lem:par3}
    We consider a situation of type \parallelogramWithOneDoubleEdge{}:
    \begin{center}
        \begin{tikzpicture}[baseline={([yshift=-.5em]current bounding box.center)}]
			\path[draw] (0, -0.7) -- (0, 0.7);
			\path[draw, double] (-0.5, -0.5) -- (0.5, 0.5);
			\thinpoint{0, -0.5}
            \draw (0.1, -0.5) node[anchor = west] {\coordinates{$(a, 0)$}};
			\fatpoint{-0.3, -0.3}
            \draw (-0.3, -0.2) node[anchor = south east] {\quadextension{$d$}};
		\end{tikzpicture}
    \end{center}    
	We assume that three of the corners of $\Delta_v$ are $(0,0), (1, 0), (2p, 2q)$ with $\gcd(p, q) =1$ and the simple and double point conditions are as depicted. Then there is a unique algebraic lift which is defined over $k$ and its quadratic lift is $\langle 1 \rangle \in \GW(k)$.
\end{lemma}

\begin{proof}
	Any algebraic curve lift $C$ of this situation is given by an equation of the form
	\[ f(x, y) = (x-a) (x^py^q - b)(x^py^q - \conj{b}) \]
	for some $a \in k$ and $b \in k(\sqrt{d}) \setminus k$. The $a$ and $b$ are uniquely determined by the point conditions. The polynomial factors as above because the parametrization given by the stable map  encodes this splitting. 
	
	\textbf{Step 1:} We determine all interior nodes of $C$, i.e. nodes away from the toric boundary at $x = 0$ or $y = 0$. We compute the partial derivatives of $f$
	\begin{align*}
		f_x &= (x^py^q - b) (x^py^q - \conj{b}) + (x-a) \Big[px^{p-1}y^q (x^py^q - \conj{b}) + (x^py^q-b)px^{p-1}y^q \Big] \\
		f_y &= (x-a)  qx^py^{q-1}(x^py^q-\conj{b}) + (x-a)(x^py^q-b)qx^py^{q-1}
	\end{align*}
	and solve for $f_x = f_y = 0$. The equation $f_y = 0$ with $x \neq 0$ and $y \neq 0$ holds precisely for $x = a$ or $x^py^q = (\conj{b} + b)/2$. If $x = a$, then $f_x = 0$ implies that $y^q = b/a^p$ or $y^q = \conj{b}/a^q$. In the other case, meaning $x^py^q = (\conj{b} + b)/2$ we obtain no further solutions since $f_x = 0$ translates into
	\begin{align*}
		&\left( \frac{\conj{b} + b}{2} -b \right)\left( \frac{\conj{b} + b}{2} - \conj{b} \right) + (x-a) \frac{p}{x} \frac{\conj{b} + b}{2} \left( \frac{\conj{b} + b}{2} - \conj{b} \right) + (x-a) \left( \frac{\conj{b} + b}{2} - b \right) \frac{p}{x} \frac{\conj{b} + b}{2} = 0 \\
		&\Longleftrightarrow (\conj{b} - b) (b - \conj{b}) + p(\conj{b} + b)(b - \conj{b}) + p(\conj{b} - b)(\conj{b} + b) = \frac{ap}{x} \underbrace{\Big[ (\conj{b} + b)(b - \conj{b}) + (\conj{b} - b)(\conj{b} + b) \Big]}_{= 0}\\
        &\Longleftrightarrow(\overline{b}-b)(b-\overline{b})=0
	\end{align*}
	which is a contradiction.
	In summary, we found $2q$ many nodes with $x$-coordinate $x=a$ and $y$ coordinate the solutions to $y^q=b/a^p$ and $y^q=\overline{b}/a^p$.
    
	\textbf{Step 2:} In order to compute quadratic weight of $C$, we proceed to compute the second derivatives of $f$. We find that
	\begin{equation*}
		f_{yy} = (x-a) qx^p \Big[ (q-1) y^{q-2} (2x^py^q - \conj{b} - b) + y^{q-1} 2qx^py^{q-1} \Big]
	\end{equation*}
	is equal to 0 at every node because $x = a$ at every node. Therefore, $-\det \Hess_f = -f_{xy}^2$ at every node. We compute
	\begin{equation*}
		f_{xy} = (x^py^q - b) qx^py^{q-1} + qx^py^{q-1}(x^py^q - \conj{b})
		+ (x-a)px^{p-1} \Big[ \cdots \Big]
	\end{equation*}
	and plugging in the coordinates of the nodes gives
	\[ \Wel(C) = \gw{ \prod_{i} \left( \frac{q}{y_i} b (b - \conj{b}) \right)^2 \cdot \prod_{i} \left( \frac{q}{y_i} \conj{b} (b - \conj{b}) \right)^2 } = \gw{1},\]
    since $\prod_i y_i=(b/a^p)(\overline{b}/a^p)\in k$ by \cref{lm:fullproduct}.
\end{proof}

\begin{lemma} \label{lem:par4}
    Consider a situation of type \parallelogramWithTwoDoubleEdges{}:
    \begin{center}
        \begin{tikzpicture}[baseline={([yshift=-.5em]current bounding box.center)}]
			\path[draw, double] (0, -0.7) -- (0, 0.7);
			\path[draw, double] (-0.5, -0.5) -- (0.5, 0.5);
			\fatpoint{0, -0.5}
            \draw (0.1, -0.5) node[anchor = west] {\quadextension{$d_1$}, \coordinates{$(a, 0), (\conj a, 0)$}};
			\fatpoint{-0.3, -0.3}
            \draw (-0.3, -0.2) node[anchor = south east] {\quadextension{$d_2$}};
		\end{tikzpicture}
    \end{center}
	Let $\Delta_v$ be the parallelogram spanned by $(0,0), (2, 0), (2p, 2q)$ with $\gcd(p, q) = 1$ and double point conditions on both branches. Then $F_v=k$ and the quadratic weight of the unique algebraic lift is $\gw{1}$.
\end{lemma}

\begin{proof}
	The two double point conditions may a priori live over two different quadratic field extension $k(\sqrt{d_1})$ and $k(\sqrt{d_2})$. In the following we will denote conjugation in $k(\sqrt{d_1})$ by $\conj{(\cdot)}$ and conjugation in $k(\sqrt{d_2})$ by $\conjb{(\cdot)}$.
	With this notational convention the curve $C$ is defined by an equation of the form
	\[ f(x,y) = (x - a) (x - \conj{a}) (x^py^q - b) (x^py^q - \conjb{b}) \]
	for $a \in k(\sqrt{d_1}) \setminus k$ and $b \in k(\sqrt{d_2}) \setminus k$ uniquely given by the point conditions. This curve is unique and it is defined over~$k$.
    We compute the partial derivatives of $f$ in order to find the interior nodes of $C$.
	\begin{align*}
		f_x &= (x-a) (x-\conj{a}) \Big[ (x^py^q - b) px^{p-1}y^q + px^{p-1}y^q (x^py^q - \conjb{b}) \Big]\\&\phantom{=} + \Big[ (x-a) + (x-\conj{a}) \Big] (x^py^q-b)(x^py^q-\conjb{b})  \\
		f_y &= (x-a)(x-\conj{a}) \Big[ (x^py^q-b)qx^py^{q-1} + qx^py^{q-1}(x^py^q-\conjb{b}) \Big]
	\end{align*}
	The $y$-derivative vanishes for $x = a$, $x = \conj{a}$, or $x^py^q = (b+\conjb{b})/2$. We show that the latter is not a valid option: plugging $x^py^q = (b+\conjb{b})/2$ into $f_x = 0$ leads to 
	\begin{equation*}
		0 = \underbrace{(x-a)(x-\conj{a}) \Big[ (\conjb{b}-b)\frac{p}{x}(b+\conjb{b}) + \frac{p}{x}(b + \conjb{b})(b - \conjb{b}) \Big]}_{=0} + \Big[(x-a) + (x-\conj{a}) \Big] \underbrace{(\conjb{b} - b)(b - \conjb{b})} _{\neq 0}
	\end{equation*}
	which implies $x = (a + \conj{a})/2$. But any point defined by $x_0 = (a + \conj{a})/2$ and $x_0^py_0^q = (b + \conjb{b})/2$ does not lie in $C = V(f)$:
	\[f(x_0, y_0) = \left(\frac{\conj{a} - a}{2}\right)  \left(\frac{a - \conj{a}}{2}\right)  \left(\frac{\conjb{b}-b}{2}\right)  \left(\frac{b - \conjb{b}}{2}\right) \neq 0. \]
	
	We continue the search for nodes of $C$ with $x = a$, in which case the condition $f_x = 0$ yields
	\[ 0 = \underbrace{(x-\conj{a})}_{\neq 0} (x^py^q - b)(x^py^q - \conjb{b}), \]
	i.e. $x^py^q = b$ or $x^py^q = \conjb{b}$. We find the same for $x = \conj{a}$, leading to $4q$ nodes in total.
	
	We proceed to evaluate $-\det\Hess_f$ at the nodes. First, we see that the derivative $f_{yy}$ is of the form $(x-a)(x - \conj{a})\big[ \cdots \big]$, which is clearly zero for every interior node of $C$. The mixed second derivative is equal to
	\[ f_{xy} = (x-a)(x-\conj{a}) \Big[ \cdots \Big] + \Big[ (x-a) + (x-\conj{a}) \Big] \Big[ (x^py^q-b)qx^py^{q-1} + qx^py^{q-1}(x^py^q-\conjb{b}) \Big] \]
	and thus the node $\eta = (x_0, y_0)$ with $x_0 = a$ and $x_0^py_0^q = b$ contributes a factor of
	\[ \gw{ -\left((a-\conj{a})\frac{q}{y_0}b(b - \conjb{b})\right)^2 } = \gw{-y_0^2} \]
	to $\Wel(C)$. Similar for nodes with $x_0 = \conj{a}$ and/or $x_0^py_0^q = \conjb{b}$. In total we get
	\begin{align*}
		\Wel(C) = &\gw{ \prod_i -y_i^2 } 
		\gw{ \prod_i -y_i^2 } \gw{ \prod_i -y_i^2 }
		\gw{ \prod_i -y_i^2 } = \gw{1},
	\end{align*}
    since the product of the $y$-coordinate over all nodes equals $(b/a^p)(\hat{b}/a^p)(b/\bar{a}^p)(\hat{b}/\bar{a}^p)$ by \cref{lm:fullproduct}.
\end{proof}

\begin{lemma} \label{lem:par5}
    Consider a situation of type \fatPointOnParallelogram{}:
    \begin{center}
        \begin{tikzpicture}[baseline={([yshift=-.5em]current bounding box.center)}]
			\path[draw] (0, -0.7) -- (0, 0.7);
			\path[draw] (-0.5, -0.5) -- (0.5, 0.5);
			\fatpoint{0, 0}
            \draw (0, 0) node[anchor = south east] {\quadextension{$d$}, \coordinates{$(a, \eta)$}};
		\end{tikzpicture} 
    \end{center}
	Let $\Delta_v$ be the parallelogram spanned by $(0,0), (1, 0), (p, q)$ with $\gcd(p, q) = 1$ and a double point condition at the intersection of the edges of the tropical curve. There is a unique algebraic lift which is defined over $F_v = k[\iota]/(\iota^2 - d)$ and its quadratic weight is $\gw{1} \in \GW(F_v)$.
\end{lemma}

\begin{proof}
	In this case the curve $C$ has two branches, one passing through one of the point conditions, the other through its conjugate. This means that $C$ is given by the equation 
    \[f(x, y) = (x-a)(x^py^q - b)\] for suitable $a, b \in k(\sqrt{d})$, where $a$ and $b$ cannot both live in $k$. In particular, $C$ is defined over $F_v = k[\iota]/(\iota^2 - d)$. The $a$ and $b$ are uniquely defined. From here on the computation of $\Wel(C)$ is the same as in \cref{lem:par1}, except that the computation takes place in $k[\iota] / (\iota^2 - d)$, i.e. $a^2$ and $b^2$ are indeed squares.
\end{proof}

\subsection{Remaining vertex cases}

We now deal with the remaining types \triangleWithMergedEdge{}, \fourValentVertexWithMergedEdge{}, \triangleWithTwoMergedEdges{}, and \mergedTriangles{}. 

\begin{lemma}
    \label{lem:trapezoid}
    $F_v$ and the quadratic weight of the algebraic lift $C_v$ of a vertex $v$ of one of the types \triangleWithMergedEdge{}, \fourValentVertexWithMergedEdge{}, \triangleWithTwoMergedEdges{}, or \mergedTriangles{} are as stated in Tables~\ref{tab:possibilities} and~\ref{tab:possibilities2}.
\end{lemma}

\begin{proof}
    Note that all of them are \enquote{composite} in the sense that the dual polygon of these types is a Minkowski sum of smaller pieces:
    \begin{center}
        \triangleWithMergedEdge{} = triangle + edge \\
        \fourValentVertexWithMergedEdge{} = \fourValentVertex{} + edge \\
        \triangleWithTwoMergedEdges{} = triangle + edge + edge \\
        \mergedTriangles{} = triangle + triangle.
    \end{center}
    The parametrization of the plane tropical curve tells us that the defining equation $f$ of any algebraic lift $C_v = V(f)$ must factor according to these decompositions. We detail the computations for the case \triangleWithMergedEdge{}, the other cases are analogous.
    
    In case \triangleWithMergedEdge{} the local polygon~$\Delta_v$ is a trapezoid and we call the triangle $\Delta_1$ and the edge $\Delta_2$. So the defining equation $f$ factors as $f = f_1 f_2$ with $V(f_1)$ lifting $\Delta_1$ and $V(f_2)$ lifting the edge~$\Delta_2$. Call $\sigma$ the longer edge of the two parallel edges of $\Delta_v$. We can without loss of generality assume that $\sigma=[0,m]$. Then $f\vert_{\sigma}=(x-a)^p\cdot f_1\vert_\sigma$. 
    Denote $\overline{a}$ the Galois conjugate of $a$ in $k(\sqrt{d})/k$.
	Then $f_1\vert_{\sigma}= \big(x-\overline{a} \big)^p$. Note that there are two possibilities for the choice of coefficients of $f_1\vert_\sigma$ and $f_2$ since can we can swap $a$ and $\overline{a}$.
	In particular, we have found $f_2$ and it remains to find all possibilities for the remaining coefficients of $f_1$. Now recall that $f_1$ has Newton polygon the triangle and two edges of this triangle are marked, one by $\overline{a}$ and the other one by an additional point condition. So we have reduced the problem to where $\Delta_v$ is a triangle which we already dealt with in \cref{sec:triangle}. 
	
	To compute the quadratic weight $\operatorname{Wel}^{\A^1}(C_v)$ of such a curve one needs to find all nodes and compute their quadratic type. On the one hand, we have the nodes of the curve defined by $f_1$, but there are none, since the Newton polygon of $f_1$ has no interior points because we are only considering vertically stretched point conditions. On the other hand, we have nodes arising from the intersection of the curves defined by $f_1$ and~$f_2$. The latter are all split nodes, that is, their branches live over $F_v$, one does not have to adjoin a square root and the type of these nodes is~$1\in F_v^\times/(F_v^\times)^2$.
    In particular, they contribute the factor $\gw{1}$ to $\Wel(C_v)$ and thus we get that $\Wel(C_v)= \gw{1}\in\GW(F_v)$.
    
\end{proof}

\section{Contributions from edges}\label{sec-edges}
We start with some preliminary computations in the Grothendieck Witt ring and some trace computations which we will use later in this section to compute the contributions from the edges.
\begin{lemma}
    Let $a,b,d\in k^\times$ and assume that $a^2-b^2d\neq 0$ and assume $\operatorname{char}k\neq 2$. Then the following relations hold in $\GW(k)$:
    \begin{equation}\label{eq:gwrelation}
        \gw{1}+\gw{-d}=\gw{a^2-b^2d} \big(\gw{1}+\gw{-d} \big),
    \end{equation}
    \begin{equation}\label{eq:gwrealtion2}
        2 \big(\gw{1}+\gw{-d} \big) = 2 \big(\gw{2}+\gw{-2d} \big).
    \end{equation}
    Moreover, the following non-degenerate quadratic forms are isometric and thus represent the same element in $\GW(k)$ 
    \begin{equation}\label{eq:gwrealtion3}
        \left[(x,y)\longmapsto bx^2+2axy+bdy^2\right] 
        =\left[(x,y)\longmapsto bx^2-\frac{a^2-b^2d}{b}y^2\right].
    \end{equation}
\end{lemma}

\begin{proof}
     Equations~\eqref{eq:gwrelation} and~\eqref{eq:gwrealtion2} follow from the relations in the presentation of the Grothendieck-Witt ring (cf. \cref{df:gw}) together with the computations:
        \begin{align*}
            \gw{1}+\gw{-d} 
            &=\gw{a^2}+\gw{-b^2d} \\
            &=\gw{a^2-b^2d}+\gw{a^2\cdot(-db^2)\cdot(a^2-b^2d)}\\ 
            &=\gw{a^2-b^2d}+\gw{-d(a^2-b^2d)} \\
            &=\gw{a^2-b^2d} \big (\gw{1}+\gw{-d} \big),
        \end{align*}
together with  $\gw{1}+\gw{1}=\gw{2}+\gw{2}$, and $\gw{-d}+\gw{-d}=\gw{-2d}+\gw{-2d\cdot (-d)\cdot (-d)}=\gw{-2d}+\gw{-2d}$.

     \cref{eq:gwrealtion3} follows from the following equality at the level of Gram matrices:
    \[
    \begin{bmatrix} -\frac{a}{b} & 1\\1 & 0\end{bmatrix}
    \begin{bmatrix} b & a\\a & bd\end{bmatrix}
    \begin{bmatrix} -\frac{a}{b} & 1\\1 & 0\end{bmatrix}=
    \begin{bmatrix} \frac{b^2d-a^2}{b} & 0\\ 0 & b\end{bmatrix}.
    \]
\end{proof}

\begin{lemma}\label{lm:intertrace}
    Let $F$ be a field such that $2m\in F^\times$. Let $L_\circ$ be the algebra
    \[
    L_\circ\coloneqq\frac{F[x]}{\big(x^m-(a^2-b^2d) \big)},
    \]
    where $a,b,d\in F^\times$ such that $a^2-b^2d\neq0$.
    If $m$ is even, then
    \begin{align*}
    \gw{m}\Tr_{L_\circ/F} \gw{x^\frac{m}{2}+a}=&
    \begin{cases}
        \gw{a}+\gw{ad(a^2-b^2d)}+\frac{m-2}{2}\h
        & \text{ if } m\equiv 2 \mod 4,\\
        \gw{a}+\gw{ad(a^2-b^2d)}+\gw{1}+\gw{-d}+\frac{m-4}{2}\h
        & \text{ if } m\equiv 0 \mod 4.\\
    \end{cases}\\
    \gw{m}\Tr_{L_\circ/F} \gw{x^{\frac{m}{2}+1}+ax}=&
    \begin{cases}
        \gw{1}+\gw{-d}+\frac{m-2}{2}\h 
        & \hspace{3.9cm}\text{ if } m\equiv 2 \mod 4,\\
        \frac{m}{2}\h
        & \hspace{3.9cm}\text{ if } m\equiv 0 \mod 4.\\
    \end{cases}
    \end{align*}
\end{lemma}

\begin{proof}
    Let $\mathcal{B}_\circ$ be the ordered basis $\big\{1,x^\frac{m}{2},x,x^{\frac{m}{2}+1},\dots,x^{\frac{m}{2}-1},x^{m-1}\big\}$  of $L_\circ$ as an $F$-vector space.
    It follows that the matrix with entries $m_{i,j}=\tr_{L_\circ/F}\big(x^{i-1}\cdot x^{j-1} \cdot (x^\frac{m}{2}+a) \big)$ with $i,j=1,\ldots,m$ 
    is a Gram matrix of a quadratic form
    representing
 $\Tr_{L_\circ/F} \gw{x^\frac{m}{2}+a}$. 
 Recall that the algebraic trace $\tr_{L_\circ/F}(1)=m$ and that $\tr_{L_\circ/F}(x^i)=0$ if $0<i<m$.
 So this Gram matrix is equal to
    \[
    m\cdot \begin{bmatrix}
    A & 0 & 0 & \cdots & 0 \\
    0 & 0 & 0 & \cdots & B \\
    0 & 0 & 0 & \iddots & 0 \\
    \vdots & \vdots & \iddots & \iddots & \vdots \\
    0 & B & 0 & \cdots & 0
    \end{bmatrix}
    \]
    with $2\times 2$ blocks
    \[
        A=\begin{bmatrix}
            a&(a^2-b^2d)\\(a^2-b^2d)&a(a^2-b^2d)
        \end{bmatrix},\qquad
        B=(a^2-b^2d)\cdot\begin{bmatrix}
            1&a\\a&(a^2-b^2d)
        \end{bmatrix}.
    \]
    By \cref{eq:gwrealtion3}, the matrix $A$ represents the class $\gw{a}+\gw{ad(a^2-b^2d)}$ in $\GW(k)$ and the matrix $B$ represents $\gw{a^2-b^2d}\big(\gw{1}+\gw{-d} \big)$, which equals~$\gw{1}+\gw{-d}$, by \cref{eq:gwrelation}, in $\GW(k)$.
    
    The matrix with entries $\tr_{L_\circ/F} \big(x^i\cdot x^j \cdot (x^{\frac{m}{2}+1}+ax) \big)$ (that is a Gram matrix with respect to the basis $\mathcal{B}_\circ$) represents $\Tr_{L_\circ/F} \gw{x^{\frac{m}{2}+1}+ax}$ and is the matrix with $2\times 2$ blocks
    \[
    \begin{bmatrix}
        0 & 0 & \cdots & B \\
         0 & 0 & \iddots & 0 \\
        \vdots & \iddots & \iddots & \vdots \\
         B & 0 & \cdots & 0
    \end{bmatrix}.
    \]
    The statement follows by considering cases according to the parity of $m/2$.
\end{proof}

\begin{proposition}\label{lem-algebratom}
    Let $F$ be a field such that $2m\in F^\times$. Let $d\in F^\times$ and assume that $a^2-b^2d\neq0$. Let $c=a+b\iota\in F[\iota]/(\iota^2-d)$ and $\bar{c}\coloneqq a-b\iota$. 
    Let $L^\phi$ be the sub-algebra of the algebra 
    \[L\coloneqq\frac{F[\iota][u,u']}{(\iota^2-d,u^m-c, u'^m-\overline{c})}\]
    fixed by the involution of $F$-algebras $\phi\colon L\rightarrow L$ which is determined by $\iota \mapsto -\iota$, $u \mapsto u'$, and $u' \mapsto u$.
    
    Then, if $m$ is odd, we have the following equality in $\GW(F)$:
    \begin{equation}\label{eq: m odd}
        \Tr_{L^\phi/F}\gw{1}=\gw{1}+\frac{m-1}{2} \big(\gw{2}+\gw{-2d} \big)+\frac{m(m-1)}{2}\h.
    \end{equation}
    
    Else, if $m$ is even, we have the following equalities in $\GW(F)$:
    \begin{equation}\label{eq:m even 1}
        \Tr_{L^\phi/F}\gw{1}=\gw{1}+\gw{c\bar{c}}+\gw{2a}+\gw{2adc\bar{c}}+\frac{m-2}{2} \big(\gw{2}
            +\gw{-2d} \big)+\frac{m^2-m-2}{2}\h,
    \end{equation}
    \begin{equation}\label{eq:m even uv}
        \Tr_{L^{\phi}/F}\gw{uu'}=
        \frac{m}{2} \big(\gw{2}+\gw{-2d} \big)+\frac{m(m-1)}{2}\h.
    \end{equation}
\end{proposition}

\begin{proof}
Let $L_\circ$ be the $F$-algebra
\[
L_\circ\coloneqq\frac{F[x]}{(x^m-c\bar{c})}.
\]
We have that $L_\circ\subset L^\phi$ by identifying $x\in L_\circ$ with $uu'$. Moreover, the algebra $L^\phi$ is an $L_\circ$-algebra of dimension $m$. Hence, $\Tr_{L^\phi/F}=\Tr_{L_\circ/F}\circ\Tr_{L^{\phi}/L_\circ}$, which we will use for the computation.

In order to compute $\Tr_{L^\phi/L_\circ}\gw{1}$, let
\[
f_i=u^i+u'^i \in L^{\phi}, \qquad \text{for } i=0,1,\dots m,
\]
and the set $\mathcal{B}\coloneqq\{f_0,f_1,    \dots, f_{m-1}\}$.
The set $\mathcal{B}$ is a basis of the algebra $\L^\phi$ as an $L_\circ$-vector space. We have that in $L^\phi$, the product of elements of the basis $\mathcal{B}$ is given by
\begin{equation}\label{eq:proddoubletr}
    f_i\cdot f_j=
    \begin{cases}
        f_{i+j}+x^{\min\{i,j\}}f_{\vert i-j\vert}   &\text{ if } i+j<m,  \\
        2af_{i+j-m}-x^{i+j-m} f_{2m-i-j}+x^{\min\{i,j\}}f_{\vert i-j\vert}   &\text{ if }  i+j\geq m. \\
    \end{cases}
\end{equation}
Hence, the algebraic trace $\tr_{L^\phi/L_\circ}(f_i)$ vanishes for all $0<i<m$ and $\tr_{L^\phi/L_\circ}(f_0)=2m$. So \cref{eq:proddoubletr} implies that the Gram matrix $(\tr_{L^\phi/L_\circ}(f_i\cdot f_j))_{i,j=0}^{m-1}$ representing $\Tr_{L^\phi/L_\circ}\gw{1}$ equals
\begin{equation}\label{eq:GramMatrixdt}
    2m\cdot
    \begin{bmatrix}
        2 & 0&0 & 0 & \cdots & 0 & 0 \\
        0&x & 0 & 0 & \cdots & 0 & a \\
        0&0 & x^2 & 0 & \cdots & a & 0 \\
        0&0 & 0 & x^3 & \cdots & 0 & 0 \\
        \vdots&\vdots & \vdots & \vdots & \ddots & \vdots & \vdots \\
        0&0 & a & 0 & \cdots & x^{m-2} & 0 \\
        0&a & 0 & 0 & \cdots & 0 & x^{m-1}
    \end{bmatrix},
\end{equation}
with the element $x^{\frac{m}{2}}+a$ in the $(m/2+1,m/2+1)$-entry if $m$ is even.
Let $\alpha_i$ be the rank $2$ quadratic form with Gram matrix
\[
    \begin{pmatrix}
        x^i & a \\ a & x^{m-i}
    \end{pmatrix}.
\]
\cref{eq:gwrealtion3} implies that $\alpha_i$ represents the same class in $\GW(L_\circ)$ as 
\[\gw{x^i} \big(\gw{1}+\gw{x^m-a^2} \big)=\gw{x^i} \big(\gw{1}+\gw{-d} \big)\in\GW(L_\circ).\]
Therefore, \cref{eq:GramMatrixdt} implies that $\Tr_{L^\phi/L_\circ}\gw{1}$ equals
\begin{equation*}
    \Tr_{L^\phi/L_\circ}\gw{1}=
    \begin{cases}\displaystyle
        \gw{m}+\gw{2m} \big(\gw{1}+\gw{-d} \big)\sum_{i=1}^{\frac{m-1}{2}}\gw{x^i}
        & m \text{ odd,}\\[6pt] \displaystyle
        \gw{m}+\gw{2m}\gw{x^\frac{m}{2}+a}+\gw{2m} \big(\gw{1}+\gw{-d} \big)\sum_{i=1}^{\frac{m}{2}-1}\gw{x^i}
        & m \text{ even.}\\
    \end{cases}
\end{equation*}
In order to apply the subsequent trace $\Tr_{L^\phi/F}$ to this, recall from \cref{prop:traces} that:
\[
\Tr_{L_\circ/F}\gw{mx^i}=
\begin{cases}
    \gw{(c\bar{c})^i}+\frac{m-1}{2}h
    &\text{ if } m\equiv1\mod2,\\
    \gw{1}+\gw{(-1)^i(c\bar{c})^{i+1}}+\frac{m-2}{2}h
    &\text{ if } m\equiv0\mod2.\\
\end{cases}
\]

If $m$ is odd, by \cref{prop:traces}, we have that
\begin{align*}
    \Tr_{L^\phi/F}\gw{1}&=
    \displaystyle
        \Tr_{L_\circ/F}\gw{m}+\gw{2m} \big(\gw{1}+\gw{-d} \big)\sum_{i=1}^{\frac{m-1}{2}}\Tr_{L_\circ/F}\gw{x^i}\\ &\mkern-6mu\overset{\ref{prop:traces}}{=}
        \qinv{1}+\gw{2m} \big(\gw{1}+\gw{-d} \big)\sum_{i=1}^{\frac{m-1}{2}}\qinv{m(c\bar{c})^i}+\frac{m(m-1)}{2}h
        \\        &\mkern-4mu\overset{\eqref{eq:gwrelation}}{=}
        \qinv{1}+\gw{2} \big(\gw{1}+\gw{-d} \big)\sum_{i=1}^{\frac{m-1}{2}}\qinv{1}+\frac{m(m-1)}{2}h.
\end{align*}
If $m$ is even, we observe that
\begin{align*}
    \big(\gw{1}+\gw{-d} \big)\Tr_{L_\circ/F}\gw{mx}&\mkern-6mu\overset{\ref{prop:traces}}{=}
    \big(\gw{1}+\gw{-d} \big)\frac{m}{2}h=mh \\
    \big(\gw{1}+\gw{-d} \big)\Tr_{L_\circ/F}\gw{m}&\mkern-6mu\overset{\ref{prop:traces}}{=}
    \big(\gw{1}+\gw{-d} \big)\left(\gw{1}+\gw{c\bar{c}} +\frac{m-2}{2}h\right)\\ &\mkern-4mu\overset{\eqref{eq:gwrelation}}{=}
    2\big(\gw{1}+\gw{-d}\big)+(m-2)h. 
\end{align*}
Therefore, by counting the number of even and odd indices $i$, we get that
\begin{equation}\label{eq:sumtrs}
    \big(\gw{1}+\gw{-d} \big)\sum_{i=1}^{\frac{m}{2}-1}\Tr_{L_\circ/F}\gw{mx^i}{}
    =
    \begin{cases}
        \frac{m-4}{2}\big(\gw{1}+\gw{-d} \big)+\frac{m^2-3m+4}{2}h
        & \text{ if } m\equiv 0 \mod 4,\\[4pt]
        \frac{m-2}{2}\big(\gw{1}+\gw{-d} \big)+\frac{(m-1)(m-2)}{2}h
        & \text{ if } m\equiv 2 \mod 4.\\
    \end{cases}
\end{equation}

Therefore, this implies that if $m\equiv0 \mod 4$, then we have that
\begin{align*}
    \Tr_{L^\phi/F}\gw{1}&\mkern-6mu\overset{\ref{prop:traces}}{=}
    \displaystyle
        \gw{1}+\gw{c\bar{c}}+\gw{2m}\Tr_{L_\circ/F}\gw{x^\frac{m}{2}+a}+\gw{2} \big(\gw{1}+\gw{-d} \big)\sum_{i=1}^{\frac{m}{2}-1}\Tr_{L_\circ/F}\gw{mx^i}+\frac{m-2}{2}h\\ &\mkern-4mu\overset{\eqref{eq:sumtrs}}{=}
        \gw{1}+\gw{c\bar{c}}+\gw{2m}\Tr_{L_\circ/F}\gw{x^\frac{m}{2}+a}+
        \frac{m-4}{2}\big(\gw{2}+\gw{-2d} \big) +\frac{m^2-2m+2}{2}h
        \\        &\mkern-6mu\overset{\ref{lm:intertrace}}{=}
        \gw{1}+\gw{c\bar{c}}+\gw{2a}+\gw{2adc\bar{c}}+\frac{m-2}{2} \big(\gw{2}
            +\gw{-2d} \big)+\frac{m^2-m-2}{2}\h.
\end{align*}
Else, if $m\equiv2\mod 4$, then
\begin{align*}
    \Tr_{L^\phi/F}\gw{1}&\mkern-6mu\overset{\ref{prop:traces}}{=}
    \displaystyle
        \gw{1}+\gw{c\bar{c}}+\gw{2m}\Tr_{L_\circ/F}\gw{x^\frac{m}{2}+a}+\gw{2} \big(\gw{1}+\gw{-d} \big)\sum_{i=1}^{\frac{m}{2}-1}\Tr_{L_\circ/F}\gw{mx^i}+\frac{m-2}{2}h\\ &\mkern-4mu\overset{\eqref{eq:sumtrs}}{=}
        \gw{1}+\gw{c\bar{c}}+\gw{2m}\Tr_{L_\circ/F}\gw{x^\frac{m}{2}+a}+
        \frac{m-2}{2}\big(\gw{2}+\gw{-2d} \big) +\frac{m^2-2m}{2}h
        \\        &\mkern-6mu\overset{\ref{lm:intertrace}}{=}
        \gw{1}+\gw{c\bar{c}}+\gw{2a}+\gw{2adc\bar{c}}+\frac{m-2}{2} \big(\gw{2}
            +\gw{-2d} \big)+\frac{m^2-m-2}{2}\h.
\end{align*}

Lastly, in order to compute $\Tr_{L^\phi/F}\gw{uu'}$ for $m$ even, recall that $uu'=x\in L_\circ$. Hence, $\Tr_{L^\phi/L_\circ}\gw{x}$ equals $\gw{x}\cdot\Tr_{L^\phi/L_\circ}\gw{1}$, and so
\begin{align*}
    \Tr_{L^\phi/F}\gw{uu'}&=
    \displaystyle
        \Tr_{L_\circ/F} \gw{x}\cdot
        \left(\displaystyle
        \gw{m}+\gw{2m}\gw{x^\frac{m}{2}+a}+\gw{2m} \big(\gw{1}+\gw{-d} \big)\sum_{i=1}^{\frac{m}{2}-1}\gw{x^i}\right)\\ &\mkern-6mu\overset{\ref{prop:traces}}{=}
        \frac{m}{2}h+\gw{2m}\Tr_{L_\circ/F} \gw{x^{\frac{m}{2}+1}+ax}+\gw{2m} \big(\gw{1}+\gw{-d} \big)\sum_{i=1}^{\frac{m}{2}-1}\Tr_{L_\circ/F}\gw{x^{i+1}}.\\
\end{align*}
Proceeding analogously, we have that
\[
    \big(\gw{1}+\gw{-d} \big)\sum_{i=1}^{\frac{m}{2}-1}\Tr_{L_\circ/F}\gw{mx^{i+1}}=
    \begin{cases}
        \frac{m}{2}\big(\gw{1}+\gw{-d} \big)+\frac{m^2-3m}{2}h
        & \text{ if } m\equiv 0 \mod 4,\\[4pt]
        \frac{m-2}{2}\big(\gw{1}+\gw{-d} \big)+\frac{(m-1)(m-2)}{2}h
        & \text{ if } m\equiv 2 \mod 4.\\
    \end{cases}
\]
Therefore, if $m\equiv0\mod 4$, then
\begin{align*}
    \Tr_{L^\phi/F}\gw{uu'}
    \displaystyle
        &=\gw{2m}\Tr_{L_\circ/F} \gw{x^{\frac{m}{2}+1}+ax}+
        \frac{m}{2}\big(\gw{2}+\gw{-2d} \big)+\frac{m^2-2m}{2}h\\        &\mkern-6mu\overset{\ref{lm:intertrace}}{=}
        \frac{m}{2}\big(\gw{2}+\gw{-2d} \big)+\frac{m^2-m}{2}h.
\end{align*}
Else, if $m\equiv2\mod4$, then
\begin{align*}
    \Tr_{L^\phi/F}\gw{uu'}&=
    \displaystyle
        \gw{2m}\Tr_{L_\circ/F} \gw{x^{\frac{m}{2}+1}+ax}+
        \frac{m-2}{2}\big(\gw{2}+\gw{-2d} \big)+\frac{m^2-2m+2}{2}h\\        &\mkern-6mu\overset{\ref{lm:intertrace}}{=}
        \frac{m}{2}\big(\gw{2}+\gw{-2d} \big)+\frac{m^2-m}{2}h.
\end{align*}
\end{proof}

\subsection{Deformation patterns}
Edges of the underlying graph $\Gamma$ contribute to the enriched count because of the following phenomenon: given a bounded edge of weight $\geq 2$, any curves $C_v$ lifting the vertices $v$ at either end meet their joint toric boundary in a single point with higher intersection multiplicity. Perturbing this intersection point reveals more nodes. 
The contributions coming from non-twin edges or refined points conditions given by $k$-rational points have been computed in \cite{JPP23}. We extend these results to the general case in the following lemmas.

\begin{remark}
    In \cite{JPP23} one considers plane tropical curves and not parametrized ones. In the setting of plane tropical curves one needs to consider equivalence classes of \emph{extended edges}. However, this agrees with just considering edges of the underlying graph $\Gamma$. 
    To be more precise, we forget the marked points momentarily and consider \emph{edges} of the underlying graph after forgetting. We speak of \emph{marked edges} when we consider edges that previously contained a marked point (and had been subdivided into two edges by the marked point accordingly).
\end{remark}

Deformation patterns give a contribution for each edge of the dual subdivision. More precisely, if $F$ is the algebra defined by the possible algebraic curves for the polygons $\Delta_v$ in the dual subdivision for $v$ ranging over all vertices, then every edge $e$ of $\Gamma$ defines an $F$-algebra $L_{e}$. Refined point conditions define similar algebras $M_e$ for each marked edge $e$ of the underlying graph $\Gamma$.

\begin{proposition}\label{prop-deformpatternandrefinedpoint}
\begin{enumerate}
    \item Let $e$ be an edge of weight $m$ which does not belong to a twin tree. Then the algebra defining all deformation patterns for $e$ is given by
	\begin{equation*}
		L_{e} = \begin{cases}
            \begin{aligned}
                \frac{F[u]}{(u^m - 1)}  &&& \text{if $m$ is odd}, \\[1ex]
			\frac{F[u]}{(u^m - c)}  &&& \text{if $m$ is even}.
            \end{aligned}			
		\end{cases}
	\end{equation*}
     for some $c\in F^\times$.
     This defines a curve $C_e$ over $L_e$ consisting of all possible deformation patterns.
     The quadratic weight $\Wel(C_e) \in \GW(L_e)$ of the curve $C_e$ (defined over $L_e$) is 
     \[\Wel(C_e) = \begin{cases}
         \gw{1} & \text{if $m$ is odd,} \\
         \gw{-u} & \text{if $m$ is even.}
     \end{cases} \]

     \item  Let $e$ be an edge of weight $m>1$ which is part of a twin tree, so $e$ is necessarily marked by a point in $k(\sqrt{d})$. Let $e'$ be its twin edge (necessarily also of weight $m$). Denote by $L_{e,e'}$ the $F$-algebra defined by all deformation patterns.
     Then $L_{e,e'}$ is the subalgebra $L^\phi$  of
    	\begin{equation*}
    		L = \begin{cases}
    			\frac{F[\iota, u,u']}{\big(\iota^2 - d, u^m - 1,u'^m-1 \big)}  & \text{if $m$ is odd}, \\
    			\frac{F[\iota, u,u']}{\big(\iota^2 - d, u^m - c, u'^m-\overline{c} \big)} & \text{if $m$ is even}
    		\end{cases}
    	\end{equation*}
     fixed under the involution 
     \[
         \phi\colon L \longrightarrow L \colon \iota \longmapsto -\iota \colon u\longmapsto u' \colon  u'\longmapsto u.
     \]
     Here, $c=a+b\iota \in F[\iota]/(\iota^2 - d) \cong F(\sqrt{d})$ and $\overline{c}=a-b\iota$ is its Galois conjugate. All deformation patterns define $C_{e,e'}$ which is a nodal curve defined over $L_{e,e'}$.
    The quadratic weight $\Wel(C_{e,e'}) \in \GW(L_{e,e'})$ of $C_{e,e'}$ is given by 
    \[ \Wel(C_{e,e'})= \begin{cases}
        \gw{1} & \text{if $m$ is odd,} \\
        \gw{uu'} & \text{if $m$ is even.}
    \end{cases} \]
    
    \item Let $e$ be a (not necessarily marked) edge in a twin tree of weight $m=1$ and let $e'$ be its twin edge. Then $L_{e,e'}=k$ and $\Wel(C_{e,e'})=\gw{1}\in \GW(L_{e,e'})$.
\end{enumerate}
\end{proposition}

\begin{proof}
The first part follows from \cite[Proposition~5.10]{JPP23}. The only case to pay attention to is that here, the edge can be merged with another edge and end in a vertex of type \fourValentVertex{}, while in \cite[Proposition~5.10]{JPP23}, all end vertices were $3$-valent. It turns out however that this does not change the main arguments, as we point out in the following.

Assume an edge of weight $m$ (which is not part of a twin tree) ends at a vertex of type \fourValentVertex{}. 
To reveal the hidden nodes, we pass to a tropical refinement as in Lemma 3.9 \cite{Shu04}. 
A vertex of type \fourValentVertex{} in a tropical stable map passing through vertically stretched point conditions corresponds to a curve given by an equation of the form
$$f(x,y)=y+a\cdot (x-\varepsilon)^p\cdot (x-\overline{\varepsilon})^q,$$
where $p$ and $q$ are the weights of the two parallel edges adjacent to the vertex of type \fourValentVertex{}. 
We can assume that the valuation of $\varepsilon, \overline{\varepsilon},a$ are all $0$, while the valuations of the other monomials (including the ones corresponding to the other end vertices of our edge and its parallel edge) are bigger than~$0$. As in  \cite[Lemma 3.9]{Shu04}, we change coordinates by substituting $x$ by $x+\varepsilon$. The valuation of the coefficient of the factor $(x+\varepsilon-\overline{\varepsilon})$ is still $0$, and so we obtain a dual subdivision which locally looks as depicted in Figure \ref{fig-refinement}.
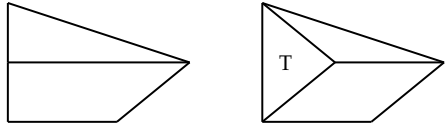
\begin{figure}[t]
    \centering

\tikzset{every picture/.style={line width=0.75pt}} 

\begin{tikzpicture}[x=0.75pt,y=0.75pt,yscale=-1,xscale=1]

\draw    (293.33,420) -- (330,390) ;
\draw    (238.33,360) -- (330,390) ;
\draw    (238.33,420) -- (275,390) ;
\draw    (165,420) -- (201.67,390) ;
\draw    (110,390) -- (110,420) ;
\draw    (275,390) -- (330,390) ;
\draw    (238.33,360) -- (275,390) ;
\draw    (110,420) -- (165,420) ;
\draw    (110,360) -- (201.67,390) ;
\draw    (110,390) -- (110,360) ;
\draw    (110,390) -- (201.67,390) ;
\draw    (238.33,420) -- (293.33,420) ;
\draw    (238.33,420) -- (238.33,360) ;

\draw (245.33,384.67) node [anchor=north west][inner sep=0.75pt]  [font=\scriptsize] [align=left] {T};
\end{tikzpicture}

    \caption{Refinement of an edge adjacent to a vertex of type \fourValentVertex{}.}
    \label{fig-refinement}
\end{figure}
The triangle $T$ we observe is now as in \cite[Lemma 3.9]{Shu04} and we can continue with it as there, yielding the same result for $L_{e}$ and $\Wel(C_e)$ as in \cite[Proposition~5.10]{JPP23}.

For the second part, we basically do the same twice and everything has to be Galois conjugate, that is invariant under $\phi$. 

To compute $\Wel(C_{e,e'})$, we take the product of the quadratic weights associated with the two conjugate deformation patterns obtained in this manner, that is we get $\gw{u}\cdot \gw{u'}$ when $m$ is even and $\gw{1}\cdot \gw{1}$ when $m$ is odd.

The third part follows immediately (2) with $m=1$.
\end{proof}

\subsection{Refined point conditions}

The next propositions will define the algebras of the refined point conditions.
\begin{proposition}\label{prop:algebrarefinedpointconditions}
    Let $e$ be a marked edge of weight $m$. 
    \begin{enumerate}
        \item If $e$  does not belong to a twin tree, then 
        	\begin{equation*}
        		M_{e} = 
        			\frac{L[w]}{(w^m - \alpha)},
        	\end{equation*}
         for some $\alpha\neq 0$ determined by the point conditions. 

         \item If $e$ is part of a twin tree, let $e'$ be its twin edge (necessarily also of weight $m$), and then $M_{e,e'}$ is the sub-algebra of 
        	\begin{equation*}
        		M = \frac{L[\iota, w,w']}{\big( \iota^2 - d, w^m-\alpha,w'^m-\overline{\alpha} \big)}
        	\end{equation*}
         fixed by the involution of $L$-algebras $\phi\colon M\rightarrow M$ determined by $\iota \mapsto -\iota$, $w \mapsto w'$, and $w' \mapsto w$.
    \end{enumerate}
\end{proposition}

\begin{proof}
    This follows from \cite[(3.8)]{Shu06b}; see also \cite[Section 3.4]{JPP23}. Specifically, one must adjoin an $m$-th root of an element determined by the refined point condition. For twin edges, this process must be performed twice, ensuring that all constructions remain invariant under the Galois conjugation induced by the point condition defined over a quadratic field extension. Consequently, this induces an involution on the algebra of solutions.
\end{proof}

\subsection{Contribution from edges}

For this subsection let $F$ again be a finite étale $k$-algebra. Further let $L_e$, $L_{e,e'}$, $M_e$ and $M_{e,e'}$ be as in Proposition \ref{prop-deformpatternandrefinedpoint} and Proposition \ref{prop:algebrarefinedpointconditions}. The following proposition explains the definition of $\mult^{\AA^1}(e)$.

\begin{proposition}[Non-twin edge multiplicity]
\label{prop:nontwinedgemult}
    For any marked edge $e$ which is not part of a twin tree we get $F\subset L_e\subset M_e$ and $\Tr_{M_e/F}\big(\Wel(C_e)\big)=\mult^{\AA^1}(e)$ (see Definition \ref{def:twinedgemult}). More precisely:
    \begin{enumerate}
        \item When the weight $m$ of $e$ is odd, then $\Wel(C_e)=\gw{1}$ and
        $\Tr_{M_e/F}(\gw{1})=\frac{m^2-1}{2}\h+\gw{1}$.
        \item When the weight $m$ of $e$ is even, then $\Wel(C_e)=\gw{u}$ (with $u\in L_e$ as in Proposition \ref{prop-deformpatternandrefinedpoint}) and $\Tr_{M_e/F}(\gw{u})=\frac{m^2}{2}\h$.
    \end{enumerate}
\end{proposition}
\begin{proof}
    \begin{enumerate}
        \item Assume $m$ is odd and recall from \cref{prop-deformpatternandrefinedpoint} that $L_e=\quotient{F[u]}{(u^m-1)}$ and from and \cref{prop:algebrarefinedpointconditions} that $M_e=\quotient{L_e[w]}{(w^m-\alpha)}$. So
        \begin{align*}
            \Tr_{M_e/F}(\gw{1})&=\Tr_{L_e/F}\circ\Tr_{M_e/L_e}(\gw{1}) \\
            &\overset{\ref{prop:traces}}{=} \Tr_{L_e/F}\left(\qinv{m}+\frac{m-1}{2}\h\right)\\
           &\overset{\ref{prop:traces}}{=} \gw{m^2}+m\cdot \frac{m-1}{2}\h+ \frac{m-1}{2}\h \\
           &=\gw{1}+\frac{m^2-1}{2}\h .
        \end{align*}
        
        \item Assume $m$ is even and recall from \cref{prop-deformpatternandrefinedpoint} that $L_e=\quotient{F[u]}{(u^m-ac)}$ and from \cref{prop:algebrarefinedpointconditions} that $M_e=\quotient{L_e[w]}{(w^m-\alpha)}$. So
         \begin{align*}
            \Tr_{M_e/F}(\gw{u})&=\Tr_{L_e/F}\circ\Tr_{M_e/L_e}(\gw{u}) \\
            &\overset{\ref{prop:traces}}{=} \Tr_{L_e/F}\left(\qinv{m u}+\qinv{m u \alpha}+\frac{m-2}{2}\h\right)\\
           &\overset{\ref{prop:traces}}{=} \frac{m}{2}\h+\frac{m}{2}\h + m\cdot \frac{m-2}{2}\h \\
           &=\frac{m^2}{2}\h .
        \end{align*}
    \end{enumerate}
\end{proof}

The following proposition explains the definition of the twin edge multiplicity $\mult^{\A^1}(e, e')$ from \cref{def:twinedgemult}.  

\begin{proposition}[Twin edge multiplicities]
    \label{prop:twinedgemult}
    Let $e$ and $e'$ be marked twin edges, with a point in $k(\sqrt{d})$. Then we get $F\subset L_{e,e'}\subset M_{e,e'}$ and $\Tr_{M_{e, e'} / F} \big(\Wel(C_{e,e'}) \big) = \mult^{\AA^1}(e,e')$ (see Definition \ref{def:twinedgemult}). More precisely:
    \begin{enumerate}
        \item When the weight $m$ of $e$ and $e'$ is odd then
        $$\Tr_{M_{e,e'}/F} \big(\Wel(C_{e,e'}) \big)=\Tr_{M_{e,e'}/F}(\gw{1})=\gw{1}+\frac{m^2-1}{2} \big(\gw{1}+\gw{-d} \big)+\frac{m^4-m^2}{2}\h\in \GW(F).$$
        
        \item When the weight $m$ of $e$ and $e'$ is even then
            $$\Tr_{M_{e,e'}/F} \big(\Wel(C_{e,e'}) \big)=\Tr_{M_{e,e'}/F} \big(\gw{uu'} \big)=\frac{m^2}{2} \big(\gw{1}+\gw{-d} \big)+ \frac{m^4-m^2}{2}\h\in \GW(F).$$
    \end{enumerate}
\end{proposition}
\begin{proof}
    \begin{enumerate}
        \item If $m$ is odd then $\Wel(C_{e, e'}) = \gw{1}$ and
        \begin{align*}
            \Tr_{M_{e,e'}/F}(\gw{1})&=\Tr_{L_{e,e'}/F}\circ\Tr_{M_{e,e'}/L_{e,e'}}(\gw{1})\\
            &\overset{\eqref{eq: m odd}}{=}\Tr_{L_{e,e'}/F}\left(\gw{1}+\frac{m-1}{2} \big(\gw{2}+\gw{-2d} \big)+\frac{m(m-1)}{2}\h\right)\\
            &\overset{\eqref{eq: m odd}}{=} \gw{1}+2\cdot \frac{m-1}{2} \big(\gw{2}+\gw{-2d} \big)+\frac{(m-1)^2}{2} \big(\gw{1}+\gw{-d} \big)+\frac{m^4-m^2}{2}\h\\
            &\overset{\eqref{eq:gwrealtion2}}{=}\gw{1}+\left(\frac{(m-1)^2}{2}+m-1\right) \big(\gw{1}+\gw{-d} \big)+\frac{m^4-m^2}{2}\h\\
            &=\gw{1}+\frac{m^2-1}{2} \big(\gw{1}+\gw{-d} \big)+\frac{m^4-m^2}{2}\h.
        \end{align*}
        \item If $m$ is even, then recall that $\Wel(C_{e, e'}) = \gw{uu'}$ with $uu'\in L_{e,e'}$.
        \begin{align*}
            \Tr_{M_{e,e'}/F} \big(\gw{uu'} \big) &= \Tr_{L_{e,e'}/F}\circ \Tr_{M_{e,e'}/L_{e,e'}}\big(\gw{uu'}\big)\\
            &\overset{\eqref{eq:m even 1}}{=} \Tr_{L_{e,e'}/F} \Big(\gw{uu'}+\gw{(a^2-b^2d)uu'}+\gw{2auu'}+\gw{2ad(a^2-b^2d)uu'}\\
            &\qquad+\frac{m-2}{2} \big(\gw{2uu'}+\gw{-2duu'} \big) +\frac{m^2-m-2}{2}\h\Big)\\
            &\overset{\eqref{eq:m even uv}}{=} \frac{m}{2} \big(\gw{2}+\gw{-2d} \big)+\frac{m}{2}\Big(\gw{2(a^2-b^2d)}+\gw{-2d(a^2-b^2d)} \Big)\\
            &\qquad+\frac{m}{2} \big(\gw{a}+\gw{-ad} \big)+\frac{m}{2} \Big(\gw{ad(a^2-b^2d)}+\gw{-a(a^2-b^2d)} \Big)\\
            &\qquad+\frac{m(m-2)}{4} \big(\gw{1}+\gw{-d}+\gw{-d}+\gw{1} \big)+\frac{m^4-m^2-2m}{2}\h\\
            &\overset{\eqref{eq:gwrelation}}{=} \frac{m}{2} \big(\gw{2}+\gw{-2d} \big)+\frac{m}{2}\big(\gw{2}+\gw{-2d} \big) +\frac{m}{2} \big(\gw{a}+\gw{-ad} \big) \\
            &\qquad+\frac{m}{2} \big(\gw{ad}+\gw{-a} \big)+\frac{m(m-2)}{2} \big(\gw{1}+\gw{-d} \big)+\frac{m^4-m^2-2m}{2}\h
            \\
            &\overset{\eqref{eq:gwrealtion2}}{=} \frac{m^2}{2} \big(\gw{1}+\gw{-d} \big)+\frac{m^4-m^2}{2}\h.
         \end{align*}
    \end{enumerate}
\end{proof}

\section{Twin trees}\label{sec-twintree}

Let $\cT$ be a twin tree. We use the same notation as in \cref{def-multtwintree}, i.e. we call $q_1, \ldots, q_t$ the double points on $\cT$ corresponding to field extensions $k \subset k(\sqrt{d_i})$.
In this section we prove that $\mult^{\AA^1} (\cT)$ as defined in \cref{def-multtwintree} is the contribution of the algebraic lifts of $\cT$ to $\mult^{\AA^1}(\Gamma,f)$.
We have already seen in \cref{lem-structuredoublepart} that $\cT$ has a unique root vertex of type \fourValentVertex{}, so in particular $\mult^{\AA^1}(\cT)$ is well defined.

\begin{lemma} \label{lem-thirdend}
    Consider a vertex $v$ of type \allDoubleVertex{} as in \cref{lem:fatvertex}: 
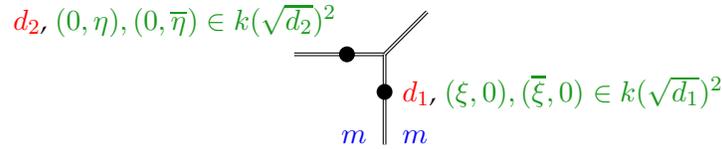
    \begin{center}
        \begin{tikzpicture}[baseline={([yshift=-.5em]current bounding box.center)}]
			\path[draw, double] (0,0) --+ (180:1.2);
			\path[draw, double] (0,0) --+ (270:1.2);
			\path[draw, double] (0,0) --+ (45:0.8);
            \draw (-0.1, -1.1) node[anchor = east] {\edgeweight{$m$}};
            \draw (0.1, -1.1) node[anchor = west] {\edgeweight{$m$}};
   
			\fatpoint{-0.5, 0}
            \draw (-0.5, 0.1) node[anchor = south east] {\quadextension{$d_2$}, \coordinates{$(0, \eta), (0, \conj{\eta}) \in k(\sqrt{d_2})^2$}};
			\fatpoint{0, -0.5}
            \draw (0.1, -0.5) node[anchor = west] {\quadextension{$d_1$}, \coordinates{$(\xi, 0), (\conj{\xi}, 0) \in k(\sqrt{d_1})^2$}};
		\end{tikzpicture}
    \end{center}
    Let $C_v$ be an algebraic lift of this situation. 
    Then the points of $C_v$ corresponding to the third end of $v$ are defined over $k[\iota_1, \iota_2]/(\iota_1^2 - d_1, \iota_2^2 - d_2)$ and are conjugate under the involution $\conj{\cdot} : \iota_i \mapsto -\iota_i$.
\end{lemma}

\begin{proof}
    From the proof of \cref{lem:fatvertex} we know that the algebraic lift $C_v$ of the local situation at $v$ is given by $C_v = V(f_1 \conj{f_1})$ with $f_1(x, y) = y - a(x - \xi)^m$ and $a = \eta / (-\xi)^m$. Therefore, the two intersection points of $C_v$ with the third toric boundary divisor are conjugate as claimed and we only need to show that they are defined over $k[\iota_1, \iota_2]/(\iota_1^2 - d_1, \iota_2^2 - d_2)$.
    
    To this end note that $f_1$ defines a curve in the complete toric surface given by the fan
    
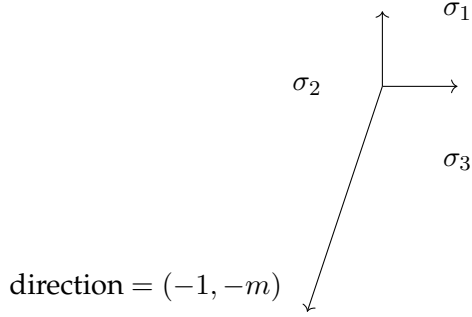
\begin{center}
        \begin{tikzpicture}
            \path[draw, ->] (0,0) -- (0, 1);
            \path[draw, ->] (0,0) -- (1, 0);
            \path[draw, ->] (0,0) -- (-1, -3);  
            \draw (-1.2, -3) node[anchor = south east] {direction $= (-1, -m)$};
    
            \draw (1, 1) node {$\sigma_1$};
            \draw (-1, 0) node {$\sigma_2$};
            \draw (1, -1) node {$\sigma_3$};        
        \end{tikzpicture}
    \end{center}
    and $V(f_1)$ is the affine part in $U_1 = \Spec k[\sigma_1^\vee \cap \ZZ^2] = \Spec k[x, y]$. We compute the change of charts to pass to 
    \[U_2 = \Spec k[\sigma_2^\vee \cap \ZZ^2] = \Spec k[x^{-1}, x^{-m}y] \cong \Spec k[u, v].\]    
    The intersection $U_1 \cap U_2$ is precisely $\Spec k[x, x^{-1}, y]$, i.e. on $U_1 \cap U_2$ we can have $x \neq 0$. Therefore,
    \[ y - a(x - \xi)^m = x^m \big( x^{-m}y - a (1 - \xi x^{-1})^m \big) \]
    and dropping the leading $x^m$ gives the presentation $g(u, v) = v - a(1 - \xi u)^m$ on $U_2$. 

    In the chart $U_2$ the toric boundary at $v = 0$ corresponds to $y = 0$ in $U_1$. So the \enquote{third} toric boundary divisor, which we are after, is the one given by $u = 0$. The point of intersection is $g(0, v) = 0$, i.e. $v = a$. Since $a \in k[\iota_1, \iota_2]/(\iota_1^2 - d_1, \iota_2^2 - d_2)$ the proof is complete.
\end{proof}

We now use \cref{lem-thirdend} in order to inductively determine the $k$-algebra $F_\cT$ of all algebraic lifts of the vertices in $\cT$.

\begin{lemma} \label{lem-algebratwintree}
    The $k$-algebra of algebraic lifts of the vertices in $\cT$ is given by
    \[ F_{\cT}  
    = \frac{k[\iota_{1,2}, \iota_{2,3}, \ldots \iota_{t-1,t}]}{\big( \iota_{i, i+1}^2 - d_id_{i+1} \bigmid i = 1, \ldots, t-1 \big)}. \]
\end{lemma}

\begin{proof}
    To prove the claim on $F_\cT$ we proceed by induction on $t$.
    Let $t = 1$, i.e. $\cT$ consists of a single unbounded twin edge, marked with the double point $q_1$ and $\cT$ has no vertex of type \allDoubleVertex{}. Then $F_\cT = k$ by \cref{lem:type_C}. This coincides with the claim since the list of new variables $\iota_{i, i+1}$ is empty in this case. 
    
    Now let $t \geq 2$. Then there exists at least one vertex $v_0 \in \cT$ with the following properties (compare \cref{fig-twintree}):
    \begin{itemize}
        \item only one of the (twin) edges incident to $v_0$ is bounded, say $(e_0, e_0')$, and 
        \item at $v_0$ two of the three directions are fixed by double points. Without loss of generality we assume these are $q_1$ and $q_2$.
    \end{itemize} 
    We consider the \enquote{truncated} twin tree $\cT'$ obtained from $\cT$ by cutting the bounded twin edge $(e_0, e_0')$. 
    This new tree is completely fixed by the double points $q_3, \ldots, q_t$ plus one more double point. 
    
    \textbf{Case 1:} $(e_0, e_0')$ was marked. Without loss of generality we assume that $q_2$ lies on $(e_0, e_0')$. Then $\cT'$ is fixed by $q_2, \ldots, q_t$ and by induction we get
    \begin{align*}
        F_\cT = F_{v_0} \otimes_k F_{\cT'} &= \frac{k[\iota_{1,2}]}{(\iota_{1,2}^2 - d_1d_2)} \otimes_k \frac{k[ \iota_{2,3}, \ldots, \iota_{t-1,t}]}{\big( \iota_{i, i+1}^2 - d_id_{i+1} \bigmid i = 2, \ldots, t-1 \big)} \\
        &= \frac{k[\iota_{1,2}, \iota_{2,3}, \ldots \iota_{t-1,t}]}{\big( \iota_{i, i+1}^2 - d_id_{i+1} \bigmid i = 1, \ldots, t-1 \big)}.
    \end{align*}

    \textbf{Case 2:} $(e_0, e_0')$ was not marked. In this case we use \cref{lem-thirdend} on $v_0$ to obtain a new double point $q'$ defined over the field extension 
    \[k(\sqrt{d_1d_2}) \subset k(\sqrt{d_1d_2})(\sqrt{d_2}) = k(\sqrt{d_1}, \sqrt{d_2}). \]
    We may now consider $\cT'$ as being fixed by the points $q', q_3, \ldots, q_t$ and the $k(\sqrt{d_1d_2})$-algebra defined by $\cT'$ is by induction
    \[ \frac{\Big(\frac{k[\iota_{1,2}]}{(\iota_{1,2}^2 - d_1d_2)}\Big) [\iota_{2, 3}, \ldots, \iota_{t-1, t}]}{\big( \iota_{i, i+1}^2 - d_id_{i+1} \bigmid i = 2, \ldots, t-1 \big)} = \frac{k[\iota_{1,2}, \iota_{2,3}, \ldots \iota_{t-1,t}]}{\big( \iota_{i, i+1}^2 - d_id_{i+1} \bigmid i = 1, \ldots, t-1 \big)}. \]
    But this means that the $k$-algebra $F_\cT$ is obtained from $k$ by first adjoining $\iota_{1,2}$ in order to pass to $F_{v_0}$ and then adjoining the remaining $\iota_{i, i+1}$ and hence the claim follows in this case as well.
\end{proof}

The quadratic weight of an algebraic lift of $\cT$ is given by
\[ \Wel(C_\cT) = \underbrace{\prod_{v \text{ vertex}} \Wel(C_v)}_{\in \GW(F_\cT)} \cdot \big(\text{contribution from deformation patterns} \big).  \]
For now, we focus on the contribution from vertices which we will temporarily denote by 
\begin{equation}
    \label{eq:WelCV}
    \Wel(C_V) = \prod_{v \text{ vertex}} \Wel(C_v).
\end{equation}
Recall from \cref{def-multtwintree} that $m_\circ$ is the weight of the root of $\cT$ plus the number of unbounded elevators in $\cT$.

\begin{lemma} \label{lem-WelCT}
    Let $d_\circ$ be any of the $d_i$. Then
    $\Wel (C_V) = \gw{d_\circ^{m_\circ}} \in \GW(F_\cT)$. 
\end{lemma}

\begin{proof}
    Since we are working with vertically stretched point configurations we know that every (twin) elevator edge $(e, e')$ in $\cT$ is marked and hence has an associated $d(e, e') \in \{d_1, \ldots, d_t\}$. 
    If $t = 1$ then by \cref{lem:type_C} we have $\Wel(C_V) = \gw{1}$. 
    If $t \geq 2$, then there are vertices of type \allDoubleVertex{}.
    By \cref{lem:fatvertex}, the quadratic weight $\Wel(C_v)$ at a vertex $v$ of type \allDoubleVertex{} is $\gw{d(e, e')^{m(e,e')}}$ for the unique elevator $(e, e')$ incident to $v$. We rewrite $\Wel(C_V)$ as a product over elevators $(e, e')$ in $\cT$. For this note that each $(e, e')$ is precisely of one of the following types:
    \begin{itemize}        
        \item Bounded elevator between two vertices of type \allDoubleVertex{}. The associated $\gw{d(e, e')^{m(e,e')}}$ will appear twice as a factor in $\Wel(C_V)$ and hence cancel. 
        \item Unbounded elevator incident to one vertex of type \allDoubleVertex{}. Then $m(e, e') = 1$ and we get a factor of $\gw{d(e, e')}$.
        \item The elevator connecting $\cT$ to its root vertex of type \fourValentVertex{}. Denote the weight of this elevator by $m_{\mathrm{root}}$ and its corresponding $d$ by $d_{\mathrm{root}}$.
    \end{itemize}
    In summary we see that 
    \begin{equation} \label{eq-WelCV_almost}
        \Wel(C_V) = \gw{(d_{\mathrm{root}})^{m_{\mathrm{root}}}} \cdot \prod_{(e,e') \text{ unbounded elevator}} \gw{d(e,e')}.
    \end{equation}

    We note that in $\GW(F_\cT)$ it holds $\gw{d_i} = \gw{d_j}$ for all $i$ and $j$: indeed, $d_id_{i+1}$ is a square in $F_\cT$, i.e. $\gw{d_i} = \gw{d_i\cdot d_id_{i+1}} = \gw{d_{i+1}}$, and iterating this shows $\gw{d_i} = \gw{d_j}$. In particular, in the right hand side of \cref{eq-WelCV_almost} we can pick any of the $d_i$ as our $d_\circ$ and obtain the expression from the claim.
\end{proof}

\begin{lemma} \label{lem-tracetwintree}
The trace of the quantity in Equation (\ref{eq:WelCV}) equals:
\[\Tr_{F_\cT / k} \big(\Wel(C_V) \big) = \gw{2^{t-1}} \cdot \sum_{\substack{I\subset \{1,\dots,t\}\\ \vert I\vert\equiv m_{\circ}\mod 2}} \qinv{\prod_{i\in I} d_i}.\]
\end{lemma}

\begin{proof}
    By \cref{lem-WelCT} we know $\Wel(C_V) = \gw{d_\circ^{m_\circ}}$. By iterated application of \cref{prop:traces} we compute
    \begin{align*}
        &\phantom{=} \Tr_{F_\cT / k} \gw{d_\circ^{m_\circ}} \\
        &= \gw{d_\circ^{m_\circ}} \cdot \Tr_{F_\cT / k} \gw{1} \\
        &= \gw{d_\circ^{m_\circ}} \cdot \Tr_{k(\sqrt{d_1d_2}) / k} \cdots \Tr_{k(\sqrt{d_1d_2}, \ldots, \sqrt{d_{t-1}d_t}) / k(\sqrt{d_1d_2}, \ldots, \sqrt{d_{t-2}d_{t-1}})} \gw{1}  \\
        &= \gw{d_\circ^{m_\circ}} \cdot \big( \gw{2} + \gw{2d_{t-1}d_t} \big) \cdot
        \Tr_{k(\sqrt{d_1d_2}) / k} \cdots \Tr_{k(\sqrt{d_1d_2}, \ldots, \sqrt{d_{t-2}d_{t-1}}) / k(\sqrt{d_1d_2}, \ldots, \sqrt{d_{t-3}d_{t-2}})} \gw{1} \\
        &= \cdots \\
        &= \gw{d_\circ^{m_\circ}} \cdot \prod_{i = 1}^{t-1} \big( \gw{2} + \gw{2d_{i}d_{i+1}} \big) \\
        &= \gw{d_\circ^{m_\circ}} \cdot \gw{2^{t-1}} \cdot \prod_{i = 1}^{t-1} \big( \gw{1} + \gw{d_{i}d_{i+1}} \big).
    \end{align*}
    The equality 
    \begin{equation} \label{eq-product_expanded}
        \prod_{i = 1}^{t-1} \big( \gw{1} + \gw{d_{i}d_{i+1}} \big) = \sum_{\substack{I\subset \{1,\dots,t\}\\ \vert I\vert\equiv 0\mod 2}} \qinv{\prod_{i\in I} d_i}.
    \end{equation}
    was already verified in the proof of \cref{lem-A1mult-twintree}. 
    Finally, the lemma follows by considering the parity of $m_\circ$: if $m_\circ$ is even, then $\gw{d_\circ^{m_\circ}} = \gw{1}$ and we are done. If $m_\circ$ is odd, then $\gw{d_\circ^{m_\circ}} = \gw{d_\circ}$ and multiplying this with the right hand side of \cref{eq-product_expanded} effectively changes the index sets $I$ which contain $\circ$ to $I \setminus \{\circ\}$ and the $I$ which do not contain $\circ$ to $I \cup \{\circ\}$. The result is a summation over all $J \subseteq \{1, \ldots, t\}$ of odd size, as predicted.
\end{proof}

In order to account for the remaining factors in $\mult^{\AA^1}(\cT)$ we consider the deformation patterns of the twin edges in $\cT$. To do so, we define the following tower of $k$-algebras:
\begin{enumerate}
    \item $F_{\cT}$ as determined in \cref{lem-algebratwintree}.
    \item $L_\cT = \bigotimes_{\text{twin edges } (e, e')} L_{e, e'}$ with $L_{e,e'}$ being the deformation patterns from \cref{prop-deformpatternandrefinedpoint} and the tensor product being taken over $F_\cT$.
    \item $M_\cT = \bigotimes_{\text{marked twin edges } (e, e')} M_{e, e'}$ with $M_{e, e'}$ being the refined condition to pass through a point from \cref{prop:algebrarefinedpointconditions} and the tensor product being taken over $L_\cT$.
\end{enumerate}

\begin{lemma} \label{lem-prooftwintreemult}
    Let $C_\cT$ be the algebraic lift of $\cT$ defined over $M_\cT$.
    Then
    \[ \Tr_{M_\cT / k} \big(\Wel(C_\cT) \big) = \mult^{\AA^1} (\cT). \]
\end{lemma}

\begin{proof}
    \[ \Wel(C_\cT) = \underbrace{\prod_{v \text{ vertex}} \Wel(C_v)}_{\in \GW(F_\cT)} \cdot \underbrace{\prod_{(e, e') \text{ twin edges}} \Wel(C_{e, e'})}_{\in \GW(L_\cT)}  \]
    \begin{align*}
        \Tr_{M_\cT / k} \big(\Wel(C_\cT) \big) &= \Tr_{F_\cT / k} \circ \Tr_{M_\cT / F_\cT} \big(\Wel(C_\cT) \big) \\
        &= \Tr_{F_\cT / k} \bigg(\prod_{v \text{ vertex}} \Wel(C_v) \cdot \Tr_{M_\cT / F_\cT} \Big(\prod_{(e, e')} \Wel(C_{e, e'}) \Big) \bigg)\\
        &\overset{\ref{prop:twinedgemult}}{=} \Tr_{F_\cT/k}\Big(\prod_{v\text{ vertex}}\Wel(C_v)\cdot \underbrace{\prod_{(e, e')}\mult^{\AA^1}(e,e')}_{\in \GW(k)}\Big)\\
        &=\prod_{(e, e')}\mult^{\AA^1}(e,e')\cdot \Tr_{F_\cT/k}\Big(\underbrace{\prod_{v\text{ vertex}}\Wel(C_v)}_{= \Wel(C_V)}\Big) \\
        &\overset{\ref{lem-tracetwintree}}{=} \prod_{(e, e')}\mult^{\AA^1}(e,e')\cdot \gw{2^{t-1}} \cdot \sum_{\substack{I\subset \{1,\dots,t\}\\ \vert I\vert\equiv m_{\circ}\mod 2}} \qinv{\prod_{i\in I} d_i} \\
        &\overset{\ref{def-multtwintree}}{=} \mult^{\AA^1}(\cT).
    \end{align*}
    If a twin edge $(e, e')$ in $\cT$ is unmarked, then it cannot be an elevator. Hence its weight is 1 and therefore $\mult^{\AA^1}(e, e') = \gw{1}$ by \cref{prop-deformpatternandrefinedpoint}. 
\end{proof}

\section{The correspondence theorem}
\label{sec-corres}

In this section we assemble the local computations from Sections~\ref{section:computations local pieces},~\ref{sec-edges}, and~\ref{sec-twintree} to prove our Correspondence Theorem \ref{thm-main} from the introduction, see \cref{thm:correspondence} below. However, before we present the proof, we return to the examples in \cref{fig-exfloordecomposed} and give examples for the algebras $F$, $L$, and $M$ which play the central role in the proof below. Heuristically, $F$ is the $k$-algebra defined by the different choices for the local pieces, $L$ is the $F$-algebra defined by the different choices of deformation patterns and $M$ is the $L$-algebra defined by the different choices for the refined point conditions. In particular, we have $k\subset F\subset L\subset M$. All these choices define one algebraic lift $C_M$ of the given tropical stable map $(\Gamma,f)$ over $M$ and 
the correspondence theorem \ref{thm:correspondence} basically computes $\Tr_{M/k}(\Wel(C_M))$ and shows that it is equal to the expression for $\mult^{\AA^1}(\Gamma,f)$ derived in \cref{cor-A1mult-twintrees}.

\begin{figure}[t]
    \centering
    \begin{tikzpicture}
	\tikzset{every path/.style={
		draw, line width=\edgewidth}}
	\path[double] (-1, 0) -- (0, 0) -- (0.7, 0.7);
	\fatpoint{-0.7, 0}
	\draw (-0.7, 0.1) node[anchor = south] {$\quadextension{d_1}$};
	\path[double] (0, 0) -- (0, -1);
	\fatpoint{0, -0.5}
	\draw (0.1, -0.5) node[anchor = west] {$\quadextension{d_2}$};
	\path (-1, -1) -- (0, -1) -- (1, -3) -- (2, -4) -- (3, -4) -- (4, -3);
	\fatpoint{1, -3}
	\draw (0.9, -3) node[anchor = north east] {$v_3, \quadextension{d_3}$};
	\path (1, -3) -- (1, -5);
	\path (2, -4) -- (2, -8.5);
	\path (3, -4) -- (3, -8.5);
	\path (-1, -5) -- (1, -5) -- (2.5, -6.5) -- (4, -5);
	\fatpoint{2, -6}
	\draw (1.9, -6) node[anchor = north east] {$v_4, \quadextension{d_4}$};
	\path[double] (2.5, -6.5) -- (2.5, -8.5);
	\fatpoint{2.5, -7.5}
	\draw (3.1, -7.5) node[anchor = west] {$\quadextension{d_5}$};
	\thinpoint{3, -8.2}
	
	\draw[thick, dashed] (0,-0.1) ellipse (1.1 and 1.3);
	\draw (1, 1) node {$\mathcal{T}_1$};
	\draw[thick, dashed] (2.5,-7.5) ellipse (0.4 and 1.5);
	\draw (3, -9) node {$\mathcal{T}_2$};
\end{tikzpicture}
\begin{tikzpicture}
	\tikzset{every path/.style={
			draw, line width=\edgewidth}}
	\path (-1, 0) -- (0, 0) -- (1, -2) -- (1.5, -3.5) -- (4, -3.5);
	\path[double] (0, 1) -- (0, 0);
	\fatpoint{0, 0.5}
	\draw (-0.1, 0.5) node[anchor = east] {$\quadextension{d_1}$};
	\path (1, 1) -- (1, -2);
	\thinpoint{1, -1}
	\path (1.5, -3.5) -- (1.5, -5);
	\fatpoint{1.5, -3.5}
	\draw (1.4, -3.5) node[anchor = east] {$v_2, \quadextension{d_2}$};
	\draw (1.4, -4.25) node[anchor = east] {$e, \edgeweight{3}$};
	\path (-1, -5) -- (1.5, -5) -- (2, -6.5) -- (3, -7.5) -- (4, -7.5); 
	\fatpoint{2, -6.5}
	\draw (1.9, -6.5) node[anchor = east] {$v_3, \quadextension{d_3}$};
	\path (2, -6.5) -- (2, -8.5);
	\path[double] (3, -7.5) -- (3, -8.5);
	\fatpoint{3, -8}
	\draw (3.1, -8) node[anchor = west] {$\quadextension{d_4}$};
	
	\draw[thick, dashed] (0, 0.4) ellipse (0.7 and 0.8);
	\draw (0.5, 1.4) node {$\mathcal{T}_1$};
	
	\draw[thick, dashed] (3, -8) ellipse (0.7 and 0.8);
	\draw (4, -8.5) node {$\mathcal{T}_2$};
	
\end{tikzpicture}
\begin{tikzpicture}
	\tikzset{every path/.style={
			draw, line width=\edgewidth}}
	\path (-1, 0) -- (0, 0) -- (2, -2) -- (4, -2);
	\path (0, 1) -- (0, 0);
	\thinpoint{0, 0.5}
	\path (1, 1) -- (1, -3.5);
	\fatpoint{1, -1}
	\draw (0.9, -1) node[anchor = east] {$v_1, \quadextension{d_1}$};
	\path (-1, -3.5) -- (1, -3.5) -- (2, -4.5) -- (2.5, -6) -- (3.5, -8) -- (4.5, -8);
	\path (2 - \doubleedgesep, -2) -- (2 - \doubleedgesep, -4.5);
	\path (2 + \doubleedgesep, 1) -- (2 + \doubleedgesep, -4.5);
	\fatpoint{2, -3.25}
	\draw (2.1, -1.9) node[anchor = south west] {$v'_2$};
	\draw (2.1, -3.25) node[anchor = west] {$\quadextension{d_2}$};
	\draw (2.1, -4.5) node[anchor = west] {$v_2$};
	\path (2.5, -6) -- (2.5, -9);
	\fatpoint{2.5, -6}
	\draw (2.6, -6) node[anchor = west] {$v_3, \quadextension{d_3}$};
	\path[double] (3.5, -8) -- (3.5, -9);
	\fatpoint{3.5, -8.5};
	\draw (3.6, -8.5) node[anchor = west] {$\quadextension{d_4}$};
	
	\draw[thick, dashed] (3.5, -8.5) ellipse (0.7 and 0.8);
	\draw (4, -7.5) node {$\mathcal{T}$};
\end{tikzpicture}
    \caption{The tropical curves from \cref{fig-exfloordecomposed} with additional notation.}
    \label{fig-exalgebras}
\end{figure}
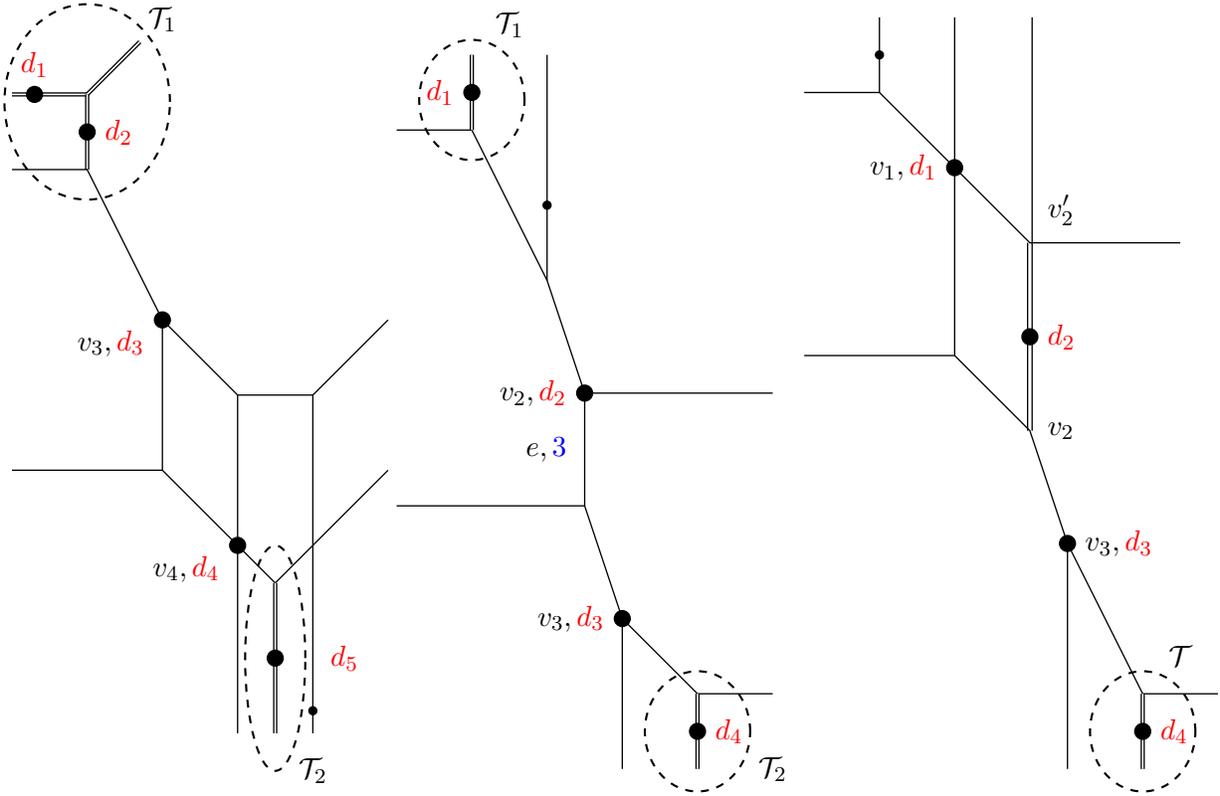

\begin{example} 
    Consider the first tropical curve form \cref{fig-exfloordecomposed} with the additional notation as in \cref{fig-exalgebras}. The algebras locally at the twin trees and vertices are as follows:
    \begin{align*}
        F_{\cT_1} &= \quotient{k[\iota_{1,2}]}{ (\iota_{1,2}^2 - d_1d_2) } & &\text{\cref{lem-algebratwintree}} \\
        F_{\cT_2} &= k & &\text{\cref{lem-algebratwintree}} \\
        F_{v_3} &= \operatorname{Fix} \left( \quotient{k[\iota, u]}{(\iota^2 - d_3, u - \eta/\conj{\eta})} \right) \cong k & &\text{\cref{lem:fatpoint_thintriangle}}\\
        F_{v_4} &= \quotient{k[\iota_4]}{(\iota_4^2 - d_4)} & &\text{\cref{lem:par5}}
    \end{align*}
    and all other $F_v = k$. In total we see that the algebra $F$ for this curve is 
    $$ F = \quotient{k[\iota_{1, 2}, \iota_4]}{(\iota_{1,2}^2 - d_1d_2, \iota_4^2 - d_4)}. $$
    Since all edges are of weight 1 in this curve, all deformation patterns and refined conditions to pass through points are trivial, i.e. $M = L = F$. The quadratic weight of a lift of the tropical curve is a product over twin trees and vertices. By \cref{lem-WelCT} we get the trace of $\gw{d_1} = \gw{d_2}$ as a factor and the contributions from every other vertex are trivial. The result is 
    \[ \Tr_{M / k} (\gw{d_1}) = \Tr_{F, k}(\gw{d_1} )= \beta_4 \cdot \big( \gw{2d_1} + \gw{2d_2} \big), \]
    i.e. this is precisely the multiplicity which was computed in \cref{fig-exfloordecomposed}.
    
    For the second curve the local algebras are as follows: 
    \begin{align*}
        F_{\cT_1} &= F_{\cT_2} = k  & &\text{\cref{lem-algebratwintree}}\\
        F_{v_2} &= \operatorname{Fix} \left( \quotient{k[\iota_2, u]}{(\iota_2^2 - d_2, u^3 - \eta/\conj{\eta})} \right) & &\text{\cref{lem:fatpoint_thintriangle}}\\
        F_{v_3} &= k & &\text{\cref{lem:fatpoint_thintriangle}}
    \end{align*}
    and the other $F_v = k$. Hence in total we have that $F = F_{v_2}$ is a rank 3 algebra over $k$. Since $e$ is the only edge with weight larger than 1, we get $L = L_e = F[u] / (u^3 - 1)$ and $M = L$ since $e$ is unmarked. For the quadratic weight of the tropical curve we obtain
    \[ \Tr_{M/ k} (\gw{1}) = 3^{\AA^1} \gamma(3, d_2) . \]

    For the third curve we have the following local algebras:
    \begin{align*}
        F_\cT &= k  & &\text{\cref{lem-algebratwintree}} \\
        F_{v_1} &= \quotient{k[\iota_1]}{(\iota_1^2 - d_1)} & &\text{\cref{lem:par5}} \\
        F_{v_2} &= k & &\text{\cref{lem:type_C}} \\
        F_{v'_2} &= \quotient{k[\iota_2]}{(\iota_2^2 - d_2)} & &\text{\cref{lem:trapezoid}}\\
        F_{v_3} &= k & &\text{\cref{lem:fatpoint_thintriangle}}  
    \end{align*}
    and all other $F_v = k$. In total we have 
    $$F = \quotient{k[\iota_1, \iota_2]}{(\iota_1^2 - d_1, \iota_2^2 - d_2)}$$
    and since all edge weights are 1, we have $M = L = F$. In conclusion, the quadratic weight for the tropical curve is 
    \[ \Tr_{M/k} (\gw{1}) = \beta_1 \beta_2. \]
\end{example}

We are now ready to prove the main result of this paper.

\begin{theorem}[Correspondence theorem for vertically stretched point conditions]\label{thm:correspondence}
    Assume $k$ is a perfect field of characteristic~$0$ or greater than the diameter of $\Delta$. Furthermore, assume that the characteristic is not~$2$ or $3$. Let $K=\Puiseux{k}$. Assume that $S$ is the toric del Pezzo surface over $K$ associated with a given $\Delta$. Let $r,s\in\NN$ such that $r+2s=\#(\partial\Delta\cap \ZZ^2)-1$. 

    Let $\mathcal{P}$ be a point configuration consisting of $r$ $K$-rational points and $s$ points with residue field $K(\sqrt{d_1})$, $\ldots$, $K(\sqrt{d_s})$ which tropicalizes to a generic vertically stretched configuration of $r$ simple and $s$ double points.
    Let $(\Gamma,f)$ be a degree-$\Delta$ rational tropical stable map passing through this configuration. Then the multiplicity $\mult^{\AA^1}(\Gamma,f)$ is equal to the sum of traces of quadratic weights of all log stable curves tropicalizing to $(\Gamma,f)$ which pass through $\mathcal{P}$.
    More precisely,
    \[\mult^{\AA^1}(\Gamma,f)=\sum_{\substack{\text{$u$ through $\mathcal{P}$} \\ \trop(u)=(\Gamma,f)} }\Tr_{\kappa(u)/K}(\Wel(u))\in\GW(K)\cong\GW(k)\]
    where $\kappa(u)$ is the field of definition of $u$.

    From this we deduce that $$N^{\AA^1}_\Delta(r,(d_1,\ldots,d_s))=N^{\AA^1,\trop}_\Delta(r,(d_1,\ldots,d_s)):=\sum_{(\Gamma,f)}\mult^{\AA^1}(\Gamma,f),$$
    in $\GW(k)$
    where the sum runs over all $(\Gamma,f)$ passing through a fixed generic vertically stretched configuration of $r$ simple and $s$ double points.
\end{theorem}
\begin{proof}
    
For the next paragraph we add the field or algebra $A$ over which we work to the notation of $N^{\AA^1}_{\Delta}(r,(d_1,\ldots,d_s);A)$ to indicate that we are working in $\GW(A)$.

    The field of Puiseux series $K=\Puiseux{k}$ has the same characteristic as $k$ and it is perfect. In particular, $N^{\A^1}_\Delta(r,(d_1,\ldots,d_s);K)$ is independent of the chosen point configuration by \cref{thm:KLSW}. Hence, we can choose a point configuration that specializes to a generic point configuration in $k$. Then under the canonical isomorphism $\GW(K)\cong \GW(k)$, we have that $N^{\A^1}_\Delta(r,(d_1,\ldots,d_s);K)$ is mapped to $N^{\A^1}_\Delta(r,(d_1,\ldots,d_s);k)$ by Lemma \ref{lm:GWPuiseux}. 
    
    Consequently, to show that
    $N^{\A^1}_\Delta(r,(d_1,\ldots,d_s);k)=N^\trop_\Delta(r,(d_1,\ldots,d_s);k)$, it suffices to show that $N^{\A^1}_\Delta(r,(d_1,\ldots,d_s);K)=N^\trop_\Delta(r,(d_1,\ldots,d_s);k)$.

    To show this, we choose a point configuration of $r$ points defined over $K$ and $s$ pairs of Galois conjugate points defined over $K(\sqrt{d_i})$ for $i=1,\ldots,s$ which tropicalizes to a vertically stretched point configuration of $r$ simple and $s$ double points. 
    We find all log stable maps defined over $\Puiseux{\overline{k}}$ through this point configuration that tropicalize to the given tropical stable map $(\Gamma,f)$. For this we generalize \cite{Shu06b}. 
    The first step is to find all stable maps on the local pieces $\Delta_v$ in the dual subdivision of the tropical curve, as computed in \cref{section:computations local pieces}. 
 
    By Theorem \ref{theorem:multofatropstablemap}, the different choices for the log stable maps extending the local pieces are determined by the different choices for the deformation patterns (Proposition \ref{prop-deformpatternandrefinedpoint}) and the refined point conditions (Proposition \ref{prop:algebrarefinedpointconditions}). 
    For vertically stretched point configurations only the vertical edges can have weight greater than $1$ and thus, they are the only ones contributing.

    The different choices of log stable maps define a finite étale $k$-algebra $M$.
    This $k$-algebra $M$ is isomorphic to a product $M\cong L_1\times \cdots \times L_\lambda$ where the $L_i$ are finite separable field extensions of $k$. This means that there are $\lambda$ log stable maps $C_1,\ldots,C_\lambda$ with fields of definitions $\Puiseux{L_1},\ldots,\Puiseux{L_\lambda}$, respectively. The statement of the theorem follows once we show that the multiplicity $\mult^{\AA^1}(\Gamma,f)$ is mapped to
    \[\sum_{i=1}^\lambda \Tr_{\Puiseux{L_i}/\Puiseux{k}}\big(\Wel(C_i)\big)\]
    under the canonical isomorphism $\GW(k)\cong \GW(K)$
    for the given tropical stable map $(\Gamma,f)$ through the given point configuration.
    This equals 
    \[\Tr_{\Puiseux{M}/\Puiseux{k}}\big(\Wel(C_M) \big)\]
    with $C_M=C_1\times\cdots\times C_\lambda$ defined over $\Puiseux{M}$.
    By Theorem \ref{theorem:multofatropstablemap},
    $\Wel(C_M)\in \GW \big(\Puiseux{M} \big)$ maps to 
    $$\prod_{\text{vertices }v}\Wel(C_v)\cdot \prod_{\text{edges }e}\Wel(C_{e})\cdot \prod_{\text{twin edges }(e,e')}\Wel(C_{e,e'})\in \GW(M)$$ 
    under the canonical isomorphism $\GW \big(\Puiseux{M} \big)\cong \GW(M)$ where the $C_v$ are the local pieces corresponding to the vertices $v$ of $f(\Gamma)$ and the $C_e$ and $C_{e,e'}$ the deformation patterns corresponding to the non-twin edges $e$ and twin edges $(e,e')$ of $\Gamma$, respectively. 
    
    Therefore, 
    $\Tr_{\Puiseux{M}/\Puiseux{k}} \big(\Wel(C_M) \big)\in \GW \big(\Puiseux{k} \big)$  can be computed as 
\begin{equation}
    \label{eq:computationquadraticweight}
\Tr_{M/k}\left(\prod_{\text{vertices }v}\Wel(C_v)\cdot \prod_{\text{edges }e}\Wel(C_{e})\cdot\prod_{\text{twin edges }e,e'}\Wel(C_{e,e'})\right)\in \GW(k).
    \end{equation}
     The goal is to identify this expression with the formula for the quadratically enriched multiplicity $\mult^{\AA^1}(\Gamma,f)$ derived in \cref{cor-A1mult-twintrees}.

Recall that the algebra $M$ is constructed in two steps. First, we consider the $k$-algebras $F_v$ corresponding to the local pieces associated with the vertices $v$ of the image plane tropical curve~$f(\Gamma)$. These algebras combine to form a finite étale $k$-algebra $F$. 
It is important to note that $F$ is not simply the tensor product over the $F_v$. This is because the point conditions may appear at multiple vertices $v$, but it is unnecessary to adjoin the corresponding variable more than once. 
In Tables~\ref{tab:possibilities} and~\ref{tab:possibilities2}, we list all possible $k$-algebras $F_v$ for the different types of vertices $v$. From this table, we deduce that $F$ can be obtained from $k$ as follows:
\begin{itemize}
    \item For a vertex $v$ of type \fatPointOnVertex{} the $k$-algebra $F_v$ is the fixed subalgebra $L^\phi$ as in Table \ref{tab:possibilities}. By the primitive element theorem, there exist an element $w_v\in L$ such that $L^\phi=k(w_v)$. For each vertex $v$ of type \fatPointOnVertex{}, we adjoin the element $w_v$ to $k$.
    \item Let $R=\{i:  d_i\in \fourValentVertex\cup  \triangleWithMergedEdge{}\cup \fourValentVertexWithMergedEdge{}{}\cup\triangleWithTwoMergedEdges{} \cup \mergedTriangles{} \}$. 
    For each $i\in R$ we have to adjoin $\iota_i$ to $k$ satisfying $\iota_i^2=d_i$.
    \item Let $\mathcal{T}$ be a twin tree with double points corresponding to $d_1,\ldots,d_t$. Then
    $F$ contains $\iota_{12},\ldots,\iota_{t-1,t}$ for variables $\iota_{i,i+1}$ satisfying $\iota_{i,i+1}^2=d_id_{i+1}$ for each twin tree $\mathcal{T}$. In other words $F$ contains $F_\mathcal{T}$ for each twin tree $\mathcal{T}$, as in \cref{lem-algebratwintree}.
\end{itemize}
Then, one obtains $M$ from $F$ from the different choices for the deformations patterns and refined point conditions.
It follows from Proposition \ref{prop-deformpatternandrefinedpoint} and Proposition \ref{prop:algebrarefinedpointconditions}
that
\[M=\bigotimes_{e \text{ marked}}M_e\otimes\bigotimes_{e,e' \text{ marked}}M_{e,e'}\otimes\bigotimes_{e \text{ not marked}}L_e\otimes\bigotimes_{e,e' \text{ not marked}}L_{e,e'}\]
where all the tensor products are over $F$.

In particular, we have that $k\subset F\subset M$ and thus, \cref{eq:computationquadraticweight} equals
\begin{align*}
\begin{split}
     \Tr_{F/k}\circ \Tr_{M/F}\left(\prod_{\text{vertices }v}\Wel(C_v)\cdot \prod_{\text{edges }e}\Wel(C_e)\cdot\prod_{\text{twin edges }e,e'}\Wel(C_{e,e'})\right). 
    \end{split}
\end{align*}

Since each factor $\Wel(C_v)$ belongs to $\GW(k)$ (cf. \cref{tab:possibilities,tab:possibilities2}), the relation \cref{eq:traceofsomethingink} in Proposition \ref{prop:traces} implies that this trace form is equal to
\begin{align*}
        \left(\prod_{\text{vertices }v} \Wel(C_v)\right)\cdot \left(\Tr_{F/k}\circ \Tr_{M/F}\left(\prod_{\text{edges }e}\Wel(C_e)\cdot\prod_{\text{twin edges }e,e'}\Wel(C_{e,e'})\right)\right).
\end{align*}

 Item \ref{eq:traceoftensorproduct} in Proposition \ref{prop:traces} implies that
\begin{align*}
    &\Tr_{M/F}\left( \prod_{\text{edges }e}\Wel(C_e)\cdot\prod_{\text{twin edges }(e,e')}\Wel(C_{e,e'})\right)\\
    &\qquad=\prod_{e \text{ marked}}\Tr_{M_e/F} \big(\Wel(C_e) \big)\cdot\prod_{(e,e') \text{ marked}}\Tr_{M_{e,e'}/F} \big(\Wel(C_{e,e'}) \big)\\
    &\qquad\cdot\prod_{e \text{ not marked}}\Tr_{L_e/F} \big(\Wel(C_e) \big)\cdot\prod_{(e,e') \text{ not marked}}\Tr_{L_{e,e'}/F} \big(\Wel(C_{e,e'}) \big).
\end{align*}

We compute these factors separately. In \cref{prop:nontwinedgemult,prop:twinedgemult}, we proved that for a marked edge~$e$ and for a marked twin edge~$(e,e')$, respectively, the trace forms equal
\[
    \Tr_{M_e/F} \big(\Wel(C_e) \big)=\mult^{\AA^1}(e),
    \qquad\text{ and }\qquad
    \Tr_{M_{e,e'}/F} \big(\Wel(C_{e,e'}) \big)=\operatorname{mult}^{\AA^1}(e,e').
\]

Let us remark that edges of weight one contribute trivially to this computation. Indeed, for a non-twin edge of weight one we have $L_{e}=F$ and $\Wel(C_e)=\gw{1}$. 
For a twin-edge of weight one we have $L_{e,e'}=F$ and $\Wel(C_{e,e'})=\gw{1}$. 
For a vertically stretched point configuration, this already covers almost all unmarked (twin-)edges. The only unmarked edges of weight greater than one are the vertical non-twin edges adjacent to a vertex of type \fatPointOnVertex{}. For such a non-marked non-twin edge $e$ of weight $m$ one has that
\begin{equation*}
    \Tr_{L_e/F} \big(\Wel(C_e) \big)=\begin{cases}
        \Tr_{L_e/F}(\gw{1})=\gw{m}+\frac{m-1}{2}\h & \text{for $m$ odd,}\\
        \Tr_{L_e/F}(\gw{-u})=\frac{m}{2}\h & \text{for $m$ even.}
    \end{cases}
\end{equation*}
by Proposition \ref{prop:traces} 
and Proposition \ref{prop-deformpatternandrefinedpoint}.
Note that this equals $m^{\AA^1}$ as defined in Notation \ref{not-vertexedgeweight}. 
Hence, we get that
\begin{align*}
   &\Tr_{M/F}\left(\prod_{\text{twin edges $(e,e')$}}\Wel(C_{e,e'})\cdot \prod_{\text{edges $e$}}\Wel(C_e)\right)\\
   &\qquad=\prod_{(e,e')}\mult^{\AA^1}(e,e')\cdot \prod_{e\in E_{\neg\fatPointOnVertex{}}} \mult^{\AA^1}(e)\cdot \prod_{e\in E_{ \fatPointOnVertex{}}}m_e^{\AA^1}\in\GW(F)
\end{align*}
where $E_{\fatPointOnVertex{}}$ (resp., $E_{\neg\fatPointOnVertex{}}$) is the set of non-twin edges adjacent (resp., non-adjacent) to a vertex of type~$\fatPointOnVertex{}$ and in the last product $m_e$ is the weight of $e$.
Using the relation in \cref{eq:traceofsomethingink} of \cref{prop:traces}, it follows that the trace form in \cref{eq:computationquadraticweight} equals
\[\left(\prod_{\text{vertices } v}\Wel(C_v)\right)\cdot \left(\prod_{\text{twin edges }(e,e')}\mult^{\AA^1}(e,e')\cdot \prod_{e\in E_{\neg\fatPointOnVertex{}}} \mult^{\AA^1}(e)\cdot \prod_{e\in E_{\fatPointOnVertex{}}}m_e^{\AA^1} \right)\cdot \Tr_{F/k}(\gw{1}).\]
Next we compute $\Tr_{F/k}(\gw{1})$.
We have that the trace form
\begin{align*}
    \begin{split}
        \Tr_{F/k}(\gw{1})&= \prod_{v\in V_{\fatPointOnVertex{}}} \Tr_{F_v/k}(\gw{1}) \cdot \prod_{\text{twin trees $\mathcal{T}$}} \Tr_{F_\mathcal{T}/k}(\gw{1})\cdot \prod_{i\in R}\Tr_{k[\iota_i]/k}(\gw{1}).
    \end{split}
\end{align*}
Recall from \cref{tab:possibilities} that if $v$ does not belong to a twin tree then $\Wel(C_v)=\gw{1}$. It follows that that 
\[\prod_{\text{vertices } v}\Wel(C_v)=\prod_{\text{twin trees }\mathcal{T}}\Wel(C_\mathcal{T})\overset{\cref{lem-tracetwintree}}{=}\prod_{\text{twin tree $\mathcal{T}$}}\gw{(d_\circ({\mathcal{T})})^{m_o(\mathcal{T})}}\]
where $m_\circ(\mathcal{T})$ is the associated weight  of $\mathcal{T}$ and $d_\circ(\mathcal{T})$ is one of its markings.
Substituting this equality in  \cref{eq:computationquadraticweight} yields
\begin{equation} \label{eq:finaleqcorrespondence}
    \begin{aligned}
        &\prod_{\text{twin tree $\mathcal{T}$}}\gw{(d_o(\mathcal{T}))^{m_o(\mathcal{T})}}\cdot \prod_{v\in V_{\fatPointOnVertex{}}} \Tr_{F_v/k}(\gw{1}) \cdot \prod_{\text{twin trees $\mathcal{T}$}} \Tr_{F_\mathcal{T}/k}(\gw{1})\cdot \prod_{R}\Tr_{k[\iota_i]/k}(\gw{1})\\
        \cdot &\prod_{\text{twin edges }e,e'}\mult^{\AA^1}(e,e')\cdot \prod_{e\in  E_{\neg\fatPointOnVertex{}}} \mult^{\AA^1}(e)\cdot \prod_{e\in E_{\fatPointOnVertex{}}}m_e^{\AA^1}.
    \end{aligned}
\end{equation}

Now we are ready to identify the factors in \eqref{eq:finaleqcorrespondence} with the factors of $\mult^{\AA^1}(\Gamma,f)$. Recall from \cref{rem-A1mult-twintree} that $\mult^{\AA^1}(\Gamma,f)$ can be written as
    \[\mult^{\AA^1}(\Gamma,f)=\prod_{i=1}^T
\mult^{\AA^1}(\mathcal{T}_i)\cdot \prod_{v \in V_{\fatPointOnVertex{}}} \gamma(v) m^{\AA^1} \cdot \prod_{e \in E_{\neg\fatPointOnVertex{}}} \mult^{\AA^1}(e) \cdot \prod_{i\in R} \beta_i.\]
\begin{itemize}
    \item For every twin tree $\mathcal{T}$, we have that \cref{lem-prooftwintreemult} implies
    \[
        \gw{(d_o(\mathcal{T}))^{m_o(\mathcal{T})}}\cdot \Tr_{F_\mathcal{T}/k}(\gw{1})\cdot  \prod_{e,e'\in \mathcal{T}}\mult^{\AA^1}(e,e')=\mult^{\AA^1}(\mathcal{T})
    \]
    where $\mult^{\AA^1}(\mathcal{T})$ is as in \cref{def-multtwintree}.
    
    \item For each vertex $v$ of type \fatPointOnVertex{} there is a corresponding $e \in E_\fatPointOnVertex{}$ and we combine their corresponding factors to find $\gamma(m_v, d_v) m^{\AA^1}$. Indeed, if $e$ is of odd weight $m_v$ and the double point condition on $v$ is defined over $k(\sqrt{d_v})$, we have that
    \begin{align*}
        &m_v^{\AA^1}\cdot \Tr_{F_v/k}(\gw{1}) \\
        &\mkern-2mu\overset{\ref{lemma:algebratypeA}}{=}{} \left(\gw{m_v} + \frac{m_v-1}{2}\h\right)\cdot \left(\gw{m_v}+\frac{m_v-1}{2}\left(\gw{2m_v}+\gw{-2d_vm_v}\right)\right) \\
        &={} \gamma(m_v, d_v) m_v^{\AA^1}.
    \end{align*}
    If $m_v$ is even then
    \[m_v^{\AA^1}\cdot \Tr_{F_v/k}(\gw{1})=\frac{m_v}{2}\h\cdot \dim_kF_v= \gamma(m_v, d_v) m_v^{\AA^1}.\]
    
    \item For $e \in E_{\neg\fatPointOnVertex{}}$, the necessary factor $\mult^{\AA^1}(e)$ already appears in \cref{eq:finaleqcorrespondence}.

    \item For each $i\in R$, the trace form
    \[\Tr_{k[\iota_i]/k}(\gw{1})=\beta_i.\]
\end{itemize}
\end{proof}

\section{Floor diagrams}\label{sec-floor}

Floor diagrams are an efficient combinatorial tool to compute tropical counts of curves (see e.g.\ \cite{BM08, Bru15, FM09}). In \cite{JPMPR23} these techniques have been modified to handle the quadratically enriched count with $k$-rational point conditions as defined in \cite{JPP23}. In this section we describe the adaptations needed in order to use floor diagram techniques to compute $N^{\trop}_\Delta(r, (d_1, \ldots, d_s))$. The main result of this section is \cref{thm-flooreqtrop_intro} from the introduction, see \cref{thm-flooreqtrop} below. The floor diagram approach was used to obtain the computational results which we present in \cref{sec-computations}.

Fix a smooth del Pezzo degree $\Delta$ as in Definition \ref{def-deg}. First, we consider the count of rational tropical stable maps satisfying $n(\Delta)$ vertically stretched simple point conditions, i.e.\ we do not have double points for the moment. This count (no matter whether we count with the complex, real or quadratically enriched multiplicity) can be shown to equal the count of floor diagrams as follows. 

Recall that if we count tropical stable maps with the real multiplicity, this equals the Welschinger invariant, i.e.\ the signed count of rational curves of degree $\Delta$ satisfying $n(\Delta)$ \emph{real} point conditions. If we count with quadratically enriched multiplicity, this equals the quadratically enriched count of rational curves of degree $\Delta$ satisfying $n(\Delta)$ \emph{$k$-rational} point conditions. 

\begin{definition}\label{def-floordiagram}
    A rational \emph{floor diagram} of degree $\Delta$ is a bipartite tree with vertices of colors black and white, such that there are $n(\Delta)$ vertices and they are linearly ordered. Moreover, the tree may also have ends and the edges are weighted. All ends are required to be of weight $1$. By the linear order on the vertices, we can view the edges as being oriented from the smaller to the bigger vertex and similarly, the ends are oriented and can thus be viewed as incoming resp. outgoing ends.
    
    The \emph{divergence} of a vertex is the sum of the weights of the incoming edges minus the sum of the weight of the outgoing edges.
    We require that black vertices are always of divergence $0$.
    Vertices may be decorated with a set of integer numbers that we call \emph{leaks}, which sum up to the divergence.
    Depending on $\Delta$, we require a suitable condition on the divergence and on the leaks to be satisfied at every white vertex. These conditions are as follows. 

    \begin{enumerate}
    
        \item If $\Delta$ corresponds to degree $d$ in $\mathbb{P}^2$, we have $d$ white vertices which are all of divergence $1$. (Also all leaks are required to be $1$, which makes it irrelevant to keep track of them.)
        \item If $\Delta$ corresponds to bidegree $(a_1,a_2)$ in $\mathbb{P}^1\times \mathbb{P}^1$, we have $a_1$ white vertices which are all of divergence $0$. (All leaks are $0$.)
        \item If $\Delta$ corresponds to the blowup of the plane in a point as the third polygon in Figure \ref{fig-degrees}, then $a_1$ white vertices are of divergence $0$ and the remaining $d-a_1$ white vertices are of divergence~$1$. (Leaks are also either $1$ or $0$.)
        \item If $\Delta$ corresponds to the blowup of the plane in two points as the fourth polygon in Figure~\ref{fig-degrees}, we can assume without restriction that $a_1\geq a_2$. Then for the $d$ white vertices we must have a total of $a_2$ leaks of $-1$, $a_1-a_2$ leak $0$ and  $d-a_1$ leak $1$. (Here it is possible to combine a leak of $-1$ and a leak of $1$ to one white vertex of total divergence $0$, but one must keep track of the two canceling leaks.) 
        \item If $\Delta$ corresponds to the blowup of the plane in three points as the fifth polygon in Figure~\ref{fig-degrees}, we can assume without restriction that $a_1\geq a_2$. Then for the $d-a_3$ white vertices, we must have $a_2$ leaks of $-1$, $a_1-a_2$ leaks $0$ and  $d-a_1-a_3$ leaks of $1$. (As before,  leaks can be combined and it is therefore necessary to keep track of them.)
    \end{enumerate}
    Furthermore, there are conditions on the ends of a floor diagram which depend on $\Delta$ as follows. 

    \begin{enumerate}
        \item If $\Delta$ corresponds to degree $d$ in $\mathbb{P}^2$, we have $d$ incoming ends and no outgoing ends.
        \item If $\Delta$ corresponds to bidegree $(a_1,a_2)$ in $\mathbb{P}^1\times \mathbb{P}^1$, we have $a_2$ incoming and $a_2$ outgoing ends.
        \item If $\Delta$ corresponds to the blowup of the plane in a point as the third polygon in Figure \ref{fig-degrees}, we have $d-a_1$ incoming ends and no outgoing end.
        \item If $\Delta$ corresponds to the blowup of the plane in two points as the fourth polygon in Figure~\ref{fig-degrees}, and $a_1\geq a_2$, we have $d-a_1-a_2$ incoming ends and no outgoing ends.
        \item If $\Delta$ corresponds to the blowup of the plane in three points as the fifth polygon in Figure~\ref{fig-degrees}, and  $a_1\geq a_2$, we have  $d-a_1-a_2$ incoming ends and $a_3$ outgoing ends.
        \end{enumerate}
\end{definition}

\begin{example}
Figure \ref{fig-floordiagramnormal} depicts a rational floor diagram of degree $4$. Even though we use vertically stretched point conditions for our tropical stable maps, we choose to depict the vertices of this floor diagram in a horizontal arrangement, to better fit the picture.

    \begin{figure}[t]
        \centering

\tikzset{every picture/.style={line width=0.75pt}} 

\begin{tikzpicture}[x=0.75pt,y=0.75pt,yscale=-1,xscale=1]

\draw   (200,226.58) .. controls (200,224.7) and (201.53,223.17) .. (203.42,223.17) .. controls (205.3,223.17) and (206.83,224.7) .. (206.83,226.58) .. controls (206.83,228.47) and (205.3,230) .. (203.42,230) .. controls (201.53,230) and (200,228.47) .. (200,226.58) -- cycle ;
\draw  [fill={rgb, 255:red, 0; green, 0; blue, 0 }  ,fill opacity=1 ] (140,226.58) .. controls (140,224.7) and (141.53,223.17) .. (143.42,223.17) .. controls (145.3,223.17) and (146.83,224.7) .. (146.83,226.58) .. controls (146.83,228.47) and (145.3,230) .. (143.42,230) .. controls (141.53,230) and (140,228.47) .. (140,226.58) -- cycle ;
\draw  [fill={rgb, 255:red, 0; green, 0; blue, 0 }  ,fill opacity=1 ] (113.17,226.58) .. controls (113.17,224.7) and (114.7,223.17) .. (116.58,223.17) .. controls (118.47,223.17) and (120,224.7) .. (120,226.58) .. controls (120,228.47) and (118.47,230) .. (116.58,230) .. controls (114.7,230) and (113.17,228.47) .. (113.17,226.58) -- cycle ;
\draw  [fill={rgb, 255:red, 0; green, 0; blue, 0 }  ,fill opacity=1 ] (170,226.58) .. controls (170,224.7) and (171.53,223.17) .. (173.42,223.17) .. controls (175.3,223.17) and (176.83,224.7) .. (176.83,226.58) .. controls (176.83,228.47) and (175.3,230) .. (173.42,230) .. controls (171.53,230) and (170,228.47) .. (170,226.58) -- cycle ;
\draw  [fill={rgb, 255:red, 0; green, 0; blue, 0 }  ,fill opacity=1 ] (230,226.58) .. controls (230,224.7) and (231.53,223.17) .. (233.42,223.17) .. controls (235.3,223.17) and (236.83,224.7) .. (236.83,226.58) .. controls (236.83,228.47) and (235.3,230) .. (233.42,230) .. controls (231.53,230) and (230,228.47) .. (230,226.58) -- cycle ;
\draw  [fill={rgb, 255:red, 0; green, 0; blue, 0 }  ,fill opacity=1 ] (260,226.58) .. controls (260,224.7) and (261.53,223.17) .. (263.42,223.17) .. controls (265.3,223.17) and (266.83,224.7) .. (266.83,226.58) .. controls (266.83,228.47) and (265.3,230) .. (263.42,230) .. controls (261.53,230) and (260,228.47) .. (260,226.58) -- cycle ;
\draw   (290,226.58) .. controls (290,224.7) and (291.53,223.17) .. (293.42,223.17) .. controls (295.3,223.17) and (296.83,224.7) .. (296.83,226.58) .. controls (296.83,228.47) and (295.3,230) .. (293.42,230) .. controls (291.53,230) and (290,228.47) .. (290,226.58) -- cycle ;
\draw   (380,226.58) .. controls (380,224.7) and (381.53,223.17) .. (383.42,223.17) .. controls (385.3,223.17) and (386.83,224.7) .. (386.83,226.58) .. controls (386.83,228.47) and (385.3,230) .. (383.42,230) .. controls (381.53,230) and (380,228.47) .. (380,226.58) -- cycle ;
\draw  [fill={rgb, 255:red, 0; green, 0; blue, 0 }  ,fill opacity=1 ] (320,226.58) .. controls (320,224.7) and (321.53,223.17) .. (323.42,223.17) .. controls (325.3,223.17) and (326.83,224.7) .. (326.83,226.58) .. controls (326.83,228.47) and (325.3,230) .. (323.42,230) .. controls (321.53,230) and (320,228.47) .. (320,226.58) -- cycle ;
\draw  [fill={rgb, 255:red, 0; green, 0; blue, 0 }  ,fill opacity=1 ] (350,226.58) .. controls (350,224.7) and (351.53,223.17) .. (353.42,223.17) .. controls (355.3,223.17) and (356.83,224.7) .. (356.83,226.58) .. controls (356.83,228.47) and (355.3,230) .. (353.42,230) .. controls (351.53,230) and (350,228.47) .. (350,226.58) -- cycle ;
\draw   (410,226.58) .. controls (410,224.7) and (411.53,223.17) .. (413.42,223.17) .. controls (415.3,223.17) and (416.83,224.7) .. (416.83,226.58) .. controls (416.83,228.47) and (415.3,230) .. (413.42,230) .. controls (411.53,230) and (410,228.47) .. (410,226.58) -- cycle ;
\draw    (176.83,225.08) -- (200,225.08) ;
\draw    (102.5,227.59) -- (116.5,227.58) ;
\draw    (93.33,227.34) .. controls (132.5,197.09) and (153.25,253.84) .. (200.25,227.97) ;
\draw    (82.5,224.84) .. controls (110.17,207.34) and (137,212.84) .. (173.42,226.58) ;
\draw    (117.75,227.34) .. controls (163.5,251.34) and (235,241.84) .. (291.5,228.84) ;
\draw    (232,226.34) .. controls (254.5,227.84) and (255.75,207.09) .. (290.58,225.42) ;
\draw    (79.75,218.09) .. controls (92.75,203.59) and (201.42,209.83) .. (233.42,226.58) ;
\draw    (207,227.84) .. controls (246.25,248.59) and (252.5,217.34) .. (290,226.58) ;
\draw    (296.83,226.58) .. controls (338.38,231.92) and (350.38,203.34) .. (380,226.58) ;
\draw    (353.42,226.58) .. controls (375.88,237.09) and (391.13,235.59) .. (410.75,228.08) ;
\draw    (296.08,228.5) .. controls (320.63,236.59) and (330.75,235.6) .. (350.38,228.09) ;

\end{tikzpicture}

        \caption{A rational floor diagram of degree $4$. All weights of edges are $1$.}
        \label{fig-floordiagramnormal}
    \end{figure}
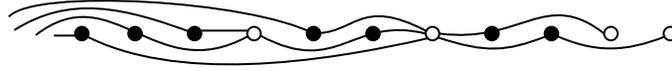
\end{example}

\begin{definition}\label{def-countoffloordwithoutmerging}
    Let $e$ be an edge of weight $w$ of a floor diagram.
    We define the \emph{quadratically enriched multiplicity of $e$} to be 
	\[\mult^{\mathbb{A}^1}(e)= w^{\AA^1} = 
	\begin{cases} 
		  \qinv{w}+\frac{w-1}{2}\cdot h & \mbox{if $w$ is odd,}\\
		\frac{w}{2}\cdot h & \mbox{if }w\mbox{ is even}.
	\end{cases}\]
    We define the \emph{quadratically enriched multiplicity} of a floor diagram $\cD$ to be
	$$\mult^{\mathbb{A}^1}(\cD)=\prod_e \mult^{\mathbb{A}^1}(e),$$
	where the product is taken over all bounded edges of $\cD$.
    We define the weighted count of rational floor diagrams of degree $\Delta$ to be
	\[
        N_\Delta^{\AA^1,\floor} 
        := \sum_\cD \mult^{\mathbb{A}^1}(\cD) 
    \]
    where the sum goes over all rational floor diagrams of degree $\Delta$.
\end{definition}

\begin{remark}\label{rem-arithmultfloor}
    The quadratically enriched multiplicity of a floor diagram interpolates between its complex and its real multiplicity in the sense that if we specialize to the corresponding field (via rank resp.\ signature), we get the complex resp.\ real multiplicity of a floor diagram, see e.g.\ \cite{BM08}. Furthermore, the quadratically enriched multiplicity is in fact determined by the real and the complex multiplicity:
    \[ \mult^{\AA^1}(\cD) = \begin{cases}
        \frac{\mult_\CC(\cD)}{2} \h & \text{if $\mult_\CC(\cD)$ is even,} \\
        \frac{\mult_\CC(\cD) - 1}{2} \h + \gw{\mult_\RR(\cD)} & \text{else.}
    \end{cases} \]
\end{remark}

\begin{theorem}[\cite{JPMPR23}, Theorem~C]\label{thm-tropcurvesfloorthinpoints}
The count of rational floor diagrams of degree $\Delta$ equals the count of rational tropical curves of degree $\Delta$ satisfying $k$-rational point conditions only, counted with their quadratically enriched multiplicity:
    \[N_\Delta^{\AA^1,\floor} =N_\Delta^{\AA^1,\trop} .\]
    In particular, there is a bijection from the set of rational floor diagrams of degree $\Delta$ to the set of rational tropical stable maps of degree $\Delta$ satisfying $n(\Delta)$ simple point conditions.
\end{theorem}
When specializing the result above to the complex or the real numbers, we have the well-known correspondence between the complex resp.\ real count of floor diagrams and the complex resp.\ real count of tropical curves \cite{BM08}. Because of Remark \ref{rem-arithmultfloor}, the statement above is actually a straight forward consequence of the complex and real correspondence.
We shortly explain the basic idea that leads to the correspondence: consider the set of rational tropical stable maps of degree $\Delta$ satisfying $n(\Delta)$ vertically stretched point conditions. One can show that any such curve is \emph{floor decomposed}, i.e.\ if we remove all vertical edges (elevator edges), the remaining connected components consist either of a single marked point or of a \emph{floor}, i.e.\ a path of edges leading from one non-vertical end to another. Each floor passes through one of the point conditions, and each elevator edge is adjacent to a marked point. We color the points on elevators black and the points on floors white, and then we shrink each floor to a vertex. The result is a floor diagram. The multiplicity of a floor diagram equals the multiplicity of the tropical stable map from which it is produced. Vice versa, for each floor diagram $\cD$ we can find precisely one tropical stable map satisfying our conditions that yields $\cD$ with the above procedure of coloring points and shrinking floors.
 
\begin{remark}\label{rem-labelsforfloordiagrams}
    For practical reasons, if we want to come up with a list of all rational floor diagrams of degree $\Delta$, one can start by drawing all possibilities how the white vertices can be connected, and then consider possibilities to distribute black points on the elevators (following the linear order). The latter are often called the \emph{markings} of a diagram on the white vertices. The multiplicity of a floor diagram does not depend on the linear order of the black vertices, and so one can sum the multiplicity times the number of markings to produce a total count of floor diagrams. 
\end{remark}

Now we introduce double points by merging adjacent simple points.
This is similar to the extension of floor diagrams considered in Section 2.2 \cite{BJP22}, and to Section 3.2 in \cite{BM08}.

The statement of the following lemma follows from the definition of floor diagram:
\begin{lemma}
    Consider all rational floor diagrams of degree $\Delta$. If we merge two neighboring points, the following possibilities arise:
    \begin{itemize}
        \item two floors (i.e.\ white vertices) can merge,
        \item a white vertex and an adjacent black vertex can merge,
        \item a white vertex and a non-adjacent black vertex can merge,
        \item two black vertices can merge.
    \end{itemize}
\end{lemma}

 Assume two black vertices adjacent to the same white vertex are merged, and assume that each black vertex is adjacent to another white vertex and that this pair of white vertices is also merged. Continuing like this, we assume that there is an involution on the underlying tree which exchanges two identical parts passing through the same pairs of points. Note that this is a phenomenon which cannot occur in a floor diagram without merged points: an automorphism of a floor diagram is necessarily trivial since it has to respect the linear order of the vertices. However, once we merge points, then we only have a partial order to be respected. In this case, non-trivial automorphisms can show up.

\begin{definition}
   Assume a floor diagram with merged points has an involution which exchanges two identical parts passing through the same pairs of points, i.e.\ an automorphism respecting the partial order we obtain from merging the points.
    Then these two identical parts are called a \emph{twin tree} of the floor diagram.
\end{definition}

Note that the smallest possible twin tree is given by two merged black points adjacent to ends of the same direction, and adjacent to the same white vertex.

\begin{lemma}\label{lem-twintreeforfloordiagrams}
    Given a floor diagram and its corresponding rational tropical stable map. Merging pairs of points to double points in the tropical stable map results in a piece which is a twin tree if and only if the corresponding piece of the floor diagram is a twin tree.
\end{lemma}
\begin{proof}
    By Definition \ref{def-twintree}, a twin tree of a tropical stable map ends at a vertex of type \fourValentVertex{}. By Proposition \ref{prop-vertextypesinfloordecomposedcurve}, its adjacent double edges are elevator edges. In the corresponding floor diagram, we must thus have two black vertices adjacent to the same white vertex which are merged. The involution on the tropical stable map yields an involution on the floor diagram when shrinking floors to white vertices. Vice versa, a twin tree in a floor diagram comes from a twin tree of the corresponding tropical stable map: two such merged black vertices come from a double point on a double elevator. If there are further adjacent merged white vertices, these come from adjacent floors, which, because of the involution on the floor diagram, must be identical: the weights of the adjacent incoming and outgoing elevators determine the directions of each edge of the floor, and they are the same for each of the two identical parts which are switched by the involution. 
\end{proof}

\begin{example}
Figure \ref{fig-floordiagrammerged} shows the rational floor diagram of Figure \ref{fig-floordiagramnormal} but with the points $4$ and $5$, $6$ and $7$, $8$ and $9$ and $10$ and $11$  merged. The edges connecting the merged point $6$ and $7$ with the merged point $8$ and $9$ and then $10$ and $11$ form a twin tree of this floor diagram.
 
    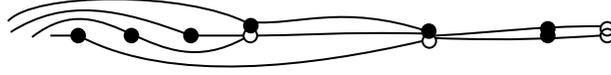
\begin{figure}[t]
        \centering
       
        \tikzset{every picture/.style={line width=0.75pt}} 
        
        \begin{tikzpicture}[x=0.75pt,y=0.75pt,yscale=-1,xscale=1]
            
            \draw   (200,226.58) .. controls (200,224.7) and (201.53,223.17) .. (203.42,223.17) .. controls (205.3,223.17) and (206.83,224.7) .. (206.83,226.58) .. controls (206.83,228.47) and (205.3,230) .. (203.42,230) .. controls (201.53,230) and (200,228.47) .. (200,226.58) -- cycle ;
            \draw  [fill={rgb, 255:red, 0; green, 0; blue, 0 }  ,fill opacity=1 ] (140,226.58) .. controls (140,224.7) and (141.53,223.17) .. (143.42,223.17) .. controls (145.3,223.17) and (146.83,224.7) .. (146.83,226.58) .. controls (146.83,228.47) and (145.3,230) .. (143.42,230) .. controls (141.53,230) and (140,228.47) .. (140,226.58) -- cycle ;
            \draw  [fill={rgb, 255:red, 0; green, 0; blue, 0 }  ,fill opacity=1 ] (113.17,226.58) .. controls (113.17,224.7) and (114.7,223.17) .. (116.58,223.17) .. controls (118.47,223.17) and (120,224.7) .. (120,226.58) .. controls (120,228.47) and (118.47,230) .. (116.58,230) .. controls (114.7,230) and (113.17,228.47) .. (113.17,226.58) -- cycle ;
            \draw  [fill={rgb, 255:red, 0; green, 0; blue, 0 }  ,fill opacity=1 ] (170,226.58) .. controls (170,224.7) and (171.53,223.17) .. (173.42,223.17) .. controls (175.3,223.17) and (176.83,224.7) .. (176.83,226.58) .. controls (176.83,228.47) and (175.3,230) .. (173.42,230) .. controls (171.53,230) and (170,228.47) .. (170,226.58) -- cycle ;
            \draw  [fill={rgb, 255:red, 0; green, 0; blue, 0 }  ,fill opacity=1 ] (200.25,221.58) .. controls (200.25,219.7) and (201.78,218.17) .. (203.67,218.17) .. controls (205.55,218.17) and (207.08,219.7) .. (207.08,221.58) .. controls (207.08,223.47) and (205.55,225) .. (203.67,225) .. controls (201.78,225) and (200.25,223.47) .. (200.25,221.58) -- cycle ;
            \draw  [fill={rgb, 255:red, 0; green, 0; blue, 0 }  ,fill opacity=1 ] (290,224.33) .. controls (290,222.45) and (291.53,220.92) .. (293.42,220.92) .. controls (295.3,220.92) and (296.83,222.45) .. (296.83,224.33) .. controls (296.83,226.22) and (295.3,227.75) .. (293.42,227.75) .. controls (291.53,227.75) and (290,226.22) .. (290,224.33) -- cycle ;
            \draw   (290.25,229.33) .. controls (290.25,227.45) and (291.78,225.92) .. (293.67,225.92) .. controls (295.55,225.92) and (297.08,227.45) .. (297.08,229.33) .. controls (297.08,231.22) and (295.55,232.75) .. (293.67,232.75) .. controls (291.78,232.75) and (290.25,231.22) .. (290.25,229.33) -- cycle ;
            \draw   (380,226.58) .. controls (380,224.7) and (381.53,223.17) .. (383.42,223.17) .. controls (385.3,223.17) and (386.83,224.7) .. (386.83,226.58) .. controls (386.83,228.47) and (385.3,230) .. (383.42,230) .. controls (381.53,230) and (380,228.47) .. (380,226.58) -- cycle ;
            \draw  [fill={rgb, 255:red, 0; green, 0; blue, 0 }  ,fill opacity=1 ] (350,223.17) .. controls (350,221.28) and (351.53,219.75) .. (353.42,219.75) .. controls (355.3,219.75) and (356.83,221.28) .. (356.83,223.17) .. controls (356.83,225.05) and (355.3,226.58) .. (353.42,226.58) .. controls (351.53,226.58) and (350,225.05) .. (350,223.17) -- cycle ;
            \draw  [fill={rgb, 255:red, 0; green, 0; blue, 0 }  ,fill opacity=1 ] (350,226.58) .. controls (350,224.7) and (351.53,223.17) .. (353.42,223.17) .. controls (355.3,223.17) and (356.83,224.7) .. (356.83,226.58) .. controls (356.83,228.47) and (355.3,230) .. (353.42,230) .. controls (351.53,230) and (350,228.47) .. (350,226.58) -- cycle ;
            \draw   (380,223.17) .. controls (380,221.28) and (381.53,219.75) .. (383.42,219.75) .. controls (385.3,219.75) and (386.83,221.28) .. (386.83,223.17) .. controls (386.83,225.05) and (385.3,226.58) .. (383.42,226.58) .. controls (381.53,226.58) and (380,225.05) .. (380,223.17) -- cycle ;
            \draw    (176.83,226.58) -- (200,226.58) ;
            \draw    (102.58,226.59) -- (116.58,226.58) ;
            \draw    (93.33,227.34) .. controls (132.5,197.09) and (153.25,253.84) .. (200.25,227.97) ;
            \draw    (82.5,224.84) .. controls (110.17,207.34) and (137,212.84) .. (173.42,226.58) ;
            \draw    (117.75,227.09) .. controls (163.5,251.09) and (233.75,242.33) .. (290.25,229.33) ;
            \draw    (204.88,219.85) .. controls (227.38,221.35) and (249.88,210.6) .. (291.63,224.09) ;
            \draw    (81.13,219.34) .. controls (94.13,204.84) and (170.67,203.83) .. (202.67,220.58) ;
            \draw    (206.83,226.58) .. controls (224,226.59) and (250.38,224.6) .. (292.38,225.59) ;
            \draw    (296.83,226.58) .. controls (330.25,224.84) and (346,221.84) .. (380.25,221.84) ;
            \draw    (296.75,227.84) .. controls (335.5,229.34) and (360,227.84) .. (380,226.58) ;
        
        \end{tikzpicture}

        \caption{The points $4$ and $5$, $6$ and $7$, $8$ and $9$ and $10$ and $11$ have been merged in the rational floor diagram of degree $4$ in Figure \ref{fig-floordiagramnormal}. The edges connecting the merged point $6$ and $7$ with the merged point $8$ and $9$ and then $10$ and $11$ form a twin tree of this floor diagram.}
        \label{fig-floordiagrammerged}
    \end{figure}

For more examples, see Table \ref{tab:Felixfantasticpicture}.  
\end{example}

\begin{definition}
    Given a rational floor diagram $\cD$ of degree $\Delta$ with $s$ pairs of merged points, we assign the following quadratically enriched multiplicity:

    First, note that the partial order of the vertices induces a linear order on the pairs of merged points, which we can therefore denote by $q_1,\ldots,q_s$. 

    Let $\mathcal{T}$ be a twin tree of $\cD$. Assume that the merged points $q_1,\ldots,q_t$ are on $\mathcal{T}$. 
     If edges of higher weight appear in a twin tree, they must be bounded. 
     We define the quadratically enriched multiplicity of the twin tree to be

    \begin{equation*}
            \mult^{\AA^1}(\mathcal{T})=\prod_{(e, e') \text{ twin edge in $\cT$}} \mult^{\AA^1}(e,e') \cdot 
            \qinv{2^{t-1}}\cdot
            \sum_{\substack{I\subset \{1,\dots,t\}\\ \vert I\vert\equiv m_{\circ}\mod 2}} \qinv{\prod_{i\in I} d_i}  
            \in \GW(k).
        \end{equation*}

    Let  $V_{\fatPointOnVertex{}}$ be the subset of merged points where a white vertex and a neighboring black vertex are merged. Let $R\subset \{1,\ldots,s\}$ be the subset of merged points which are neither in $V_{\fatPointOnVertex{}}$ nor contained in a twin tree.   
    We define the quadratically enriched multiplicity of $\cD$ as
    
    $$\mult^{\AA^1}(\cD)=\prod_{\mathcal{T}}
    \mult^{\AA^1}(\mathcal{T})\cdot \prod_{ v \in V_{\fatPointOnVertex{}}} \gamma(m_v, d_v) \cdot
    \prod_{i\in R} \beta_i\cdot
    \prod_e (m_e)^{\AA^1},
    $$  
    where the first product goes over all twin trees $\mathcal{T}$ of $\cD$,
    and the last product goes over all 
    bounded edges of $\cD$ which are not contained in a twin tree and $m_e$ denotes their weight. As always, we let $m_v$ be the weight of the elevator incident to the type \fatPointOnVertex{} vertex $v$ and $d_v$ defining the point condition on $v$.

   We define  $N^{\A^1,\floor}_\Delta(r,(d_1,\ldots,d_s))$ to be the weighted count of rational floor diagrams of degree $\Delta$, with $s$ merged points, counted with the quadratically enriched multiplicity
   \[N^{\A^1,\floor}_\Delta(r,(d_1,\ldots,d_s))\coloneqq \sum_{\mathcal{D}} \mult^{\AA^1}(\mathcal{D}).\]
\end{definition}

\begin{example}
    Consider the rational floor diagram $\cD$ with of degree $4$ with $4$ pairs of merged points in Figure \ref{fig-floordiagrammerged}.
    We assume that the four merged points correspond to the extensions $k\subset k(\sqrt{d_1}), \ldots, k\subset k(\sqrt{d_4})$, ordered from left to right (i.e.\ in accordance with the linear order of the vertices of the floor diagram).
    Then $\cD$ has quadratically enriched multiplicity 
    $$\mult^{\AA^1}(\cD) = \big(\langle 2d_3\rangle+\langle 2d_{4}\rangle \big) \cdot \gamma(1, d_2) \cdot \beta_1 =  \big(\langle 2d_3\rangle+\langle 2d_{4}\rangle \big) \cdot \beta_1.$$
    As the adjacent edge from which the black vertex was moved is of weight $1$, the $\gamma$-expression is equal to $\langle 1\rangle$ and can be dropped in the product.
    
    For more examples, see Table \ref{tab:Felixfantasticpicture}.
    
\end{example}

The following theorem states that enriched counts of tropical curves in the presence of double point conditions can be computed with a floor diagram count. This is \ref{thm-flooreqtrop_intro} from the introduction.

\begin{theorem}\label{thm-flooreqtrop}
    The count of rational floor diagrams of degree $\Delta$ with $s$ merged points with quadratically enriched multiplicity equals the count of rational tropical stable maps of degree $\Delta$ satisfying $r$ simple and $s$ double point conditions with quadratically enriched multiplicity, where $r+2s=n(\Delta)$:
    $$N^{\A^1,\floor}_\Delta(r,(d_1,\ldots,d_s)) =N^{\A^1,\trop}_\Delta(r,(d_1,\ldots,d_s)).$$
\end{theorem}

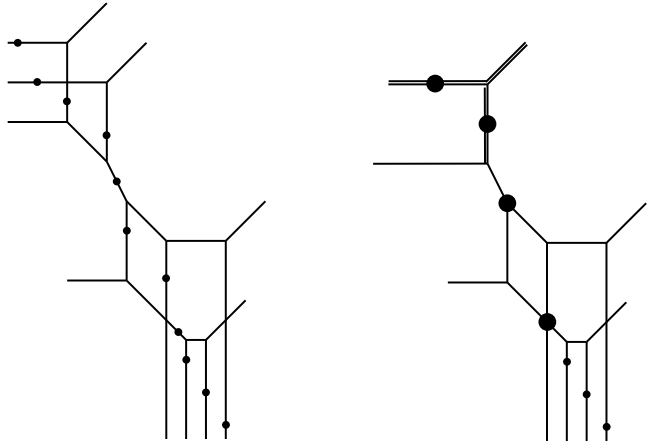
\begin{figure}[t]
        \centering

\tikzset{every picture/.style={line width=0.75pt}} 

\begin{tikzpicture}[x=0.75pt,y=0.75pt,yscale=-1,xscale=1]

\draw    (120,250) -- (150,250) ;
\draw    (150,250) -- (150,290) ;
\draw    (170,230) -- (150,250) ;
\draw    (190,250) -- (170,270) ;
\draw    (120,270) -- (170,270) ;
\draw    (170,270) -- (170,310) ;
\draw    (120,290) -- (150,290) ;
\draw    (150,290) -- (170,310) ;
\draw    (170,310) -- (180,330) ;
\draw    (180,330) -- (180,370) ;
\draw    (180,330) -- (200,350) ;
\draw    (200,350) -- (200,450) ;
\draw    (150,370) -- (180,370) ;
\draw    (200,350) -- (230,350) ;
\draw    (180,370) -- (210,400) ;
\draw    (210,400) -- (220,400) ;
\draw    (210,400) -- (210,450) ;
\draw    (220,400) -- (220,450) ;
\draw    (240,380) -- (220,400) ;
\draw    (250,330) -- (230,350) ;
\draw    (230,350) -- (230,450) ;
\draw  [fill={rgb, 255:red, 0; green, 0; blue, 0 }  ,fill opacity=1 ] (123.6,250.03) .. controls (123.6,249.21) and (124.26,248.55) .. (125.08,248.55) .. controls (125.9,248.55) and (126.57,249.21) .. (126.57,250.03) .. controls (126.57,250.85) and (125.9,251.51) .. (125.08,251.51) .. controls (124.26,251.51) and (123.6,250.85) .. (123.6,250.03) -- cycle ;
\draw  [fill={rgb, 255:red, 0; green, 0; blue, 0 }  ,fill opacity=1 ] (228.61,442.86) .. controls (228.61,442.05) and (229.27,441.38) .. (230.09,441.38) .. controls (230.91,441.38) and (231.57,442.05) .. (231.57,442.86) .. controls (231.57,443.68) and (230.91,444.35) .. (230.09,444.35) .. controls (229.27,444.35) and (228.61,443.68) .. (228.61,442.86) -- cycle ;
\draw  [fill={rgb, 255:red, 0; green, 0; blue, 0 }  ,fill opacity=1 ] (218.52,426.48) .. controls (218.52,425.66) and (219.18,425) .. (220,425) .. controls (220.82,425) and (221.48,425.66) .. (221.48,426.48) .. controls (221.48,427.3) and (220.82,427.97) .. (220,427.97) .. controls (219.18,427.97) and (218.52,427.3) .. (218.52,426.48) -- cycle ;
\draw  [fill={rgb, 255:red, 0; green, 0; blue, 0 }  ,fill opacity=1 ] (208.61,410.01) .. controls (208.61,409.19) and (209.27,408.52) .. (210.09,408.52) .. controls (210.91,408.52) and (211.57,409.19) .. (211.57,410.01) .. controls (211.57,410.83) and (210.91,411.49) .. (210.09,411.49) .. controls (209.27,411.49) and (208.61,410.83) .. (208.61,410.01) -- cycle ;
\draw  [fill={rgb, 255:red, 0; green, 0; blue, 0 }  ,fill opacity=1 ] (204.61,396.01) .. controls (204.61,395.19) and (205.27,394.52) .. (206.09,394.52) .. controls (206.91,394.52) and (207.57,395.19) .. (207.57,396.01) .. controls (207.57,396.83) and (206.91,397.49) .. (206.09,397.49) .. controls (205.27,397.49) and (204.61,396.83) .. (204.61,396.01) -- cycle ;
\draw  [fill={rgb, 255:red, 0; green, 0; blue, 0 }  ,fill opacity=1 ] (198.4,368.86) .. controls (198.4,368.05) and (199.06,367.38) .. (199.88,367.38) .. controls (200.7,367.38) and (201.37,368.05) .. (201.37,368.86) .. controls (201.37,369.68) and (200.7,370.35) .. (199.88,370.35) .. controls (199.06,370.35) and (198.4,369.68) .. (198.4,368.86) -- cycle ;
\draw  [fill={rgb, 255:red, 0; green, 0; blue, 0 }  ,fill opacity=1 ] (178.61,344.86) .. controls (178.61,344.05) and (179.27,343.38) .. (180.09,343.38) .. controls (180.91,343.38) and (181.57,344.05) .. (181.57,344.86) .. controls (181.57,345.68) and (180.91,346.35) .. (180.09,346.35) .. controls (179.27,346.35) and (178.61,345.68) .. (178.61,344.86) -- cycle ;
\draw  [fill={rgb, 255:red, 0; green, 0; blue, 0 }  ,fill opacity=1 ] (173.52,320) .. controls (173.52,319.18) and (174.18,318.52) .. (175,318.52) .. controls (175.82,318.52) and (176.48,319.18) .. (176.48,320) .. controls (176.48,320.82) and (175.82,321.48) .. (175,321.48) .. controls (174.18,321.48) and (173.52,320.82) .. (173.52,320) -- cycle ;
\draw  [fill={rgb, 255:red, 0; green, 0; blue, 0 }  ,fill opacity=1 ] (168.38,296.7) .. controls (168.38,295.88) and (169.04,295.21) .. (169.86,295.21) .. controls (170.68,295.21) and (171.34,295.88) .. (171.34,296.7) .. controls (171.34,297.51) and (170.68,298.18) .. (169.86,298.18) .. controls (169.04,298.18) and (168.38,297.51) .. (168.38,296.7) -- cycle ;
\draw  [fill={rgb, 255:red, 0; green, 0; blue, 0 }  ,fill opacity=1 ] (133.42,269.82) .. controls (133.42,269) and (134.08,268.33) .. (134.9,268.33) .. controls (135.72,268.33) and (136.38,269) .. (136.38,269.82) .. controls (136.38,270.64) and (135.72,271.3) .. (134.9,271.3) .. controls (134.08,271.3) and (133.42,270.64) .. (133.42,269.82) -- cycle ;
\draw    (304.31,311.08) -- (362,311) ;
\draw    (360.62,272.62) -- (360.75,311) ;
\draw    (381.21,249.77) -- (361.21,269.77) ;
\draw    (382,251) -- (362,271) ;
\draw    (312,271) -- (362,271) ;
\draw    (362,271) -- (362,311) ;
\draw    (312.12,269.41) -- (361.75,269.35) ;
\draw    (362,311) -- (372,331) ;
\draw    (372,331) -- (372,371) ;
\draw    (372,331) -- (392,351) ;
\draw    (392,351) -- (392,451) ;
\draw    (342,371) -- (372,371) ;
\draw    (392,351) -- (422,351) ;
\draw    (372,371) -- (402,401) ;
\draw    (402,401) -- (412,401) ;
\draw    (402,401) -- (402,451) ;
\draw    (412,401) -- (412,451) ;
\draw    (432,381) -- (412,401) ;
\draw    (442,331) -- (422,351) ;
\draw    (422,351) -- (422,451) ;
\draw  [fill={rgb, 255:red, 0; green, 0; blue, 0 }  ,fill opacity=1 ] (420.61,443.86) .. controls (420.61,443.05) and (421.27,442.38) .. (422.09,442.38) .. controls (422.91,442.38) and (423.57,443.05) .. (423.57,443.86) .. controls (423.57,444.68) and (422.91,445.35) .. (422.09,445.35) .. controls (421.27,445.35) and (420.61,444.68) .. (420.61,443.86) -- cycle ;
\draw  [fill={rgb, 255:red, 0; green, 0; blue, 0 }  ,fill opacity=1 ] (410.52,427.48) .. controls (410.52,426.66) and (411.18,426) .. (412,426) .. controls (412.82,426) and (413.48,426.66) .. (413.48,427.48) .. controls (413.48,428.3) and (412.82,428.97) .. (412,428.97) .. controls (411.18,428.97) and (410.52,428.3) .. (410.52,427.48) -- cycle ;
\draw  [fill={rgb, 255:red, 0; green, 0; blue, 0 }  ,fill opacity=1 ] (400.61,411.01) .. controls (400.61,410.19) and (401.27,409.52) .. (402.09,409.52) .. controls (402.91,409.52) and (403.57,410.19) .. (403.57,411.01) .. controls (403.57,411.83) and (402.91,412.49) .. (402.09,412.49) .. controls (401.27,412.49) and (400.61,411.83) .. (400.61,411.01) -- cycle ;
\draw  [fill={rgb, 255:red, 0; green, 0; blue, 0 }  ,fill opacity=1 ] (331.63,270.6) .. controls (331.63,268.41) and (333.41,266.63) .. (335.6,266.63) .. controls (337.8,266.63) and (339.57,268.41) .. (339.57,270.6) .. controls (339.57,272.8) and (337.8,274.58) .. (335.6,274.58) .. controls (333.41,274.58) and (331.63,272.8) .. (331.63,270.6) -- cycle ;
\draw  [fill={rgb, 255:red, 0; green, 0; blue, 0 }  ,fill opacity=1 ] (358.03,291) .. controls (358.03,288.81) and (359.81,287.03) .. (362,287.03) .. controls (364.19,287.03) and (365.97,288.81) .. (365.97,291) .. controls (365.97,293.19) and (364.19,294.97) .. (362,294.97) .. controls (359.81,294.97) and (358.03,293.19) .. (358.03,291) -- cycle ;
\draw  [fill={rgb, 255:red, 0; green, 0; blue, 0 }  ,fill opacity=1 ] (368.03,331) .. controls (368.03,328.81) and (369.81,327.03) .. (372,327.03) .. controls (374.19,327.03) and (375.97,328.81) .. (375.97,331) .. controls (375.97,333.19) and (374.19,334.97) .. (372,334.97) .. controls (369.81,334.97) and (368.03,333.19) .. (368.03,331) -- cycle ;
\draw  [fill={rgb, 255:red, 0; green, 0; blue, 0 }  ,fill opacity=1 ] (388.2,390.92) .. controls (388.2,388.73) and (389.98,386.95) .. (392.17,386.95) .. controls (394.37,386.95) and (396.15,388.73) .. (396.15,390.92) .. controls (396.15,393.12) and (394.37,394.9) .. (392.17,394.9) .. controls (389.98,394.9) and (388.2,393.12) .. (388.2,390.92) -- cycle ;
\draw  [fill={rgb, 255:red, 0; green, 0; blue, 0 }  ,fill opacity=1 ] (148.38,279.58) .. controls (148.38,278.77) and (149.04,278.1) .. (149.86,278.1) .. controls (150.68,278.1) and (151.34,278.77) .. (151.34,279.58) .. controls (151.34,280.4) and (150.68,281.07) .. (149.86,281.07) .. controls (149.04,281.07) and (148.38,280.4) .. (148.38,279.58) -- cycle ;

\end{tikzpicture}
        \caption{A rational tropical stable map of degree $4$ satisfying $11$ simple point conditions, which (after a homeomorphism) can be viewed as vertically stretched, and a rational tropical stable map of degree $4$ satisfying $4$ double and $3$ simple point conditions. The picture on the left yields the floor diagram from Figure \ref{fig-floordiagramnormal}, the one on the right from Figure \ref{fig-floordiagrammerged}. Just as the latter floor diagram is obtained from the previous by merging points, the picture on the right is obtained from the one on the left by merging points.}
        \label{fig:curvesforfloordiagrams}
\end{figure}

\begin{proof}
    By Theorem \ref{thm-tropcurvesfloorthinpoints}, there is a bijection of sets from the set of floor diagrams of degree $\Delta$ to the set of rational tropical stable maps of degree $\Delta$ satisfying $n(\Delta)$ vertically stretched simple point conditions. Any rational tropical stable map of degree $\Delta$ satisfying $r$ simple and $s$ double vertically stretched point conditions is obtained from a rational tropical stable maps of degree $\Delta$ satisfying $n(\Delta)$ vertically stretched simple point conditions by merging points. In the same way, a floor diagram with merged points is obtained from a floor diagram. Vice versa, by \enquote{unmerging} points we can produce a rational tropical stable map of degree $\Delta$ satisfying $n(\Delta)$ vertically stretched simple point conditions from a rational tropical stable map of degree $\Delta$ satisfying $r$ simple and $s$ double vertically stretched point conditions, resp.\ a floor diagram from a floor diagram with merged points. For an example, see Figure \ref{fig:curvesforfloordiagrams}. 
   From a rational tropical stable map of degree $\Delta$ satisfying $r$ simple and $s$ double vertically stretched point conditions, we can form a floor diagram with merged points by shrinking floors as usual, and marking double points as merged points with colors according to what part of the map they fix, floors or elevator edges. Of course, we have to take the parametrization by the abstract tropical curve into account (e.g.\ in the situation of twin trees).
    It follows that there is a bijection between the set of rational tropical stable maps of degree $\Delta$ satisfying $r$ simple and $s$ double vertically stretched point conditions and the set of floor diagrams with merged points. By definition and by Lemma \ref{lem-twintreeforfloordiagrams}, each floor diagram with merged points counts towards $N^{\A^1,\floor}_\Delta(r,(d_1,\ldots,d_s))$ with the same quadratically enriched multiplicity as its corresponding rational tropical stable map counts towards $N^{\A^1,\trop}_\Delta(r,(d_1,\ldots,d_s))$. The statement follows. 
\end{proof}

\section{Computational results}\label{sec-computations}

In this section we present quadratically enriched counts for low degrees. These numbers were obtained by manual computation using our floor diagram approach as introduced in \cref{sec-floor}. To write down our results, we use the following notation. Let $\Delta$ be a degree taken from \cref{fig-degrees} and $n  (\Delta)= r + 2s$. Each of the $s$ non-$k$-rational points is assumed to live in a quadratic extension~$k(\sqrt{d_i})$ for some $d_i \in k$. Set $\beta_i = \gw{2} + \gw{2d_i}$ and write 
\[\beta_s^{(l)} = \sum_{1 \leq i_1 < i_2 < \cdots < i_l \leq s} \beta_{i_1}\beta_{i_2} \cdots \beta_{i_l} \]
for the $l$-th elementary symmetric polynomial evaluated at the $\beta_i$. We write all counts in this section in the form 
\[c_{s+1} \h + c_1 \beta_s^{(1)} + c_2 \beta_s^{(2)} + \cdots + c_s \beta_s^{(s)} + c_0 \gw{1} \]
with $c_0,\ldots,c_{s+1}\in \mathbb{Z}$. Recall that one can always do this by \cite{BrugalleRauWickelgren} (see subsection \ref{sec-generalformula}).

The base case of all point conditions being $k$-rational, i.e. $s = 0$ and $r = n$, was determined in \cite{JPP23} and we repeat the results here for completeness. Also note that the count for $s'$ can be computed from the count for $s$ conjugate point conditions with $s > s'$ by \cref{prop-specializetokrational}.

For all these examples in the last row, we also provide a formula for $N^{\mathbb{A}^1}_{\Delta}(\sigma)$ for general $\sigma$, extending beyond the multiquadratic cases.
Note that these in fact hold since $N^{\mathbb{A}^1}_{\Delta}(\sigma)$ is a Witt invariant in the sense of Serre by \cite{BrugalleRauWickelgren} (see the introduction for an explanation of this phenomenon).

The general formula which holds for all $\sigma$ will be expressed in terms of $\mathbf{a}_n$, where $$\mathbf{a}_n = a_n(\Tr_{k(\sigma)/k}(\gw{1})),$$ with $a_n$ denoting the power structure on the Grothendieck-Witt ring from \cite{PajwaniPal}.

\subsection{Curves in $\PP^2$} 
The quadratically enriched count of rational plane cubics was computed in \cite{KassLevineSolomonWickelgren}. It equals $2\h+\Tr_{k(\sigma)/k}(\gw{1})$ for a generic configuration of points defined over arbitrary finite field extensions. For point conditions over quadratic field extensions, this quadratically enriched count coincides with the corresponding quadratically enriched tropical counts. We repeat them here for the sake of completeness and provide quadratically enriched counts for rational plane quartics in \cref{tab-computationsP2}.
The numbers in this table can be cross-checked and coincide with the ones obtained from \cite{JP24}. 

\begin{table}[t]
    \begin{tabular}{c|c|c|l}
        $d$ & $n(\Delta)$ & $(r, s)$ & $N_\Delta^{\AA^1, \trop}(r, d_1, \ldots, d_s)$ \\
        \hline
        \multirow{5}{*}{3} & \multirow{5}{*}{8} & $(8, 0)$ & $2\h + 8 \gw{1}$ \\
        & & $(6, 1)$ & $2\h + \beta_1^{(1)} + 6\gw{1}$ \\
        & & $(4, 2)$ & $2\h + \beta_2^{(1)} + 4\gw{1}$ \\
        & & $(2, 3)$ & $2\h + \beta_3^{(1)} + 2\gw{1}$ \\
        & & $(0, 4)$ & $2\h + \beta_4^{(1)}$ \\
        & & $\sigma$ & $2\h+\mathbf{a}_1$\\
        \hline
        \multirow{6}{*}{4} & \multirow{6}{*}{11} & $(11, 0)$ & $190 h + 240\gw{1}$ \\
        & & $(9, 1)$ & $190 h + 48\beta_1^{(1)} + 144\gw{1}$ \\
        & & $(7, 2)$ & $190 h + 8\beta_2^{(2)} + 32 \beta_2^{(1)} +80\gw{1}$ \\
        & & $(5, 3)$ & $190 h + \beta_3^{(3)} + 6 \beta_3^{(2)} + 20 \beta_3^{(1)} + 40\gw{1}$ \\
        & & $(3, 4)$ & $190 h + \beta_4^{(3)} + 4 \beta_4^{(2)} + 12 \beta_4^{(1)} + 16\gw{1}$ \\
        & & $(1, 5)$ & $190 h + \beta_5^{(3)} + 2 \beta_5^{(2)} + 8 \beta_5^{(1)}$\\
        & & $\sigma$ & $190 h + \mathbf{a}_3-\mathbf{a}_2+2\mathbf{a}_1-2\mathbf{a}_0$
    \end{tabular}
    \caption{Counts of curves in $\PP^2$.}
    \label{tab-computationsP2}
\end{table}

\begin{figure}[hb]
    \def\scalefloordiagram{0.8}
\newcommand{\fourlegsdown}{\filldraw (0,0) ellipse (0.8 and 0.2);
    \path[draw] (-0.6, 0) -- (-0.6, -0.5);
    \path[draw] (-0.2, 0) -- (-0.2, -0.5);
    \path[draw] (0.2, 0) -- (0.2, -0.5);
    \path[draw] (0.6, 0) -- (0.6, -0.5);}
\newcommand{\threelegsdown}{\filldraw (0,0) ellipse (0.6 and 0.2);
    \path[draw] (-0.4, 0) -- (-0.4, -0.5);
    \path[draw] (0, 0) -- (0, -0.5);
    \path[draw] (0.4, 0) -- (0.4, -0.5);}
\newcommand{\twolegsdown}{\filldraw (0,0) ellipse (0.4 and 0.2);
    \path[draw] (-0.2, 0) -- (-0.2, -0.5);
    \path[draw] (0.2, 0) -- (0.2, -0.5);}

\begin{tikzpicture}[scale = \scalefloordiagram]
    \fourlegsdown
    
    \path[draw] (0,0) -- (0, 1);
    \draw (0.1, 0.5) node[anchor = west] {$3$};
    
    \filldraw (0, 1) ellipse (0.2 and 0.2);
    \path[draw] (0,1) -- (0, 2);
    \draw (0.1, 1.5) node[anchor = west] {$2$};
    
    \filldraw (0, 2) ellipse (0.2 and 0.2);
    \path[draw] (0,2) -- (0, 3);
    
    \filldraw (0, 3) ellipse (0.2 and 0.2);
\end{tikzpicture}
\quad
\begin{tikzpicture}[scale = \scalefloordiagram]
    \fourlegsdown
    
    \path[draw] (0,0) -- (0, 1);
    \draw (0.1, 0.5) node[anchor = west] {$3$};
    
    \filldraw (0, 1) ellipse (0.5 and 0.2);
    \path[draw] (0.3,1) -- (0.3, 2);
    \path[draw] (-0.3,1) -- (-0.3, 3);
    
    \filldraw (0.3, 2) ellipse (0.2 and 0.2);
    
    \filldraw (-0.3, 3) ellipse (0.2 and 0.2);
\end{tikzpicture}
\quad
\begin{tikzpicture}[scale = \scalefloordiagram]
    \fourlegsdown
    
    \path[draw] (-0.3,0) -- (-0.3, 1);
    \draw (-0.3, 0.5) node[anchor = east] {$2$};
    \path[draw] (0.3,0) -- (0.3, 3);
    
    \filldraw (-0.3, 1) ellipse (0.2 and 0.2);
    \path[draw] (-0.3,1) -- (-0.3, 2);
    
    \filldraw (0.3, 3) ellipse (0.2 and 0.2);
    
    \filldraw (-0.3, 2) ellipse (0.2 and 0.2);
\end{tikzpicture}
\quad
\begin{tikzpicture}[scale = \scalefloordiagram]
    \fourlegsdown
    
    \path[draw] (-0.3,0) -- (-0.3, 1);
    \draw (-0.3, 0.5) node[anchor = east] {$2$};
    \path[draw] (0.3,0) -- (0.3, 2);
    
    \filldraw (-0.3, 1) ellipse (0.2 and 0.2);
    \path[draw] (-0.3,1) -- (-0.3, 3);
    
    \filldraw (0.3, 2) ellipse (0.2 and 0.2);
    
    \filldraw (-0.3, 3) ellipse (0.2 and 0.2);
\end{tikzpicture}
\quad
\begin{tikzpicture}[scale = \scalefloordiagram]
    \fourlegsdown
    
    \path[draw] (-0.3,0) -- (-0.3, 2);
    \draw (-0.3, 1) node[anchor = east] {$2$};
    \path[draw] (0.3,0) -- (0.3, 1);
    
    \filldraw (-0.3, 2) ellipse (0.2 and 0.2);
    \path[draw] (-0.3,2) -- (-0.3, 3);
    
    \filldraw (0.3, 1) ellipse (0.2 and 0.2);
    
    \filldraw (-0.3, 3) ellipse (0.2 and 0.2);
\end{tikzpicture}
\quad
\begin{tikzpicture}[scale = \scalefloordiagram]
    \fourlegsdown
    
    \path[draw] (-0.4,0) -- (-0.4, 3);
    \path[draw] (0,0) -- (0, 2);
    \path[draw] (0.4,0) -- (0.4, 1);
    
    \filldraw (-0.4, 3) ellipse (0.2 and 0.2);	
    \filldraw (0, 2) ellipse (0.2 and 0.2);	
    \filldraw (0.4, 1) ellipse (0.2 and 0.2);
\end{tikzpicture}
\quad
\begin{tikzpicture}[scale = \scalefloordiagram]
    \threelegsdown
    
    \path[draw] (0.4,0) -- (0.4, 1);
    \draw (0.3, 0.5) node[anchor = east] {$2$};
    
    \filldraw (0.6, 1) ellipse (0.4 and 0.2);
    \path[draw] (0.6,1) -- (0.6, 2);
    \path[draw] (0.8,1) -- (0.8, -0.5);
    \draw (0.5, 1.5) node[anchor = east] {$2$};
    
    \filldraw (0.6, 2) ellipse (0.2 and 0.2);
    \path[draw] (0.6,2) -- (0.6, 3);
    
    \filldraw (0.6, 3) ellipse (0.2 and 0.2);
\end{tikzpicture}
\quad
\begin{tikzpicture}[scale = \scalefloordiagram]
    \threelegsdown
    
    \path[draw] (0.4,0) -- (0.4, 1);
    \draw (0.3, 0.5) node[anchor = east] {$2$};
    
    \filldraw (0.4, 1) ellipse (0.2 and 0.2);
    \path[draw] (0.4,1) -- (0.4, 2);
    
    \filldraw (0.6, 2) ellipse (0.4 and 0.2);
    \path[draw] (0.8,2) -- (0.8, -0.5);
    \path[draw] (0.6,2) -- (0.6, 3);
    
    \filldraw (0.6, 3) ellipse (0.2 and 0.2);
\end{tikzpicture}
\\[1em]
\begin{tikzpicture}[scale = \scalefloordiagram]
    \twolegsdown
    
    \path[draw] (0.2,0) -- (0.2, 1);
    
    \filldraw (0.5, 1) ellipse (0.5 and 0.2);
    \path[draw] (0.5,1) -- (0.5, 2);
    \path[draw] (0.6,1) -- (0.6, -0.5);
    \path[draw] (0.8,1) -- (0.8, -0.5);
    \draw (0.4, 1.5) node[anchor = east] {$2$};
    
    \filldraw (0.5, 2) ellipse (0.2 and 0.2);
    \path[draw] (0.5,2) -- (0.5, 3);
    
    \filldraw (0.5, 3) ellipse (0.2 and 0.2);
\end{tikzpicture}
\quad
\begin{tikzpicture}[scale = \scalefloordiagram]
    \twolegsdown
    
    \path[draw] (0.2,0) -- (0.2, 1);
    
    \filldraw (0.4, 1) ellipse (0.4 and 0.2);
    \path[draw] (0.5,1) -- (0.5, 2);
    \path[draw] (0.6,1) -- (0.6, -0.5);
    
    \filldraw (0.7, 2) ellipse (0.4 and 0.2);
    \path[draw] (0.7,2) -- (0.7, 3);
    \path[draw] (0.9,2) -- (0.9, -0.5);
    
    \filldraw (0.7, 3) ellipse (0.2 and 0.2);
\end{tikzpicture}
\quad
\begin{tikzpicture}[scale = \scalefloordiagram]
    \threelegsdown
    
    \path[draw] (0.4,0) -- (0.4, 1);
    \draw (0.3, 0.5) node[anchor = east] {$2$};
    
    \filldraw (0.6, 1) ellipse (0.4 and 0.2);
    \path[draw] (0.4,1) -- (0.4, 3);
    \path[draw] (0.8,1) -- (0.8, -0.5);
    \path[draw] (0.8,1) -- (0.8, 2);
    
    \filldraw (0.8, 2) ellipse (0.2 and 0.2);
        
    \filldraw (0.4, 3) ellipse (0.2 and 0.2);
\end{tikzpicture}
\quad
\begin{tikzpicture}[scale = \scalefloordiagram]
    \twolegsdown
    
    \path[draw] (0.2,0) -- (0.2, 1);
    
    \filldraw (0.5, 1) ellipse (0.5 and 0.2);
    \path[draw] (0.3,1) -- (0.3, 3);
    \path[draw] (0.6,1) -- (0.6, -0.5);
    \path[draw] (0.8,1) -- (0.8, -0.5);
    \path[draw] (0.7,1) -- (0.7, 2);
    
    \filldraw (0.7, 2) ellipse (0.2 and 0.2);
    
    \filldraw (0.3, 3) ellipse (0.2 and 0.2);
\end{tikzpicture}
\quad
\begin{tikzpicture}[scale = \scalefloordiagram]
    \threelegsdown
    \path[draw] (-0.3, 0) -- (-0.3, 2);
    \path[draw] (0.2, 0) -- (0.2, 1);
    
    \filldraw (-0.5, 2) ellipse (0.4 and 0.2);
    \path[draw] (-0.7, 2) -- (-0.7, -0.5);
    \path[draw] (-0.5, 2) -- (-0.5, 3);
    
    \filldraw (-0.5, 3) ellipse (0.2 and 0.2);
    \filldraw (0.2, 1) ellipse (0.2 and 0.2);
\end{tikzpicture}
\quad 
\begin{tikzpicture}[scale = \scalefloordiagram]
    \threelegsdown
    \path[draw] (-0.3, 0) -- (-0.3, 1);
    \path[draw] (0.2, 0) -- (0.2, 2);
    
    \filldraw (-0.5, 1) ellipse (0.4 and 0.2);
    \path[draw] (-0.7, 1) -- (-0.7, -0.5);
    \path[draw] (-0.5, 1) -- (-0.5, 3);
    
    \filldraw (-0.5, 3) ellipse (0.2 and 0.2);
    \filldraw (0.2, 2) ellipse (0.2 and 0.2);
\end{tikzpicture}
\quad
\begin{tikzpicture}[scale = \scalefloordiagram]
    \threelegsdown
    \path[draw] (-0.3, 0) -- (-0.3, 1);
    \path[draw] (0.2, 0) -- (0.2, 3);
    
    \filldraw (-0.5, 1) ellipse (0.4 and 0.2);
    \path[draw] (-0.7, 1) -- (-0.7, -0.5);
    \path[draw] (-0.5, 1) -- (-0.5, 2);
    
    \filldraw (-0.5, 2) ellipse (0.2 and 0.2);
    \filldraw (0.2, 3) ellipse (0.2 and 0.2);
\end{tikzpicture}
\quad
\begin{tikzpicture}[scale = \scalefloordiagram]
    \filldraw (0, 0) ellipse (0.4 and 0.2);
    \path[draw] (0.2, 0) -- (0.2, -0.5);
    \path[draw] (-0.2, 0) -- (-0.2, -0.5);
    \path[draw] (0, 0) -- (0, 2);
    
    \filldraw (0.9, 1) ellipse (0.4 and 0.2);
    \path[draw] (0.7, 1) -- (0.7, -0.5);
    \path[draw] (1.1, 1) -- (1.1, -0.5);
    \path[draw] (0.9, 1) -- (0.9, 2);
    
    \filldraw (0.45, 2) ellipse (0.6 and 0.2);
    \path[draw] (0.45, 2) -- (0.45, 3);
    
    \filldraw (0.45, 3) ellipse (0.2 and 0.2);
    
\end{tikzpicture}
    \caption{The floor diagrams we used to count plane quartics (showing only the floor vertices, the markings for the elevator vertices have to be added, see \cref{rem-labelsforfloordiagrams}). }
\end{figure}
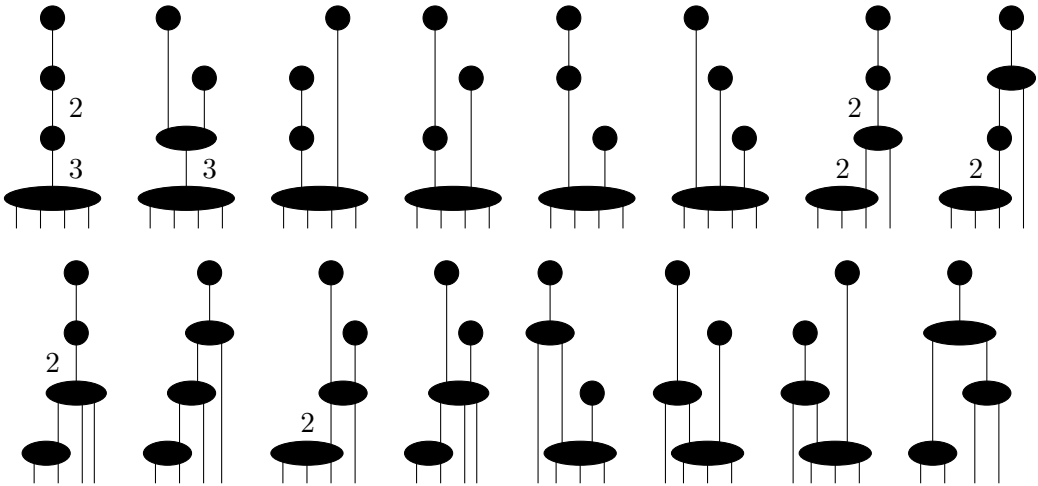

\subsection{Curves in other del Pezzo surfaces}

In $\PP^1 \times \PP^1$ with  bidegree $(2,3)$, over the complex numbers there are 96 curves through 9 points. 
In bidegree $(2, 4)$, over the complex numbers there are 640 curves through 11 points. 
In bidegree $(2, 5)$, over the complex numbers there are 3840 curves through 13 points. 
The quadratically enriched counts for these bidegrees are listed in \cref{tab-computationsP1xP1}.

\begin{table}[t]
    \begin{tabular}{c|c|c|l}
        $(a_1, a_2)$ & $n(\Delta)$ & $(r, s)$ & $N_\Delta^{\AA^1, \trop}(r, d_1, \ldots, d_s)$ \\
        \hline
        \multirow{4}{*}{$(2, 2)$} & \multirow{4}{*}{7} & $(7, 0)$ & $2 h + 8 \gw{1}$ \\
        & & $(5, 1)$ & $2 h + \beta_1^{(1)} + 6 \gw{1}$ \\
        & & $(3, 2)$ & $2 h + \beta_2^{(1)} + 4\gw{1}$ \\
        & & $(1, 3)$ & $2h + \beta_3^{(1)} + 2\gw{1}$ \\
        & & $\sigma$ & $2h+\mathbf{a}_1+\mathbf{a}_0$\\
        \hline
        \multirow{5}{*}{$(2, 3)$} & \multirow{5}{*}{9} & $(9, 0)$ &  $24 h + 48 \gw{1}$ \\
        & & $(7, 1)$ & $24 h + 8 \beta_1^{(1)} + 32 \gw{1}$ \\
        & & $(5, 2)$ & $24 h +  \beta_2^{(2)} + 6 \beta_2^{(1)} + 20 \gw{1}$ \\
        & & $(3, 3)$ & $24 h + \beta_3^{(2)} + 4 \beta_3^{(1)} + 12\gw{1}$ \\
        & & $(1, 4)$ & $24 h + \beta_4^{(2)} + 2 \beta_4^{(1)} +8 \gw{1}$ \\
        & & $\sigma$ & $24 h + \mathbf{a}_2+3\mathbf{a}_0$\\
        \hline
        \multirow{6}{*}{$(2, 4)$} & \multirow{6}{*}{11} & $(11, 0)$ & $192 h + 256 \gw{1}$ \\
        & & $(9, 1)$ & $192 h  + 48 \beta_1^{(1)} + 160 \gw{1}$ \\
        & & $(7, 2)$ & $192 h + 8 \beta_2^{(2)} + 32 \beta_2^{(1)} + 96 \gw{1}$ \\
        & & $(5, 3)$ & $192 h + \beta_3^{(3)} + 6 \beta_3^{(2)} + 20 \beta_3^{(1)} + 56 \gw{1}$ \\
        & & $(3, 4)$ & $192 h + \beta_4^{(3)} + 4\beta_4^{(2)} + 12\beta_4^{(1)} + 32 \gw{1}$ \\
        & & $(1, 5)$ & $192 h + \beta_5^{(3)} + 2\beta_5^{(2)} + 8 \beta_5^{(1)} + 16 \gw{1}$ \\
        & & $\sigma$ & $192 h + \mathbf{a}_3-\mathbf{a}_2+2\mathbf{a}_1+14\mathbf{a}_0$\\
        \hline
        \multirow{7}{*}{$(2, 5)$} & \multirow{7}{*}{13} & $(13, 0)$ & $1280 \h + 1280 \gw{1}$ \\
        & & $(11, 1)$ & $1280 \h + 256 \beta_1^{(1)} + 768 \gw{1}$ \\
        & & $(9, 2)$ & $1280 \h + 48 \beta_2^{(2)} + 160 \beta_2^{(1)} + 448 \gw{1}$ \\
        & & $(7, 3)$ & $1280 \h + 8 \beta_3^{(3)} + 32 \beta_3^{(2)} + 96 \beta_3^{(1)} + 256\gw{1}$ \\
        & & $(5, 4)$ & $1280 \h + \beta_4^{(4)} + 6 \beta_4^{(3)} + 20 \beta_4^{(2)} + 56 \beta_4^{(1)} + 144 \gw{1}$ \\
        & & $(3, 5)$ & $1280 \h + \beta_5^{(4)} + 4 \beta_5^{(3)} + 12 \beta_5^{(2)} + 32 \beta_5^{(1)} + 80\gw{1}$ \\
        & & $(1, 6)$ & $1280 \h + \beta_6^{(4)} + 2 \beta_6^{(3)} + 8 \beta_6^{(2)} + 16 \beta_6^{(1)} + 48\gw{1}$\\
        & & $\sigma$ & $1280 \h +\mathbf{a}_4-2\mathbf{a}_3+2\mathbf{a}_2+14\mathbf{a}_1+6\mathbf{a}_0$
    \end{tabular}
    \caption{Counts of curves in $\PP^1 \times \PP^1$.}
    \label{tab-computationsP1xP1}
\end{table}

These results are novel and allow us to remark that for bidegree $(2,4)$, the coefficients of the symmetric polynomials of the $\beta_i$'s coincide with the quartic plane curve counting invariants. This relation is a corollary of the wall-crossing formula for quadratic invariants in \cite{JP24} and Pick's theorem.

\begin{corollary}
    If $\Delta$ and $\Delta'$ are two smooth convex lattice polygons with the same number of interior lattice points, then for $2s\leq \min\{n(\Delta),n(\Delta')\}$ we have   
    \[
    N^{\A^1,\floor}_\Delta(r,(d_1,\ldots,d_s))-N^{\A^1,\trop}_{\Delta'}(r,(d_1,\ldots,d_s))=
    \big(W_{\Delta,s}-W_{\Delta',s} \big)\cdot \gw{1}\in\operatorname{W}(k),
    \]
    where $\operatorname{W}(k) = \GW(k) / (\h)$ is the Witt ring of $k$.
\end{corollary}

Therefore, such coefficients depend only on the number of interior lattice points of the polygon. 
Our computations also fit in this context.

Finally, for some computations in other del Pezzo surfaces see Tables~\ref{tab-computationsBl1}, \ref{tab-computationsBl2}, and~\ref{tab-computationsBl3}.

\begin{table}[t]
    \begin{tabular}{c|c|c|l}
        $(d, a_1)$ & $n(\Delta)$ & $(r, s)$ & $N_\Delta^{\AA^1, \trop}(r, d_1, \ldots, d_s)$ \\
        \hline
        \multirow{4}{*}{$(3, 1)$} & \multirow{4}{*}{7} & $(7, 0)$ & $2 h + 8 \gw{1}$ \\
        & & $(5, 1)$ & $2 h + \beta_1^{(1)} + 6 \gw{1}$ \\
        & & $(3, 2)$ & $2 h + \beta_2^{(1)} + 4\gw{1}$ \\
        & & $(1, 3)$ & $2h + \beta_3^{(1)} + 2\gw{1}$ \\
        & & $\sigma$ & $2h + \mathbf{a}_1+\mathbf{a}_0$\\
        \hline
        \multirow{5}{*}{$(4, 2)$} & \multirow{5}{*}{9} & $(9, 0)$ & $24 h + 48 \gw{1}$ \\
        & & $(7, 1)$ & $24 h + 8 \beta_1^{(1)} + 32 \gw{1}$  \\
        & & $(5, 2)$ & $24 h +  \beta_2^{(2)} + 6 \beta_2^{(1)} + 20 \gw{1}$ \\
        & & $(3, 3)$ & $24h + \beta_3^{(2)} + 4 \beta_3^{(1)} + 12\gw{1}$ \\
        & & $(1, 4)$ & $24h + \beta_4^{(2)} + 2 \beta_4^{(1)} +8 \gw{1}$ \\
        & & $\sigma$ & $24h + \mathbf{a}_2+3\mathbf{a}_0$\\
    \end{tabular}
    \caption{Counts of curves in the blow up of $\PP^2$ in one point. Here, $\Delta$ is the third polygon in \cref{fig-degrees}.}
    \label{tab-computationsBl1}
\end{table}

\begin{table}[b]
    \begin{tabular}{c|c|c|l}
        $(d, a_1, a_2)$ & $n(\Delta)$ & $(r, s)$ & $N_\Delta^{\AA^1, \trop}(r, d_1, \ldots, d_s)$ \\
        \hline
        \multirow{4}{*}{$(4, 2, 2)$} & \multirow{4}{*}{7} & $(7, 0)$ & $2 h + 8 \gw{1}$ \\
        & & $(5, 1)$ & $2 h + \beta_1^{(1)} + 6 \gw{1}$ \\
        & & $(3, 2)$ & $2 h + \beta_2^{(1)} + 4\gw{1}$ \\
        & & $(1, 3)$ & $2h + \beta_3^{(1)} + 2\gw{1}$ \\
        & & $\sigma$ & $2h+ \mathbf{a}_1+\mathbf{a}_0$\\
        \hline
        \multirow{5}{*}{$(4, 2, 1)$} & \multirow{5}{*}{8} & $(8, 0)$ & $24 h + 48 \gw{1}$ \\
        & & $(6, 1)$ & $24 h + 8 \beta_1^{(1)} + 32 \gw{1}$ \\
        & & $(4, 2)$ & $24 h +  \beta_2^{(2)} + 6 \beta_2^{(1)} + 20 \gw{1}$ \\
        & & $(2, 3)$ & $24h +  \beta_3^{(2)} + 4 \beta_3^{(1)} + 12\gw{1}$ \\
        & & $(0, 4)$ & $24h +  \beta_4^{(2)} + 2\beta_4^{(1)} +8 \gw{1}$ \\
        & & $\sigma$ & $24h +\mathbf{a}_2+\mathbf{a}_1+4\mathbf{a}_0$\\
        \hline
        \multirow{5}{*}{$(4,1,1)$} & \multirow{5}{*}{9} & $(9, 0)$ & $160 \h + 240 \gw{1}$ \\
        & & $(7, 1)$ & $160 \h + 48 \beta_1^{(1)} + 144 \gw{1}$ \\
        & & $(5, 2)$ & $160 \h + 8 \beta_2^{(2)} + 32 \beta_2^{(1)} + 80 \gw{1}$ \\
        & & $(3, 3)$ & $160 \h + \beta_3^{(3)} + 6 \beta_3^{(2)} + 20 \beta_3^{(1)} + 40 \gw{1}$ \\
        & & $(1, 4)$ & $160 \h + \beta_4^{(3)} + 4 \beta_4^{(2)} + 12 \beta_4^{(1)} + 16 \gw{1}$\\
        & & $\sigma$ & $160 \h + \mathbf{a}_3+\mathbf{a}_2+3\mathbf{a}_1+3\mathbf{a}_0$\\
    \end{tabular}
    \caption{Counts of curves in the blow up of $\PP^2$ in two points. Here, $\Delta$ is the fourth polygon in \cref{fig-degrees}.}
    \label{tab-computationsBl2}
\end{table}

\begin{table}[t]
    \begin{tabular}{c|c|c|l}
        $(d, a_1, a_2, a_3)$ & $n(\Delta)$ & $(r, s)$ & $N_\Delta^{\AA^1, \trop}(r, d_1, \ldots, d_s)$ \\
        \hline
        \multirow{3}{*}{$(3, 1, 1, 1)$} & \multirow{3}{*}{5} & $(5, 0)$ & $2 h + 8 \gw{1}$ \\
        & & $(3, 1)$ & $2 h + \beta_1^{(1)} + 6 \gw{1}$ \\
        & & $(1, 2)$ & $2 h + \beta_2^{(1)} + 4\gw{1}$ \\
        & & $\sigma$ & $2 h +\mathbf{a}_1+3\mathbf{a}_0$\\
        \hline
        \multirow{4}{*}{$(4, 1, 1, 2)$} & \multirow{4}{*}{7} & $(7, 0)$ & $24 h + 48 \gw{1}$ \\
        & & $(5, 1)$ & $24 h + 8 \beta_1^{(1)} + 32 \gw{1}$ \\
        & & $(3, 2)$ & $24 h +  \beta_2^{(2)} + 6 \beta_2^{(1)} + 20 \gw{1}$ \\
        & & $(1, 3)$ & $24h + \beta_3^{(2)} + 4 \beta_3^{(1)} + 12\gw{1}$\\
        & & $\sigma$ & $24h+ \mathbf{a}_2+2\mathbf{a}_1+6\mathbf{a}_0$\\
    \end{tabular}
    \caption{Counts of curves in the blow up of $\PP^2$ in three points. Here, $\Delta$ is the rightmost polygon in \cref{fig-degrees}.}
    \label{tab-computationsBl3}
\end{table}

\bibliographystyle{alpha} 
\bibliography{bibliographie} 
\end{document}